\documentclass[a4paper,11pt]{amsart}

\usepackage{amssymb}
\usepackage{cite}
\usepackage{url}

\newtheorem{theorem}{Theorem}[section]
\newtheorem{lemma}[theorem]{Lemma}
\newtheorem{corollary}[theorem]{Corollary}
\newtheorem{proposition}[theorem]{Proposition}

\theoremstyle{definition}
\newtheorem{remark}[theorem]{Remark}

\numberwithin{equation}{section}

\makeatletter
\@namedef{subjclassname@2020}{%
  \textup{2020} Mathematics Subject Classification}
\makeatother

\DeclareMathOperator{\RE}{Re}
\DeclareMathOperator{\IM}{Im}
\DeclareMathOperator{\Gal}{Gal}
\DeclareMathOperator{\supp}{supp}

\DeclareMathOperator{\diag}{diag}
\DeclareMathOperator{\tr}{tr}
\DeclareMathOperator{\Cond}{Cond}
\DeclareMathOperator{\Res}{Res}
\DeclareMathOperator{\Var}{Var}
\DeclareMathOperator{\poly}{poly}

\begin{document}

\title[The value-distribution of Artin $L$-functions]
{The value-distribution of Artin $L$-functions associated with cubic fields in conductor aspect}

\author[M. Mine]{Masahiro Mine}
\address{Faculty of Science and Technology\\ Sophia University\\ 7-1 Kioi-cho, Chiyoda-ku, Tokyo 102-8554, Japan}
\email{m-mine@sophia.ac.jp}

\date{}

\begin{abstract}
Arising from the factorizations of Dedekind zeta-functions of cubic fields, we obtain Artin $L$-functions of certain two-dimensional representations. 
In this paper, we study the value-distribution of such Artin $L$-functions for families of non-Galois cubic fields in conductor aspect. 
We prove that various mean values of the Artin $L$-functions are represented by integrals involving a density function which can be explicitly constructed. 
By the class number formula, the result is applied to the study on the distribution of class numbers of cubic fields. 
\end{abstract}

%% 11M06	$\zeta (s)$ and $L(s, \chi)$
%% 11M32	Multiple Dirichlet series and zeta functions and multizeta values
%% 11M35	Hurwitz and Lerch zeta functions
%% 11M36	Selberg zeta functions and regularized determinants
%% 11M41	Other Dirichlet series and zeta functions
%% 11R42	Zeta functions and $L$-functions of number fields
%% 11F66	Langlands $L$-functions


\thanks{The work of this paper was supported by Grant-in-Aid for JSPS Fellows (Grant Number JP21J00529).}

\maketitle

\section{Introduction}\label{sec:1}
Let $d$ be a non-square integer such that $d \equiv 0,1 \bmod4$, namely, a discriminant of a binary quadratic form. 
The quadratic Dirichlet $L$-function $L(s,\chi_d)$ is the Dirichlet $L$-function of the character $\chi_d(n)=(\frac{d}{n})$, where $(\frac{\cdot}{\cdot})$ indicates the Kronecker symbol. 
The distribution of  values $L(s,\chi_d)$ has been studied by various authors. 
One of the earliest results was obtained by Chowla and Erd\H{o}s \cite{ChowlaErdos1951}. 
They proved the existence of a continuous function $F_\sigma$ for $\sigma>3/4$ such that the limit formula
\begin{gather}\label{eq:07250107}
\lim_{X \to\infty} \frac{\# \left\{ 0<d<X ~\middle|~ L(\sigma,\chi_d) \leq e^a \right\}}{X/2}
=F_\sigma(a)
\end{gather}
holds for any $a \in \mathbb{R}$. 
Furthermore, $F_\sigma$ is the distribution function of a probability measure on $\mathbb{R}$, that is, $F_\sigma$ is non-decreasing over $\mathbb{R}$ and satisfies $\lim_{t \to -\infty} f(t) =0$, $\lim_{t \to \infty} f(t) =1$. 
Elliott \cite{Elliott1970, Elliott1980} studied quadratic Dirichlet $L$-functions for negative discriminants, and obtained a similar result at $\sigma=1$: 
\begin{gather}\label{eq:07250108}
\lim_{X \to\infty} \frac{\#\left\{0<-d<X ~\middle|~ {L}(1,\chi_d)\leq{e}^a \right\}}{X/2}
=F_1(a).
\end{gather}
It was also shown that $F_1$ possesses a probability density function $Q$ whose Fourier transform is represented as the infinite product 
\begin{gather*}
\widetilde{Q}(\xi)
=\prod_{p} \left[ \frac{1}{2} \left(1-\frac{1}{p}\right) \left(1-\frac{1}{p}\right)^{-i \xi} 
+\frac{1}{2} \left(1-\frac{1}{p}\right) \left(1+\frac{1}{p}\right)^{-i \xi} 
+\frac{1}{p} \right], 
\end{gather*}
where $p$ runs through all prime numbers. 
Some related results were also obtained in \cite{Elliott1971, Elliott1972, Elliott1973}. 
These studies led to the idea of comparing the value-distributions of $L(s,\chi_d)$ with a suitable random model, which brought the recent progress in the theory of value-distributions of zeta and $L$-functions; see \cite{GranvilleSoundararajan2006, GranvilleSoundararajan2003, Lamzouri2010, Lamzouri2011b, LamzouriLesterRadziwill2019} for example. 

The values of quadratic Dirichlet $L$-functions at $s=1$ are connected with the class numbers. 
Let $h_d$ denote the class number of a discriminant $d$ in the narrow sense. 
Put $\epsilon_d=(u_d+v_d \sqrt{d})/2$ for $d>0$, where $(u_d,v_d)$ is the fundamental solution of the Pell equation $u^2-dv^2=4$. 
Then we obtain $L(1,\chi_d) \sqrt{d} = h_d \log{\epsilon_d}$ for $d>0$ and $L(1,\chi_d) \sqrt{|d|}=\pi h_d$ for $d<-4$ by Dirichlet's class number formula. 
Hence, limit formulas \eqref{eq:07250107} and \eqref{eq:07250108} yield that 
\begin{gather*}
\# \left\{ 0<d<X ~\middle|~ h_d \log{\epsilon_d} \leq e^a \sqrt{d} \right\}
\sim \frac{X}{2}F_1(a), \\
\# \left\{ 0<-d<X ~\middle|~ h_d \leq \frac{e^a}{\pi} \sqrt{|d|} \right\}
\sim \frac{X}{2} F_1(a) 
\end{gather*}
as $X \to\infty$ for $a \in\mathbb{R}$. 
Here, $f(x) \sim g(x)$ stands for $f(x)/g(x) \to 1$ as $x$ tends to some limit. 
Recall that there is more classical work on the asymptotic behaviors of $h_d \log{\epsilon_d}$ and $h_d$. 
Gauss stated without proof that
\begin{gather*}
\sum_{\substack{0<d<X \\ d \equiv0 \bmod4}} h_d \log{\epsilon_d}
\sim \frac{\pi^2}{42 \zeta(3)} X^{3/2}
\quad\text{and}\quad
\sum_{\substack{0<-d<X \\ d \equiv0 \bmod4}} h_d
\sim \frac{\pi}{42\zeta(3)} X^{3/2}
\end{gather*}
as $X \to\infty$, where $\zeta(s)$ is the usual Riemann zeta-function. 
Moreover, we have more precise formulas
\begin{align}
\sum_{0<d<X} h_d \log{\epsilon_d}
&=\frac{\pi^2}{18 \zeta(3)} X^{3/2} 
+O(X \log{X}), \label{eq:07250113}\\
\sum_{0<-d<X} h_d
&=\frac{\pi}{18 \zeta(3)} X^{3/2} 
+O(X \log{X}) \label{eq:07250114}
\end{align} 
due to Siegel \cite{Siegel1944}. 
Recall further that Ayoub \cite{Ayoub1963, Ayoub1964} studied the case where $d$ varies over the set of discriminants of quadratic fields, and obtained asymptotic formulas similar to \eqref{eq:07250113} and \eqref{eq:07250114} up to the coefficients of main terms.

\subsection{Artin $L$-functions associated with cubic fields}\label{sec:1.1}
Let $K$ be a cubic field with discriminant $d_K$ which is non-Galois over $\mathbb{Q}$. 
We denote by $[K]$ the isomorphism class of cubic fields containing $K$. 
Then we define 
\begin{gather*}
L_3^\pm(X)
=\left\{[K] ~\middle|~ \text{$K$ is a non-Galois cubic field with $0<\pm d_K<X$} \right\} 
\end{gather*}
and put $N_3^\pm(X)=\# L_3^\pm(X)$. 
Throughout this paper, we write $K \in L_3^\pm(X)$ instead of $[K] \in L_3^\pm(X)$ for simplicity. 
For every $K \in K_3^\pm(X)$, the Dedekind zeta function $\zeta_K(s)$ is factorized as 
\begin{gather*}
\zeta_K(s)
=L(s,\rho_K) \zeta(s), 
\end{gather*}
where $\rho_K$ is the standard representation of the Galois group $\Gal(\widehat{K}/\mathbb{Q}) \simeq S_3$, and $L(s,\rho_K)$ is the attached Artin $L$-function. 
Here, $\widehat{K}$ indicates the Galois closure of $K$ over $\mathbb{Q}$, and $S_3$ is the symmetric group of degree $3$. 
Let $F=\mathbb{Q}(\sqrt{d_K})$. 
Then the Galois group $\Gal(\widehat{K}/F)$ is isomorphic to the cyclic group of order $3$.
We see that $\rho_K$ is the representation induced from the non-trivial character of $\Gal(\widehat{K}/F)$. 
The Artin $L$-function $L(s,\rho_K)$ is therefore holomorphic over the whole complex plane. 
Moreover, the strong Artin conjecture is true for $L(s,\rho_K)$, i.e.\ there exists a cuspidal representation $\pi$ of $GL_2(\mathbb{A}_\mathbb{Q})$ such that $L(s,\rho_K)=L(s,\pi)$ holds. 
By definition, the value $L(\sigma,\rho_K)$ is real whenever $\sigma$ is a real number. 

The purpose of this paper is to study the value-distribution of $L(s,\rho_K)$ as the cubic field $K$ varies in  the family $L_3^\pm(X)$. 
The detailed statements of the results are presented in Section \ref{sec:2}. 
In this section, we pick up two of them for comparison with the above results on quadratic Dirichlet $L$-functions. 

\begin{theorem}\label{thm:1.1}
Let $\sigma>7/8$ be a real number. 
Then there exists a non-negative $C^\infty$-function $C_\sigma$ on $\mathbb{R}$ such that 
\begin{gather}\label{eq:07250112}
\lim_{X \to\infty} \frac{\# \left\{ K \in L_3^\pm(X) ~\middle|~ L(\sigma,\rho_K) \leq e^a \right\}}{N_3^\pm(X)}
=\int_{-\infty}^{a} C_\sigma(x) \,\frac{dx}{\sqrt{2\pi}}
\end{gather}
holds for any $a \in\mathbb{R}$. 
Furthermore, the Fourier transform of $C_\sigma$ is represented as  
\begin{align*}
\widetilde{C}_\sigma(\xi)
&=\prod_{p} \frac{1}{1+p^{-1}+p^{-2}} 
\Biggl[ \frac{1}{6} \left(1-\frac{1}{p^{\sigma}}\right)^{-2i \xi} 
+\frac{1}{2} \left(1-\frac{1}{p^{2\sigma}}\right)^{-i \xi} \\
&\qquad\qquad
+\frac{1}{3} \left(1+\frac{1}{p^{\sigma}}+\frac{1}{p^{2\sigma}}\right)^{-i \xi} 
+\frac{1}{p} \left(1-\frac{1}{p^{\sigma}}\right)^{-i \xi}
+\frac{1}{p^2} \Biggr], 
\end{align*} 
where $p$ runs through all prime numbers. 
\end{theorem}

Limit formula \eqref{eq:07250112} gives analogues of \eqref{eq:07250107} and \eqref{eq:07250108} for the Artin $L$-functions. 
In this paper, we further evaluate the rate of convergence in \eqref{eq:07250112}; see Theorem \ref{thm:2.4}. 
Notice that $C_\sigma$ is similar to Elliott's density function $Q$ in view of the infinite product representations of their Fourier transforms. 
See Theorem \ref{thm:2.1} for more information about the density function $C_\sigma$. 

Put $D^+=4$ and $D^-=2\pi$. 
By the class number formula, we have 
\begin{gather}\label{eq:07250133}
L(1,\rho_K)
=D^\pm \frac{h_K R_K}{\sqrt{|d_K|}},
\end{gather}
where $h_K$ and $R_K$ denote the class number and the regulator of a cubic field $K$, respectively. 
As analogues of \eqref{eq:07250113} and \eqref{eq:07250114}, we prove the following asymptotic formulas. 

\begin{theorem}\label{thm:1.2}
There exists an absolute constant $\delta>0$ such that 
\begin{gather*}
\sum_{K \in L_3^+(X)} h_K R_K
=c X^{3/2}
+O\left(X^{3/2} \exp\left(-\delta \frac{\log{X}}{\log\log{X}}\right)\right),\\
\sum_{K \in L_3^-(X)} h_K R_K
=\frac{6}{\pi} c X^{3/2}
+O\left(X^{3/2} \exp\left(-\delta \frac{\log{X}}{\log\log{X}}\right)\right), 
\end{gather*}
where $c$ is a positive constant represented as 
\begin{gather*}
c
=\frac{\pi^2 \zeta(3)}{432} \prod_{p} (1+p^{-2}-2p^{-3}-2p^{-4}+2p^{-6}+p^{-7}-p^{-8}). 
\end{gather*}
\end{theorem}

By a standard argument using the partial summation, Theorem \ref{thm:1.2} is deduced from \eqref{eq:07250133} and some estimate on the first moment of $L(1,\rho_K)$. 
More generally, we show an asymptotic formula for the $z$-th moment
\begin{gather}\label{eq:07250226}
M_{z, \sigma}^{\pm}(X)
=\sum_{K \in L_3^\pm(X) \setminus E_\sigma(X)} L(\sigma,\rho_K)^z
\end{gather}
for $\sigma>7/8$ with $z \in \mathbb{C}$, where $E_\sigma(X)$ is a suitable subset of $L_3^\pm(X)$ satisfying at least $\# E_\sigma(X)=o(X)$ as $X \to\infty$. 
See Theorem \ref{thm:2.2} for the strict statement.

\subsection{Related topics}\label{sec:1.2}
The study of this paper is paper is partially motivated by the recent work on ``$M$-functions'' by Ihara--Matsumoto. 
We recall one of the results from \cite{IharaMatsumoto2011b} in the number field case. 
Let $F$ be $\mathbb{Q}$ or an imaginary quadratic field. 
For a prime ideal $\mathbf{f}$ of $F$, we define $X(\mathbf{f})$ as the set of all primitive Dirichlet characters $\chi$ on $F$ with conductor $\mathbf{f}$. 
Suppose that the Generalized Riemann Hypothesis (GRH) is true, that is, every Dirichlet $L$-function $L(s,\chi)$ has no zeros in the half-plane $\RE(s)>1/2$. 
Then $\log{L}(s,\chi)$ extends to a holomorphic function on $\RE(s)>1/2$. 

\begin{theorem}[Ihara--Matsumoto \cite{IharaMatsumoto2011b}]\label{thm:1.3}
Assume GRH, and denote by $|dz|$ the measure $(2\pi)^{-1} dx dy$ with $z=x+iy$. 
Then there exists a non-negative $C^\infty$-function $M_\sigma$ on $\mathbb{C}$ such that 
\begin{gather}\label{eq:07250142}
\lim_{\substack{\mathbf{f}~ {\rm prime} \\ N(\mathbf{f}) \to\infty}} 
\frac{1}{\# X(\mathbf{f})} \sum_{\chi \in X(\mathbf{f})} \Phi(\log{L}(s,\chi))
=\int_{\mathbb{C}} \Phi(z) M_\sigma(z) \,|dz|
\end{gather}
holds for any complex number $s=\sigma+i\tau$ with $\sigma>1/2$, where $\Phi$ is any continuous function on $\mathbb{C}$ satisfying $\Phi(z) \ll e^{a|z|}$ for some $a>0$. 
\end{theorem}

They also showed that a similar result is valid when $\log{L}(s,\chi)$ is replaced with the logarithmic derivative $(L'/L)(s,\chi)$. 
Furthermore, several analogous results were proved in \cite{Ihara2008, IharaMatsumoto2011a, IharaMatsumoto2014}, and so on. 
The construction of the function $M_\sigma$ was explained in \cite[Section 3]{IharaMatsumoto2011a}, and it matches a density function in the classical result of Bohr--Jessen \cite{BohrJessen1930, BohrJessen1932} on the value-distribution of the Riemann zeta-function. 
The density functions such as $M_\sigma$ were named \textit{$M$-functions} by Ihara \cite{Ihara2008}. 
Today there are a lot of variants of Theorem \ref{thm:1.3} and the corresponding $M$-functions; see the survey of Matsumoto \cite{Matsumoto2019}. 
In particular, the $M$-function for the value-distribution of quadratic Dirichlet $L$-functions was studied by Mourtada--Murty \cite{MourtadaMurty2015}. 
The density function $C_\sigma$ of Theorem \ref{thm:1.1} is regarded as a cubic analogue of Mourtada--Murty's $M$-function. 
We prove a limit formula similar to \eqref{eq:07250142}; see Theorem \ref{thm:2.5}. 

Another topic related to this paper is the work on the Artin $L$-function $L(s,\rho_K)$ due to Cho--Kim. 
They studied $L(s,\rho_K)$ not only for cubic fields but also for general $S_n$-fields of degree $n \geq2$. 
Here, a number field $K$ of degree $n$ is called an $S_n$-field if the Galois group $\Gal(\widehat{K}/\mathbb{Q})$ is isomorphic to the symmetric group $S_n$. 
Their results were often obtained under the following two conjectures: 
\begin{itemize}
\item
the strong Artin conjecture for $L(s,\rho_K)$,  
\item
the ``counting conjecture'' for $S_n$-fields; see \cite[Conjecture 3.1]{ChoKim2018}.
\end{itemize}
The truth of the former conjecture is known for $n \leq4$, and the latter conjecture is for $n \leq5$. 
Hence the results for $n=2,3,4$ are unconditional. 
In \cite{ChoKim2018}, they proved an asymptotic formula for integral moments of $\log{L}(1,\rho_K)$. 
Here we refer to the result in the cubic case. 

\begin{theorem}[Cho--Kim \cite{ChoKim2018}]\label{thm:1.4}
Let $k$ be a positive integer. 
Then we have 
\begin{gather*}
\frac{1}{N_3^\pm(X)} \sum_{K \in L_3^\pm(X)} \left(\log{L}(1,\rho_K)\right)^k
=\tilde{r}(k)
+O\left(\frac{1}{\log{X}}\right), 
\end{gather*}
where $\tilde{r}(k)$ is a positive constant which can be explicitly described. 
\end{theorem}

If we assume that $\tilde{r}(k) \ll c^{k \log\log{k}}$ is satisfied for some $c>1$, the method of moments enables us to show the existence of a continuous function $F$ satisfying 
\begin{gather*}
\lim_{X \to\infty} 
\frac{\# \left\{K \in L_3^\pm(X) ~\middle|~ L(1,\rho_K) \leq e^a\right\}}{N_3^\pm(X)}
=F(a),
\end{gather*}
where $a$ is a point of continuity of $F$. 
Then \eqref{eq:07250112} refines and generalizes this limit formula. 
Furthermore, it is remarkable that Theorem \ref{thm:1.1} can be proved without any assumptions on the constants $\tilde{r}(k)$. 
The original representation of $\tilde{r}(k)$ by Cho--Kim is described in \cite[Proposition 5.3]{ChoKim2018}, while we obtain another representation
\begin{gather*}
\tilde{r}(k)
=\int_{-\infty}^{\infty} x^k C_1(x) \,\frac{dx}{\sqrt{2\pi}}
\end{gather*}
by using the density function $C_\sigma$ of Theorem \ref{thm:1.1}. 
See also Corollary \ref{cor:6.1}. 

In addition, we recall another result due to Cho--Kim on the distribution of values $L(1,\rho_K)$. 
Let $\alpha>0$ be a real number. 
Then it was proved in \cite{ChoKim2017b} that
\begin{gather}\label{eq:07250155}
\liminf_{X \to\infty} 
\frac{\# \left\{K \in L_3^\pm(X) ~\middle|~ |L(1,\rho_K)-\alpha|<\epsilon\right\}}{N_3^\pm(X)}
>0
\end{gather}
for every $\epsilon>0$. 
This is a cubic analogue of the denseness result obtained by Mishou--Nagoshi \cite{MishouNagoshi2006b}. 
Let $\sigma>7/8$. 
Note that limit formula \eqref{eq:07250112} yields
\begin{gather}\label{eq:07250154}
\lim_{X \to\infty} \frac{\# \left\{ K \in L_3^\pm(X) ~\middle|~ |L(\sigma,\rho_K)-\alpha|<\epsilon \right\}}{N_3^\pm(X)}
=\int_{\log(\alpha-\epsilon)}^{\log(\alpha+\epsilon)} C_\sigma(x) \,\frac{dx}{\sqrt{2\pi}}
\end{gather}
if $\epsilon>0$ is sufficiently small. 
Furthermore, the support of the density function $C_\sigma$ equals to $\mathbb{R}$ for $7/8<\sigma \leq1$; see Theorem \ref{thm:2.1}.  
Thus the right-hand side of \eqref{eq:07250154} is positive. 
As a result, we recover \eqref{eq:07250155} and extend it for $7/8<\sigma \leq1$. 

\vspace{\baselineskip}

The organization of this paper is as follows. 
\begin{itemize}
\item
Section \ref{sec:2} is devoted to presenting the statement of the main results of this paper. 
\item
In Section \ref{sec:3}, we collect preliminary lemmas used later. 
\item
The proofs of the main results begin with the study of the density function $C_\sigma$ described in Theorem \ref{thm:1.1}. 
In Section \ref{sec:4.1}, we study the random Euler products attached to the Artin $L$-function $L(s,\rho_K)$. 
Then we explain the construction of the density function $C_\sigma$ in Section \ref{sec:4.2}. 
Several analytic properties of $C_\sigma$ are also proved in Section \ref{sec:4.3}. 
\item
Next, we associate the mean value of $\Phi(\log{L}(\sigma,\rho_K))$ with the integral involving the density function $C_\sigma$, where $\Phi$ is a test function as in Theorem \ref{thm:1.3}. 
We show in Section \ref{sec:5} an asymptotic formula for the $z$-th moment described in \eqref{eq:07250226}, namely, the mean value of $\Phi(\log{L}(\sigma,\rho_K))$ in the case $\Phi(x)=e^{zx}$.
\item
The proofs of the results are completed in Section \ref{sec:6}. 
We prove Theorem \ref{thm:1.1} by the asymptotic formula of the $z$-th moment with $z \in i \mathbb{R}$. 
As described before, we also deduce Theorem \ref{thm:1.2} from the case $z=1$ by using the partial summation. 
Finally, we prove an analogue of Theorem \ref{thm:1.3} for a general test function $\Phi$. 
\item
As well as the work of Ihara--Matsumoto, one can obtain similar results for logarithmic derivatives $(L'/L)(s,\rho_K)$, which are presented in the appendix. 
\end{itemize}

\section{Statement of results}\label{sec:2}
To begin with, we set up the notation for counting cubic fields with discriminants not exceeding a given quantity. 
Based on \cite{Roberts2001, TaniguchiThorne2013}, we introduce the notion of the local specifications of cubic fields as follows. 
Let 
\begin{gather*}
\mathcal{A}
=\{(111), (21), (3), (1^21), (1^3)\}
\end{gather*}
be the set of symbols. 
For a prime number $p$, we write the prime ideal decomposition of $p$ in $K$ as $(p)=\mathfrak{p}_1^{e_1} \mathfrak{p}_2^{e_2} \cdots \mathfrak{p}_r^{e_r}$. 
Then, a cubic field $K$ is said to satisfy a local specification $\mathfrak{a} \in \mathcal{A}$ at $p$ if 
\begin{itemize}
\item[$(\mathrm{a})$] 
for $\mathfrak{a}=(111)$, $p$ is totally splitting in $K$, i.e.\ $(p)=\mathfrak{p}_1 \mathfrak{p}_2 \mathfrak{p}_3$;
\item[$(\mathrm{b})$] 
for $\mathfrak{a}=(21)$, $p$ is partially splitting in $K$, i.e.\ $(p)=\mathfrak{p}_1 \mathfrak{p}_2$; 
\item[$(\mathrm{c})$] 
for $\mathfrak{a}=(3)$, $p$ remains inert in $K$, i.e.\ $(p)=\mathfrak{p}_1$; 
\item[$(\mathrm{d})$] 
for $\mathfrak{a}=(1^21)$, $p$ is partially ramified in $K$, i.e.\ $(p)=\mathfrak{p}_1^2 \mathfrak{p}_2$; 
\item[$(\mathrm{e})$] 
for $\mathfrak{a}=(1^3)$, $p$ is totally ramified in $K$, i.e.\ $(p)=\mathfrak{p}_1^3$.  
\end{itemize}
The symbol $\mathcal{S}$ is used to denote a collection of local specifications with the following data: $(\mathrm{i})$ a finite set $\supp \mathcal{S}$ consisting of prime numbers; $(\mathrm{ii})$ an element $\mathcal{S}_p \in \mathcal{A}$ for each $p \in \supp \mathcal{S}$. 
We say that a cubic field $K$ satisfies the local specifications $\mathcal{S}=(\mathcal{S}_p)_p$ if $K$ satisfies $\mathcal{S}_p$ at $p$ for every $p \in \supp \mathcal{S}$. 
Then we define
\begin{gather*}
L_3^\pm(X,\mathcal{S})
=\left\{K\in L_3^\pm(X) ~\middle|~ \text{$K$ satisfies the local specifications $\mathcal{S}$}\right\}
\end{gather*}
and put $N_3^\pm(X,\mathcal{S})=\# L_3^\pm(X,\mathcal{S})$. 
Remark that the set $\supp \mathcal{S}$ may be empty. 
We define $L_3^\pm(X,\mathcal{S})=L_3^\pm(X)$ in such a case. 
Let $p$ be a prime number and $\mathfrak{a} \in \mathcal{A}$. 
Then we define constants $C_p(\mathfrak{a})$ and $K_p(\mathfrak{a})$ as
\begin{gather}
C_p(\mathfrak{a})
=\frac{1}{1+p^{-1}+p^{-2}} 
\cdot
\begin{cases}
1/6 
& \text{if $\mathfrak{a}=(111)$},\\
1/2 
& \text{if $\mathfrak{a}=(21)$},\\
1/3 
& \text{if $\mathfrak{a}=(3)$},\\
1/p 
& \text{if $\mathfrak{a}=(1^21)$},\\
1/p^2 
& \text{if $\mathfrak{a}=(1^3)$}, 
\end{cases} \label{eq:07250323} \\
K_p(\mathfrak{a})
=\frac{1-p^{-1/3}}{(1-p^{-5/3})(1+p^{-1})} 
\cdot
\begin{cases}
1/6 \cdot (1+p^{-1/3})^3 
& \text{if $\mathfrak{a}=(111)$},\\
1/2 \cdot (1+p^{-1/3})(1+p^{-2/3}) 
& \text{if $\mathfrak{a}=(21)$},\\
1/3 \cdot (1+p^{-1}) 
& \text{if $\mathfrak{a}=(3)$},\\
1/p \cdot (1+p^{-1/3})^2 
& \text{if $\mathfrak{a}=(1^21)$},\\
1/p^2 \cdot (1+p^{-1/3}) 
& \text{if $\mathfrak{a}=(1^3)$}. 
\end{cases}\label{eq:07250324}
\end{gather}
One can check that both $\sum_{\mathfrak{a} \in \mathcal{A}} C_p(\mathfrak{a})$ and $\sum_{\mathfrak{a} \in \mathcal{A}} K_p(\mathfrak{a})$ are equal to $1$. 
We further define 
\begin{gather*}
C(\mathcal{S})
=\prod_{p \in \supp \mathcal{S}} C_p(\mathcal{S}_p)
\quad\text{and}\quad
K(\mathcal{S})
=\prod_{p \in \supp \mathcal{S}} K_p(\mathcal{S}_p), 
\end{gather*}
where the empty product is interpreted as the value $1$. 
Finally, we put 
\begin{gather*}
C^\pm(\mathcal{S})
=C^\pm C(\mathcal{S})
\quad\text{and}\quad
K^\pm(\mathcal{S})
=K^\pm K(\mathcal{S})
\end{gather*}
with $C^+=1$, $C^-=3$, $K^+=1$, $K^-=\sqrt{3}$. 
Then, Roberts \cite{Roberts2001} conjectured that for any local specifications $\mathcal{S}=(\mathcal{S}_p)_p$ the formula
\begin{gather}\label{eq:07250314}
N_3^\pm(X,\mathcal{S})
=C^\pm(\mathcal{S}) \frac{1}{12\zeta(3)} X 
+K^\pm(\mathcal{S}) \frac{4\zeta(1/3)}{5\Gamma(2/3)^3\zeta(5/3)} X^{5/6} 
+o(X^{5/6})
\end{gather}
holds as $X \to\infty$. 
This conjecture was proved to be true by Bhargava--Shankar--Tsimerman \cite{BhargavaShankarTsimerman2013} and Taniguchi--Thorne \cite{TaniguchiThorne2013} independently. 
More precisely, it was shown that
\begin{align*}
&N_3^\pm(X,\mathcal{S})
-C^\pm(\mathcal{S}) \frac{1}{12\zeta(3)} X 
-K^\pm(\mathcal{S}) \frac{4\zeta(1/3)}{5\Gamma(2/3)^3\zeta(5/3)} X^{5/6} \\
&\ll_\epsilon X^{7/9+\epsilon} \prod_{p \in \supp \mathcal{S}} p^{e_p}, 
\end{align*}
where $e_p=8/9$ for $\mathcal{S}_p=(111), (21), (3)$ and $e_p=16/9$ for $\mathcal{S}_p=(1^21), (1^3)$. 

We associate a symbol $\mathfrak{a} \in \mathcal{A}$ with a diagonal matrix $A_\mathfrak{a} \in M_2(\mathbb{C})$ by putting 
\begin{gather}\label{eq:07261556}
A_\mathfrak{a}
=
\begin{cases}
\diag(1,1) 
& \text{if $\mathfrak{a}=(111)$},\\
\diag(1,-1) 
& \text{if $\mathfrak{a}=(21)$},\\
\diag(\omega,\overline{\omega}) 
& \text{if $\mathfrak{a}=(3)$},\\
\diag(1,0) 
& \text{if $\mathfrak{a}=(1^21)$},\\
\diag(0,0) 
& \text{if $\mathfrak{a}=(1^3)$},
\end{cases}
\end{gather}
where $\omega$ is a primitive cube root of unity.
Then the Artin $L$-function $L(s,\rho_K)$ has the Euler product representation 
\begin{gather}\label{eq:07252118}
L(s,\rho_K)
=\prod_{p} \det\left(I-p^{-s} A_p(K)\right)^{-1} 
\end{gather}
for $\RE(s)>1$ by definition. 
Here, $I \in M_2(\mathbb{C})$ is the identity matrix, and $A_p(K)$ is determined by $A_p(K)=A_\mathfrak{a}$ if $K$ satisfies a local specification $\mathfrak{a}$ at $p$. 
According to formula \eqref{eq:07250314}, we define $\mathcal{X}=(\mathcal{X}_p)_p$ and $\mathcal{Y}=(\mathcal{Y}_p)_p$ as two sequences of independent random elements on the set $\{A_\mathfrak{a} \mid \mathfrak{a} \in \mathcal{A} \}$ such that 
\begin{gather}\label{eq:07252104}
\mathbb{P}\left(\mathcal{X}_p=A_\mathfrak{a}\right)
=C_p(\mathfrak{a})
\quad\text{and}\quad
\mathbb{P}\left(\mathcal{Y}_p=A_\mathfrak{a}\right)
=K_p(\mathfrak{a}), 
\end{gather}
where $C_p(\mathfrak{a})$ and $K_p(\mathfrak{a})$ are as in \eqref{eq:07250323} and \eqref{eq:07250324}, respectively. 
Here, we denote by $\mathbb{P}(E)$ the probability of an event $E$. 
The random Euler products $L(s,\mathcal{X})$ and $L(s,\mathcal{Y})$ are defined as
\begin{gather}\label{eq:07271414}
L(s,\mathcal{X})
=\prod_{p} \det\left(I-p^{-s} \mathcal{X}_p\right)^{-1} 
\quad\text{and}\quad
L(s,\mathcal{Y})
=\prod_{p} \det\left(I-p^{-s} \mathcal{Y}_p\right)^{-1}. 
\end{gather}
We show later that the former infinite product converges for $\RE(s)>1/2$, and the latter does for $\RE(s)>2/3$. 
The first main result is about the density functions for these random Euler products. 

\begin{theorem}\label{thm:2.1}
For $\sigma>1/2$, there exists a non-negative $C^\infty$-function $C_\sigma$ such that 
\begin{gather*}
\mathbb{P}\left( \log{L}(\sigma,\mathcal{X}) \in A \right)
=\int_{A} C_\sigma(x) \,\frac{dx}{\sqrt{2\pi}}
\end{gather*}
holds for all $A \in \mathcal{B}(\mathbb{R})$. 
Furthermore, it satisfies the following properties.  
\begin{itemize}
\item[$(\mathrm{i})$]
If $1/2<\sigma\leq1$, we have $\supp C_\sigma=\mathbb{R}$, that is, $C_\sigma(x)$ is not identically zero in any interval on $\mathbb{R}$. 
\item[$(\mathrm{ii})$]
If $\sigma>1$, the function $C_\sigma$ is compactly supported.  
\item[$(\mathrm{iii})$]
Let $\sigma>1/2$. 
Then the integral 
\begin{gather*}
\int_{-\infty}^{\infty} e^{ax} C_\sigma(x) \,\frac{dx}{\sqrt{2\pi}}
\end{gather*}
is finite for any $a>0$. 
\end{itemize}
Similarly, for $\sigma>2/3$, there exists a non-negative $C^\infty$-function $K_\sigma$ such that 
\begin{gather*}
\mathbb{P}\left( \log{L}(\sigma,\mathcal{Y}) \in A \right)
=\int_{A} K_\sigma(x) \,\frac{dx}{\sqrt{2\pi}}
\end{gather*}
holds for all $A \in \mathcal{B}(\mathbb{R})$. 
Furthermore, it satisfies the following properties. 
\begin{itemize}
\item[$(\mathrm{i'})$]
If $2/3<\sigma\leq1$, we have $\supp K_\sigma=\mathbb{R}$, that is, $K_\sigma(x)$ is not identically zero in any interval on $\mathbb{R}$. 
\item[$(\mathrm{ii'})$]
If $\sigma>1$, the function $K_\sigma$ is compactly supported.  
\item[$(\mathrm{iii'})$]
Let $\sigma>2/3$. 
Then the integral 
\begin{gather*}
\int_{-\infty}^{\infty} e^{ax} K_\sigma(x) \,\frac{dx}{\sqrt{2\pi}}
\end{gather*}
is finite for any $a>0$. 
\end{itemize}
\end{theorem}

Let $f$ be a function in $L^1(\mathbb{R})$ such that the integral 
\begin{gather*}
\int_{-\infty}^{\infty} e^{ax} f(x) \,\frac{dx}{\sqrt{2\pi}}
\end{gather*}
is finite for any $a>0$. 
Throughout this paper, the Fourier--Laplace transform of $f$ is defined as
\begin{gather*}
\widetilde{f}(z)
=\int_{-\infty}^{\infty} f(x) e^{izx} \,\frac{dx}{\sqrt{2\pi}}, 
\end{gather*}
which presents a holomorphic function on the whole complex plane. 
By the independence of  $\mathcal{X}=(\mathcal{X}_p)_p$, the Fourier--Laplace transform of $C_\sigma$ is represented as  
\begin{align*}
\widetilde{C}_\sigma(z)
&=\prod_{p} \left( \sum_{\mathfrak{a} \in \mathcal{A}} 
C_p(\mathfrak{a}) \det\left(I-p^{-s} A_\mathfrak{a} \right)^{-iz}\right) \\
&=\prod_{p} \frac{1}{1+p^{-1}+p^{-2}} 
\Biggl[ \frac{1}{6} \left(1-\frac{1}{p^{\sigma}}\right)^{-2iz} 
+\frac{1}{2} \left(1-\frac{1}{p^{2\sigma}}\right)^{-iz} \\
&\qquad\qquad
+\frac{1}{3} \left(1+\frac{1}{p^{\sigma}}+\frac{1}{p^{2\sigma}}\right)^{-iz} 
+\frac{1}{p} \left(1-\frac{1}{p^{\sigma}}\right)^{-iz}
+\frac{1}{p^2} \Biggr] 
\end{align*} 
for any $z \in \mathbb{C}$. 
In addition, $\widetilde{K}_\sigma(z)$ has a similar infinite product representation. 

The Grand Riemann Hypothesis (GRH) is used in the sense that every $L(s,\rho_K)$ has no zeros in the half-plane $\RE(s)>1/2$. 
We may use the following zero density estimate for $L(s,\rho_K)$ as well. 
Let $N(\alpha,T;\rho_K)$ count the number of zeros of $L(s,\rho_K)$ satisfying $\RE(\rho) \geq \alpha$ and $|\IM(\rho)| \leq T$ with multiplicity. 
Then there exists an absolute constant $1/2<\sigma_1<1$ such that the estimate
\begin{gather}\label{eq:07251505}
\sum_{K \in L_3^\pm(X)} N(\sigma_1, (\log{X})^3; \rho_K)
\ll X^{1-\delta}
\end{gather}
holds with an absolute constant $\delta>0$. 
Remark that GRH ensures the truth of \eqref{eq:07251505} for any $1/2<\sigma_1<1$. 
In Section \ref{sec:3.1}, we also see that one can choose any $\sigma_1>7/8$ unconditionally by using the zero-density estimate of Kowalski--Michel \cite{KowalskiMichel2002}. 
Let $E(X)$ be a subset of $L_3^\pm(X)$ defined as
\begin{gather}\label{eq:07252237}
E(X)
=\left\{ K \in L_3^\pm(X) ~\middle|~ 
\begin{array}{l}
\text{There exists a zero $\rho$ of $L(s,\rho_K)$ such that} \\
\text{$\RE(\rho) \geq \sigma_1$ and $|\IM(\rho)| \leq (\log{X})^{3}$.}
\end{array}
\right\}. 
\end{gather}
Then we define $E_\sigma(X)=E(X)$ for $\sigma_1<\sigma<1$ and $E_\sigma(X)=\emptyset$ for $\sigma \geq1$. 
If we assume \eqref{eq:07251505}, then we obtain the estimate $\# E_\sigma(X) = O(X^{1-\delta})$ for $\sigma>\sigma_1$. 
The second and third results are relations between the complex moment $M_{z, \sigma}^{\pm}(X)$ of \eqref{eq:07250226} and the integrals involving the density functions of Theorem \ref{thm:2.1}. 

\begin{theorem}\label{thm:2.2}
Let $\sigma_1$ be a real number for which \eqref{eq:07251505} holds. 
Then there exists an absolute constant $\delta>0$ such that 
\begin{gather}\label{eq:07251529}
M_{z, \sigma}^{\pm}(X)
=C^\pm \frac{1}{12\zeta(3)} X 
\int_{-\infty}^{\infty} e^{zx} C_\sigma(x) \,\frac{dx}{\sqrt{2\pi}} 
+O\left(X\exp\left(-\delta \frac{\log{X}}{\log\log{X}}\right)\right)
\end{gather}
holds for $\sigma>\sigma_1$ with $z \in \mathbb{C}$ satisfying $|z| \leq b_\sigma R_{\sigma}(X)$, where $b_\sigma$ is a positive constant, and $R_{\sigma}(X)$ is defined as
\begin{gather}\label{eq:07281607}
R_{\sigma}(X)
=
\begin{cases}
(\log{X}) (\log\log{X})^{-2} 
& \text{if $\sigma \geq1$}, \\
(\log{X})^{(\sigma-\sigma_1)/(1-\sigma_1)} (\log\log{X})^{-2} 
& \text{if $\sigma_1<\sigma<1$}. 
\end{cases}
\end{gather}
The implied constant in \eqref{eq:07251529} depends only on $\sigma$. 
\end{theorem}

\begin{theorem}\label{thm:2.3}
Assume GRH and the upper bound
\begin{align}\label{eq:07251548}
&N_3^\pm(X,\mathcal{S})
-C^\pm(\mathcal{S}) \frac{1}{12\zeta(3)} X 
-K^\pm(\mathcal{S}) \frac{4\zeta(1/3)}{5\Gamma(2/3)^3\zeta(5/3)} X^{5/6} \\
&\ll_\epsilon X^{\alpha+\epsilon} \prod_{p \in \supp \mathcal{S}} p^{\beta} \nonumber
\end{align}
for each $\epsilon>0$ with some constants $\alpha$ and $\beta$ such that $0<\alpha<5/6$. 
Put 
\begin{gather}\label{eq:07302054}
\sigma_2
=\frac{(5-6\alpha)+(2+2\beta)}{12(1-\alpha)}. 
\end{gather}
Then there exists a constant $\delta=\delta(\sigma)>0$ such that 
\begin{align}\label{eq:07251545}
M_{z, \sigma}^{\pm}(X)
&=C^\pm \frac{1}{12\zeta(3)} X 
\int_{-\infty}^{\infty} e^{zx} C_\sigma(x) \,\frac{dx}{\sqrt{2\pi}} \\
&\quad+K^\pm \frac{4\zeta(1/3)}{5\Gamma(2/3)^3\zeta(5/3)} X^{5/6}
\int_{-\infty}^{\infty} e^{zx} K_\sigma(x) \,\frac{dx}{\sqrt{2\pi}} 
+O\left(X^{5/6-\delta}\right) \nonumber
\end{align}
holds for $\sigma>\max(\sigma_2,2/3)$ with $z \in \mathbb{C}$ satisfying $|z| \leq \tilde{b}_\sigma \tilde{R}_{\sigma}(X)$, where $\tilde{b}_\sigma$ is a positive constant, and $\tilde{R}_{\sigma}(X)$ is defined as
\begin{gather}\label{eq:07301552}
\tilde{R}_\sigma(X)
=(\log{X})^{a(\sigma)} (\log\log{X})^{-1} 
\end{gather}
with a small constant $a(\sigma)>0$. 
The implied constant in \eqref{eq:07251545} depends only on $\sigma$.
\end{theorem}

Recall that \eqref{eq:07251548} is true for $\alpha=7/9$ and $\beta=16/9$ by Taniguchi--Thorne \cite{TaniguchiThorne2013}. 
Thus one can take $\sigma_2=53/24$ at least in Theorem \ref{thm:2.3}. 
Recently, Bhargava--Taniguchi--Thorne \cite{BhargavaTaniguchiThorne2021+} improved this estimate by $\alpha=\beta=2/3$, and we have $\sigma_2=13/12$ further. 
On the other hand, Cho--Fiorilli--Lee--S\"{o}dergren \cite{ChoFiorilliLeeSodergren2022} proved the lower bound $\alpha+\beta \geq1/2$ under GRH, when $\supp \mathcal{S}=\{p\}$. 
Note that $\alpha+\beta=1/2$ is equivalent to $\sigma_2=2/3$, which means that \eqref{eq:07251545} holds for $\sigma>2/3$. 

\begin{corollary}\label{cor:2.1}
There exists an absolute constant $\delta>0$ such that  
\begin{align*}
\sum_{K \in L_3^\pm(X)} (h_K R_K)^r
&=C^\pm \frac{1}{12\zeta(3)} \frac{2X^{r/2+1}}{(r+2) (D^\pm)^r} 
\int_{-\infty}^{\infty} e^{rx} C_1(x) \,\frac{dx}{\sqrt{2\pi}} \\
&\quad+O\left(X^{r/2+1} \exp\left(-\delta \frac{\log{X}}{\log\log{X}}\right)\right)
\end{align*}
holds for any real number $r>-2$, where the implied constant depends only on $r$. 
\end{corollary}

\begin{corollary}\label{cor:2.2}
Assume GRH and upper bound \eqref{eq:07251548} with some constants $\alpha$ and $\beta$ such that $3\alpha+\beta<5/2$. 
Then there exists an absolute constant $\delta>0$ such that  
\begin{align*}
\sum_{K \in L_3^\pm(X)} (h_K R_K)^r
&=C^\pm \frac{1}{12\zeta(3)} \frac{2X^{r/2+1}}{(r+2) (D^\pm)^r} 
\int_{-\infty}^{\infty} e^{rx} C_1(x) \,\frac{dx}{\sqrt{2\pi}} \\
&\quad+K^\pm \frac{4\zeta(1/3)}{5\Gamma(2/3)^3\zeta(5/3)} \frac{5X^{r/2+5/6}}{(3r+5) (D^\pm)^r} 
\int_{-\infty}^{\infty} e^{rx} K_1(x) \,\frac{dx}{\sqrt{2\pi}} \\
&\quad+O\left(X^{r/2+5/6} \exp\left(-\delta \frac{\log{X}}{\log\log{X}}\right)\right)
\end{align*}
holds for any real number $r>-5/3$, where the implied constant depends only on $r$. 
\end{corollary}

We have more consequences of Theorem \ref{thm:2.2} by the methods of probability theory.  
Define the quantity $D_\sigma^\pm(X; a)$ as 
\begin{gather*}
D_\sigma^\pm(X; a)
=\frac{\# \left\{ K \in L_3^\pm(X) ~\middle|~ L(\sigma,\rho_K) \leq e^a \right\}}{N_3^\pm(X)}
-\int_{-\infty}^{a} C_\sigma(x) \,\frac{dx}{\sqrt{2\pi}}
\end{gather*}
for $\sigma>1/2$ and $a \in \mathbb{R}$. 
Theorem \ref{thm:1.1} asserts that $D_\sigma^\pm(X; a) \to 0$ holds for $\sigma>7/8$ as $X \to\infty$. 
We further evaluate the decay of $D_\sigma^\pm(X; a)$ as follows. 

\begin{theorem}\label{thm:2.4}
Let $\sigma_1$ be a real number for which \eqref{eq:07251505} holds. 
Then we obtain
\begin{gather*}
\sup_{a \in \mathbb{R}} \left|D_\sigma^\pm(X; a)\right|
\ll \frac{1}{R_\sigma(X)}, 
\end{gather*}
where $R_\sigma(X)$ is as in Theorem \ref{thm:2.2}. 
\end{theorem}

\begin{remark}\label{rem:2.1}
Compared with the results in \cite{LamzouriLesterRadziwill2019, Mine2022}, it might be possible that the function $R_\sigma(X)$ of Theorems \ref{thm:2.2} and \ref{thm:2.4} is improved by
\begin{gather*}
R_\sigma(X)
=
\begin{cases}
(\log{X}) (\log\log{X})^{-1} 
& \text{if $\sigma>1$, } \\
(\log{X}) (\log\log{X} \log\log\log{X})^{-1} 
& \text{if $\sigma=1$, } \\
(\log{X})^{\sigma} 
& \text{if $\sigma_1<\sigma<1$. } 
\end{cases}
\end{gather*}
However, we recall that the fact 
\begin{gather*}
\lim_{T \to\infty} \frac{1}{T} \int_{T}^{2T} p^{-it} \,dt
=\mathbb{E}\left[ \mathcal{T}_p \right]
=0 
\end{gather*}
is used in the method of \cite{LamzouriLesterRadziwill2019}, where $\mathcal{T}_p$ is a random variable uniformly distributed on the unit circle in $\mathbb{C}$. 
In the case of Artin $L$-function $L(s,\rho_K)$, we have a similar fact
\begin{gather*}
\lim_{X \to\infty} \frac{1}{N_3^\pm(X)} \tr(A_p(K))
=\mathbb{E}\left[ \tr(\mathcal{X}_p) \right]
=\frac{p^{-1}}{1+p^{-1}+p^{-2}}, 
\end{gather*}
where $A_p(K)$ is the matrix as in \eqref{eq:07252118}, and $\mathcal{X}_p$ is a random element satisfying \eqref{eq:07252104}. 
However, the non-vanishing of $\mathbb{E}\left[ \tr(\mathcal{X}_p) \right]$ prevents us from applying the totally same method as in \cite{LamzouriLesterRadziwill2019}. 
See also \cite[condition $(1.10)$ and Lemma 2.10]{Mine2022}. 
Therefore, we adopt in this paper not the method of \cite{LamzouriLesterRadziwill2019, Mine2022} but another method partially motivated by \cite{CogdellMichel2004, Lamzouri2011b, IharaMatsumoto2011b}. 
\end{remark}

Let $\sigma>1/2$ be a fixed real number. 
We see that $L(\sigma,\rho_K)$ might be zero or negative if we do not assume GRH. 
Thus, we define 
\begin{gather}\label{eq:07252224}
A_\sigma(X)
=\left\{ K \in L_3^\pm(X) ~\middle|~ \text{$L(\alpha,\rho_K)=0$ for some $\alpha \geq \sigma$}\right\}
\end{gather}
for $1/2<\sigma<1$ and $A_\sigma(X)=\emptyset$ for $\sigma \geq1$. 
Throughout this paper, we write 
\begin{gather*}
\sideset{}{'} \sum_{K \in L_3^\pm(X)}
=\sum_{K \in L_3^\pm(X) \setminus A_\sigma(X)}. 
\end{gather*}
Note that one can define $\log{L}(\sigma,\rho_K)$ for $K \notin A_\sigma(X)$. 
Next, we define classes of test functions to describe the statement of an analogue of Theorem \ref{thm:1.3}. 
Let $C(\mathbb{R})$ be the class of all continuous functions on $\mathbb{R}$. 
Then we define the subclasses of $C(\mathbb{R})$ as 
\begin{align*}
C^{\exp}(\mathbb{R})
&=\left\{ \Phi \in{C}(\mathbb{R}) ~\middle|~ \text{$\Phi(x) \ll e^{a|x|}$ with some $a>0$} \right\}, \\
C_b(\mathbb{R})
&=\left\{ \Phi \in{C}(\mathbb{R}) ~\middle|~ \text{$\Phi$ is bounded} \right\}. 
\end{align*}
We further define
\begin{gather*}
I(\mathbb{R})
=\left\{ 1_A ~\middle|~ \text{$A$ is a continuity set of $\mathbb{R}$} \right\}, 
\end{gather*}
where $1_A$ is the indicator function of a set $A \subset \mathbb{R}$, and a Borel set $A$ is called a continuity set of $\mathbb{R}$ if its boundary $\partial A$ has Lebesgue measure zero in $\mathbb{R}$. 

\begin{theorem}\label{thm:2.5}
Let $\sigma_1$ be a real number for which \eqref{eq:07251505} holds. 
Then the limit formula
\begin{gather}\label{eq:07310245}
\lim_{X \to\infty} \frac{1}{N_3^\pm(X)} 
\sideset{}{'} \sum_{K \in L_3^\pm(X)} \Phi\left(\log{L}(\sigma,\rho_K)\right)
=\int_{-\infty}^{\infty} \Phi(u) C_\sigma(x) \,\frac{dx}{\sqrt{2\pi}}
\end{gather}
holds in the following cases: 
\begin{itemize}
\item 
$\sigma>1$ and $\Phi \in C(\mathbb{R}) \cup I(\mathbb{R})$; 
\item 
$\sigma=1$ and $\Phi \in C^{\exp}(\mathbb{R}) \cup I(\mathbb{R})$; 
\item 
$\sigma_1<\sigma<1$ and $\Phi \in C_b(\mathbb{R}) \cup I(\mathbb{R})$ without assuming GRH; 
\item 
$\sigma_1<\sigma<1$ and $\Phi \in C^{\exp}(\mathbb{R}) \cup I(\mathbb{R})$ if we assume GRH.  
\end{itemize}
\end{theorem}

We obtain several results on the logarithmic derivative $(L'/L)(s,\rho_K)$ similar to the above. 
They are presented in the appendix.

\section{Preliminaries}\label{sec:3}

\subsection{Properties of the Artin $L$-function}\label{sec:3.1}
As we mentioned in Section \ref{sec:1.1}, there exists a cuspidal representation $\pi$ such that $L(s,\rho_K)=L(s,\pi)$ holds. 
Hence the Artin $L$-function $L(s,\rho_K)$ is continued to an entire function. 
The complete $L$-function 
\begin{gather*}
\Lambda(s,\rho_K)
=|d_K|^{s/2} \gamma(s,\rho_K) L(s,\rho_K)
\end{gather*}
satisfies the functional equation $\Lambda(s,\rho_K)=\Lambda(1-s,\rho_K)$, where the gamma factor is given by 
\begin{gather*}
\gamma(s,\rho_K)
=
\begin{cases}
\pi^{-s} \Gamma\left(\frac{s}{2}\right)^2
&\quad \text{if $d_K>0$}, \\
\pi^{-s} \Gamma\left(\frac{s}{2}\right)\Gamma\left(\frac{s+1}{2}\right)
&\quad \text{if $d_K<0$}. 
\end{cases}
\end{gather*}
Note that $\gamma(s,\rho_K)$ are common over the family ${L}_3^+(X)$ or ${L}_3^-(X)$. 
Then, we obtain the following zero density estimate for $L(s,\rho_K)$. 

\begin{lemma}\label{lem:3.1}
Let $X \geq1$ and $T \geq2$. 
For any $C_0>6$, we have 
\begin{gather*}
\sum_{K \in L_3^\pm(X)} N(\alpha,T;\rho_K) 
\ll T^A X^{C_0(1-\alpha)/(2\alpha-1)}
\end{gather*}
for any $\alpha \geq3/4$, where $A>0$ is an absolute constant.  
The implied constant depends only on the choice of $C_0$. 
\end{lemma}

\begin{proof}
Let $S^\pm(X)$ be the set of all cuspidal representations $\pi$ of $GL_2(\mathbb{A}_\mathbb{Q})$ such that $L(s,\pi)=L(s,\rho_K)$ is satisfied for some $K \in L_3^\pm(X)$. 
Note the Ramanujan--Petersson conjecture is valid for any $\pi \in S^\pm(X)$ since all local roots satisfy $|\alpha_j(p)| \leq1$ by \eqref{eq:07252118}. 
Next, the conductor of $\pi \in S^\pm(X)$ satisfies $\Cond(\pi) \leq X$ due to $\Cond(\pi)=|d_K|$ if $L(s,\pi)=L(s,\rho_K)$. 
Furthermore, we obtain $\# S^\pm(X) \ll X$, which is deduced from the fact that there exists a one-to-one correspondence between $S^\pm(X)$ and $L_3^\pm(X)$. 
Finally, the gamma factors in the functional equations of $L(s,\pi)$ are common in $\pi \in S^\pm(X)$. 
From the above, the desired estimate follows directly if we apply the zero density estimate of Kowalski and Michel \cite[Theorem 2]{KowalskiMichel2002} to the family $S^\pm(X)$. 
\end{proof}

We check that estimate \eqref{eq:07251505} holds for any $7/8<\sigma_1<1$. 
Taking $C_0=6+\delta$ with $\delta=(8\sigma_1-7)/2>0$, we obtain 
\begin{gather*}
\frac{C_0(1-\sigma_1)}{2\sigma_1-1}
=1-\frac{(8+\delta) \sigma_1-(7+\delta)}{2\sigma_1-1}
\leq 1-\delta.
\end{gather*}
Hence Lemma \ref{lem:3.1} ensures the validity of \eqref{eq:07251505} with $7/8<\sigma_1<1$. 

We often consider the logarithms of $L(s,\rho_K)$ for non-real variables $s$. 
The branch of $\log{L(s,\rho_K)}$ is determined as follows. 
First, we define 
\begin{gather}\label{eq:07252233}
\log{L(s,\rho_K)}
=\sum_{p} \sum_{m=1}^{\infty} \frac{\tr(A_p(K)^m)}{m} p^{-ms}
=\sum_{n=1}^{\infty} \frac{\Lambda_K(n)}{\log{n}} n^{-s}
\end{gather}
for $\RE(s)>1$ from Euler product \eqref{eq:07252118}. 
Here, the coefficient $\Lambda_K(n)$ is calculated as $\Lambda_K(p^m)=\tr(A_p(K)^m) \log{p}$ and $\Lambda_K(n)=0$ unless $n$ is a prime power. 
Note that we have $|\Lambda_K(n)| \leq 2\Lambda(n)$ with the usual von Mangoldt function $\Lambda(n)$. 
Let $D$ denote the right half-plane $\{\sigma+i\tau \mid \sigma>1/2\}$. 
We define
\begin{gather*}
G_K
=D \setminus \bigcup_{\RE(\rho)>1/2} \{\sigma+i\IM(\rho) \mid 1/2<\sigma \leq \RE(\rho)\},
\end{gather*}
where $\rho$ runs through all possible zeros of $L(s,\rho_K)$ with $\RE(\rho)>1/2$. 
Then we extend $\log{L}(s,\rho_K)$ for $s \in G_K$ by the analytic continuation along the horizontal path from right. 
The region $G_K$ is adequate to define $\log{L(s,\rho_K)}$ as a holomorphic function, and we note that the subset $A_\sigma(X)$ defined as \eqref{eq:07252224} is represented as  
\begin{gather*}
A_\sigma(X) 
=\left\{ K \in L_3^\pm(X) ~\middle|~ \sigma \notin G_K \right\}.
\end{gather*}

\begin{lemma}\label{lem:3.2}
Assume GRH, and let $\sigma_0>1/2$ be a real number. 
Take a complex number $s=\sigma+it$ such that $\sigma \geq \sigma_0$ and $|t| \leq (\log{X})^2$. 
For any cubic field $K \in L_3^\pm(X)$, we have 
\begin{gather*}
\log{L(s,\rho_K)}
\ll \frac{(\log{X})^{2-2\sigma}}{\log\log{X}} 
+\log\log{X},
\end{gather*}
where the implied constant depends only on the choice of $\sigma_0$. 
\end{lemma}

\begin{proof}
Let $q(s,\rho_K)=|d_K| (|s|+3) (|s+1|+3)$ be the analytic conductor of $L(s,\rho_K)$ in the sense of Iwaniec--Kowalski \cite{IwaniecKowalski2004}. 
Assuming GRH, we have 
\begin{gather*}
\log{L(s,\rho_K)} 
\ll \frac{(\log{q(s,\rho_K)})^{2-2\sigma}}{(2\sigma-1) \log\log{q}(s,\rho_K)} 
+\log\log{q}(s,\rho_K)
\end{gather*}
for any $s=\sigma+it$ with $1/2<\sigma \leq5/4$ by \cite[Theorems 5.19]{IwaniecKowalski2004}. 
Thus the result for $\sigma_0 \leq \sigma \leq 5/4$ follows. 
The case $\sigma \geq 5/4$ is trivial. 
Indeed, we obtain the inequality $|\log{L(s,\rho_K)}| \leq 2\log{\zeta(5/4)}$ by \eqref{eq:07252233} and $|\Lambda_K(n)| \leq 2\Lambda(n)$. 
\end{proof}

\begin{lemma}\label{lem:3.3}
Let $\sigma_1$ be a real number for which \eqref{eq:07251505} holds. 
Take a complex number $s=\sigma+it$ such that $\sigma \geq \sigma_1+2(\log\log{X})^{-1}$ and $|t| \leq (\log{X})^2$. 
For any cubic field $K \in L_3^\pm(X) \setminus E(X)$, we have 
\begin{gather*}
\log{L(s,\rho_K)}
\ll (\log\log{X}) (\log{X})^{(1-\sigma)/(1-\sigma_1)} 
+\log\log{X}, 
\end{gather*}
where $E(X)$ is the subset defined by \eqref{eq:07252237}, and the implied constant is absolute. 
\end{lemma}

\begin{proof}
The proof is based on the method of Barban \cite[Lemma 3]{Barban1962}. 
For simplicity, we write $\kappa=(\log\log{X})^{-1}$. 
If $s=\sigma+it \in \mathbb{C}$ satisfies $\sigma \geq 1+\kappa/2$, then we obtain 
\begin{gather*}
\log{L(s,\rho_K)}
\ll \sum_{n=1}^{\infty} \frac{\Lambda(n)}{\log{n}} n^{-\sigma}
\ll \log\log{X}
\end{gather*}
by formula \eqref{eq:07252233}. 
Therefore the result follows in this case, and we let $s=\sigma+it$ with $\sigma_1+2\kappa \leq \sigma \leq 1+\kappa/2$ and $|t| \leq (\log{X})^2$ below. 
Let $z_0=\kappa^{-1}+\kappa+it \in \mathbb{C}$ and $R=\kappa^{-1}+\kappa-\sigma_1>0$. 
If $K \in L_3^\pm(X) \setminus E(X)$, then the function $(L'/L)(z,\rho_K)$ is holomorphic on $|z-z_0|<R$ since $L(z,\rho_K)$ has no zeros in this disk. 
We define 
\begin{gather*}
M(r)
=\max_{|z-z_0|=r} \left|\frac{L'}{L}(z,\rho_K)\right|
\end{gather*}
for $0<r<R$. 
Let $r_1=\kappa^{-1}-1$, $r_2=\kappa^{-1}+\kappa-\sigma$ and $r_3=\kappa^{-1}-\sigma_1$. 
Then we have $0<r_1<r_2<r_3<R$. 
As a consequence of the Hadamard three circles theorem, we obtain the inequality
\begin{gather}\label{eq:07252344}
M(r_2) 
\leq M(r_1)^{1-a} M(r_3)^a,
\end{gather}
where $a=\log(r_2/r_1) / \log(r_3/r_1)$. 
We evaluate $M(r_1)$. 
Since we have $\RE(z) \geq 1+\kappa$ on the circle $|z-z_0|=r_1$, the estimate $(L'/L)(z,\rho_K) \ll \log\log{X}$ follows. 
Therefore we obtain 
\begin{gather}\label{eq:07252345}
M(r_1)
\leq A \log\log{X},
\end{gather}
where $A \geq1$ is an absolute constant. 
Next, we let $z=x+iy$ be a complex number on the circle $|z-z_0|=r_3$. 
By \cite[Proposition 5.7 (2)]{IwaniecKowalski2004}, we have 
\begin{gather}\label{eq:07252336}
\frac{L'}{L}(z,\rho_K)
=\sum_{|z-\rho|<1} \frac{1}{z-\rho} 
+O\left(\log(|d_K|(|y|+3))\right)
\end{gather}
with an absolute implied constant, where $\rho$ runs through zeros of $L(s,\rho_K)$. 
Note that the distance between $z$ and $\rho$ is at least $\kappa$ if $K \in L_3^\pm(X) \setminus E(X)$. 
Furthermore, the number of zeros with $|z-\rho|<1$ is evaluate as $\ll \log{X}$ by \cite[Proposition 5.7 (1)]{IwaniecKowalski2004}. 
As a result, we obtain $(L'/L)(z,\rho_K) \ll (\log\log{X}) (\log{X})$ by formula \eqref{eq:07252336}. 
Therefore we derive
\begin{gather}\label{eq:07252346}
M(r_3)
\leq B (\log\log{X}) (\log{X}),
\end{gather}
where $B \geq1$ is an absolute constant. 
Inserting \eqref{eq:07252345} and \eqref{eq:07252346} to \eqref{eq:07252344}, we obtain 
\begin{gather*}
M(r_2)
\leq A^{1-a} B^a (\log\log{X})(\log{X})^a.
\end{gather*}
Note that $0<a<1$ holds since $r_1<r_2<r_3$. 
Thus we have $A^{1-a} B^a \leq AB$. 
By the definition of $a$, we further obtain 
\begin{gather*}
a
=\frac{1-\sigma}{1-\sigma_1}+O\left(\frac{1}{\log\log{X}}\right). 
\end{gather*}
Hence we arrive at the upper bound
\begin{gather}\label{eq:07260204}
\frac{L'}{L}(z,\rho_K)
\ll (\log\log{X}) (\log{X})^{(1-\sigma)/(1-\sigma_1)}
\end{gather}
in the disk $|z-z_0| \leq r_2$, where the implied constant is absolute. 
By the choice of the branch, the relation
\begin{gather*}
\log{L(s,\rho_K)}
=\log{L}(s_0,\rho_K)
-\int_{\sigma}^{1+\kappa/2} \frac{L'}{L}(x+it,\rho_K) \,dx
\end{gather*}
holds with $s_0=1+\kappa/2+it$. 
We have $\log{L}(s_0,\rho_K) \ll \log\log{X}$ as before. 
Furthermore, since the horizontal path from $s=\sigma+it$ to $s_0$ is included in $|z-z_0| \leq r_2$, we obtain
\begin{gather*}
\int_{\sigma}^{1+\kappa/2} \frac{L'}{L}(x+it,\rho_K) \,dx
\ll (\log\log{X}) (\log{X})^{(1-\sigma)/(1-\sigma_1)}
\end{gather*}
by \eqref{eq:07260204}. 
Hence the desired result follows. 
\end{proof}

\begin{lemma}\label{lem:3.4}
For any cubic field $K \in L_3^\pm(X)$, we have 
\begin{gather*}
\log{L(1,\rho_K)}
\ll \log\log{X},
\end{gather*}
where the implied constant is absolute.
\end{lemma}

\begin{proof}
Explicit upper and lower bounds for the value ${L}(1,\rho_K)=\Res_{s=1} \zeta_K(s)$ were obtained by Louboutin \cite{Louboutin2001, Louboutin2005}. 
These results yields the inequalities
\begin{gather*}
\frac{c_K \alpha_K \exp(-c_K \alpha_K/2)}{2 \log|d_K|}
\leq L(1,\rho_K)
\leq \left(\frac{e}{4} \log|d_K|\right)^2, 
\end{gather*}
where $c_K=2(\sqrt{2}-1)^2$ and $\alpha_K=\log|d_K|/\log|d_{\widehat{K}}|$. 
We know that $1/4 \leq \alpha_K \leq 1/2$ holds. 
Hence we conclude that $-\log\log{X}-4 \leq \log{L(1,\rho_K)} \leq 2\log\log{X}+2$. 
\end{proof}

\subsection{The $z$-th divisor function}\label{sec:3.2}
Let $F(w)=\log(1-w)^{-1}$ with $|w|<1$. 
For $z \in \mathbb{C}$, we define $H_r(z)$ as the coefficients of the power series
\begin{gather}\label{eq:07261602}
\exp(z F(w))
=(1-w)^{-z}
=\sum_{r=0}^{\infty} H_r(z) w^r. 
\end{gather}
Then they can be explicitly calculated as $H_0(z)=1$ and 
\begin{gather}\label{eq:07261502}
H_r(z)
=\frac{1}{r!} z(z+1)\cdots(z+r-1)
\end{gather}
for $r \geq1$. 
The $z$-th divisor function $d_z(n)$ is the multiplicative function determined by $d_z(p^r)=H_r(z)$ for a prime number $p$. 
By definition, $d_z(n)$ satisfies
\begin{gather}\label{eq:07261516}
\zeta(s)^z
=\sum_{n=1}^{\infty} d_z(n) n^{-s}
\end{gather}
for any $z \in \mathbb{C}$ and $\RE(s)>1$. 
If $k$ is a positive integer, then $d_k(n)$ is, as usual, equal to the number of representations of $n$ as the product of $k$ positive integers.

\begin{lemma}\label{lem:3.5}
Let $z \in \mathbb{C}$ and $k=\lfloor{|z|}\rfloor+1$, where $\lfloor{x}\rfloor$ indicates the largest integer less than or equal to $x$. 
Then we have $|d_z(n)| \leq d_k(n)$ for any $n \geq1$. 
\end{lemma}

\begin{proof}
Since $d_z(n)$ and $d_k(n)$ are multiplicative, it is sufficient to show the inequality in the case $n=p^r$. 
By \eqref{eq:07261502}, we have
\begin{align*}
|d_z(p^r)|
&\leq \frac{1}{r!} |z|(|z|+1)\cdots(|z|+r-1)\\
&\leq \frac{1}{r!} k(k+1)\cdots(k+r-1)
=d_k(p^r)
\end{align*}
due to $|z| \leq \lfloor{|z|}\rfloor+1 =k$. 
Hence we obtain the result. 
\end{proof}

\begin{lemma}\label{lem:3.6}
Let $\sigma_0$ be a real number. 
Let $k \in \mathbb{Z}_{\geq1}$ and $Y \geq 3$. 
Then there exists an absolute constant $C \geq1$ such that  
\begin{gather*}
\sum_{n=1}^{\infty} d_k(n) n^{-\sigma} e^{-n/Y} 
\ll (Y^{1-\sigma}+1) (C \log{Y})^{k+1}
\end{gather*}
for $\sigma \geq \sigma_0$, where the implied constant depends only on the choice of $\sigma_0$. 
\end{lemma}

\begin{proof}
If $\sigma \geq 1+(2\log{Y})^{-1}$, then we obtain
\begin{gather*}
\sum_{n=1}^{\infty} d_k(n) n^{-\sigma} e^{-n/Y}
\leq \zeta(\sigma)^k
\leq (C \log{Y})^k
\end{gather*}
by formula \eqref{eq:07261516}. 
Thus the result follows in this case. 
Let $\sigma_0 \leq \sigma \leq 1+(2\log{Y})^{-1}$. 
As a simple application of Fubini's theorem, we obtain 
\begin{gather*}
\sum_{n=1}^{\infty} d_k(n) n^{-\sigma} e^{-n/Y}
=\frac{1}{2\pi i} \int_{\RE(w)=c} \zeta(\sigma+w)^k \Gamma(w) Y^w \,dw
\end{gather*}
for any $c>\max(1-\sigma,0)$. 
Taking $c=1-\sigma+(\log{Y})^{-1}$, we have 
\begin{gather}\label{eq:07261546}
\zeta(\sigma+w)^k Y^w
\ll \zeta(\sigma+c)^k Y^{1-\sigma} 
\leq(C \log{Y})^k Y^{1-\sigma}
\end{gather}
on the line $\RE(w)=c$. 
Furthermore, the Stirling formula yields that
\begin{gather}\label{eq:07281548}
\Gamma(w)
\ll |v|^{u-1/2} \exp\left(-\frac{\pi}{2} |v|\right)
\end{gather}
holds for any $w=u+iv \in \mathbb{C}$ with $|v| \geq1$. 
Then the integral of $\Gamma(w)$ on $\RE(w)=c$ is estimated as 
\begin{align*}
\int_{\RE(w)=c} \Gamma(w) \,dw
&\ll \int_{1}^{\infty} v^{1-\sigma_0} e^{-v} \,dv 
+\max_{|v| \leq1} |\Gamma(c+iv)| \\
&\ll_{\sigma_0} 1+\max_{w \in R} |\Gamma(w)|,
\end{align*}
where $R=\{u+iv \in \mathbb{C} \mid (2\log{Y})^{-1} \leq u \leq 2-\sigma_0,~ |v| \leq1\}$. 
Applying the relation $\Gamma(w)=\Gamma(w+1)/w$, we obtain 
\begin{gather*}
\max_{w \in R} |\Gamma(w)| 
\leq \frac{\max_{w \in R} |\Gamma(w+1)|}{\min_{w \in R} |w|}. 
\end{gather*}
Let $R'=\{u+iv \in \mathbb{C} \mid 1 \leq u \leq 3-\sigma_0,~ |v| \leq1\}$. 
Then we have 
\begin{gather*}
\max_{w \in R} |\Gamma(w+1)| 
\leq \max_{w \in R'} |\Gamma(w)| 
\ll_{\sigma_0} 1.
\end{gather*}
In addition, we note that $\min_{w \in R} |w| =(2\log{Y})^{-1}$ holds by definition. 
Hence the upper bound $\max_{w \in R} |\Gamma(w)| \ll_{\sigma_0} \log{Y}$ follows, which deduces 
\begin{gather*}
\int_{\RE(w)=c} \Gamma(w) \,dw 
\ll_{\sigma_0} \log{Y}. 
\end{gather*}
By this and \eqref{eq:07261546}, we obtain 
\begin{gather*}
\sum_{n=1}^{\infty} d_k(n) n^{-\sigma} e^{-n/Y}
\ll_{\sigma_0} (C \log{Y})^k Y^{1-\sigma} \log{Y}
\leq Y^{1-\sigma} (C \log{Y})^{k+1}
\end{gather*}
as desired. 
\end{proof}

Next, we define another arithmetic function $d_z(n,\rho_K)$ for which the formula
\begin{gather}\label{eq:07261755}
L(s, \rho_K)^z
=\sum_{n=1}^{\infty} d_z(n,\rho_K) n^{-s}
\end{gather}
is satisfied. 
Let $A_\mathfrak{a}$ be the matrix defined as \eqref{eq:07261556}. 
We denote by $F(w;\mathfrak{a})$ the function $F(w;\mathfrak{a})=\log{\det(I-w A_\mathfrak{a})}^{-1}$ for $|w|<1$. 
Then we define $H_r(z;\mathfrak{a})$ as the coefficients of the power series
\begin{gather}\label{eq:07272133}
\exp(z F(w;\mathfrak{a}))
=\det(I-w A_\mathfrak{a})^{-z}
=\sum_{r=0}^{\infty} H_r(z;\mathfrak{a}) w^r 
\end{gather}
similarly to $H_r(z)$ of \eqref{eq:07261602}. 
Since $F(w;\mathfrak{a})=F(\alpha_\mathfrak{a} w)+F(\beta_\mathfrak{a} w)$ holds with the eigenvalues $(\alpha_\mathfrak{a}, \beta_\mathfrak{a})$ of the matrix $A_\mathfrak{a}$, we have 
\begin{gather}\label{eq:07251802}
H_r(z;\mathfrak{a})
=\sum_{j=0}^{r} H_j(z) H_{r-j}(z) \alpha_\mathfrak{a}^j \beta_\mathfrak{a}^{r-j}. 
\end{gather}
From the above, we define $d_z(n,\rho_K)$ as the multiplicative function in $n$ determined by $d_z(p^r,\rho_K)=H_r(z;\mathfrak{a})$ if $K$ satisfies a local specification $\mathfrak{a}$ at $p$. 
Then we have formula \eqref{eq:07261755} for any $z \in \mathbb{C}$ and $\RE(s)>1$. 

\begin{lemma}\label{lem:3.7}
Let $z \in \mathbb{C}$ and $k=\lfloor{|z|}\rfloor+1$. 
Then we have $|d_z(n, \rho_K)| \leq d_{2k}(n)$ for any cubic field $K \in L_3^\pm(X)$ and $n \geq1$. 
\end{lemma}

\begin{proof}
We consider the case $n=p^r$ as in the proof of Lemma \ref{lem:3.5}. 
Let $K \in L_3^\pm(X)$ be a cubic field satisfying a local specification $\mathfrak{a}$ at $p$. 
Recall that the eigenvalues $(\alpha_\mathfrak{a}, \beta_\mathfrak{a})$ of the matrix $A_\mathfrak{a}$ satisfy $|\alpha_\mathfrak{a}|, |\beta_\mathfrak{a}| \leq1$. 
Therefore we deduce from Lemma \ref{lem:3.5} and \eqref{eq:07251802} the inequality
\begin{gather}\label{eq:07261813}
|d_z(p^r,\rho_K)|
\leq \sum_{j=0}^{r} d_k(p^j) d_k(p^{r-j}). 
\end{gather}
By \eqref{eq:07261602} and $d_z(p^r)=H_r(z)$, we find that 
\begin{gather*}
\sum_{r=0}^{\infty} d_{2k}(p^r) w^r
= \left(\sum_{r=0}^{\infty} d_k(p^r) w^r\right)^2
= \sum_{r=0}^{\infty} \sum_{j=0}^{r} d_k(p^j) d_k(p^{r-j}) w^r
\end{gather*}
holds. 
Comparing the $r$-th coefficients, we obtain
\begin{gather*}
d_{2k}(p^r)
=\sum_{j=0}^{r} d_k(p^j) d_k(p^{r-j}). 
\end{gather*}
This yields the desired inequality by \eqref{eq:07261813}. 
\end{proof}

\subsection{Results from probability theory}\label{sec:3.3}
In this subsection, we list several lemmas on random variables and probability measures. 

\begin{lemma}\label{lem:3.8}
Let $(X_n)_n$ be a sequence of independent random variables. 
If two series
\begin{gather}\label{eq:07270148}
\sum_{n=1}^{\infty} \left|\mathbb{E}\left[X_n\right]\right|
\quad\text{and}\quad
\sum_{n=1}^{\infty} \mathbb{E}\left[|X_n|^2\right]
\end{gather}
are finite, then $X_1+\cdots+X_n$ converges almost surely as $n \to\infty$. 
\end{lemma}

\begin{proof}
Put $\mu_n=\mathbb{E}\left[X_n\right]$ for $n \geq1$. 
Then the variance $\Var(X_n)=\sigma_n^2$ is given by
\begin{gather*}
\sigma_n^2
=\mathbb{E}\left[|X_n-\mu_n|^2\right]
=\mathbb{E}\left[|X_n|^2\right]-|\mu_n|^2. 
\end{gather*}
Therefore both $\sum_{n=1}^{\infty} \mu_n$ and $\sum_{n=1}^{\infty} \sigma_n^2$ converge by the assumption that \eqref{eq:07270148} are finite. 
As is well known, this yields the convergence of the random variable $X_1+\cdots+X_n$ almost surely; see \cite[Theorem 22.6]{Billingsley1995} for example. 
\end{proof}

\begin{lemma}\label{lem:3.9}
Let $(\mu_n)_n$ be a sequence of probability measures on $(\mathbb{R}, \mathcal{B}(\mathbb{R}))$, and let $\mu$ be a probability measure on $(\mathbb{R}, \mathcal{B}(\mathbb{R}))$. 
Denote by $\phi_n$ and $\phi$ the characteristic functions
\begin{gather*}
\phi_n(\xi)
=\int_{\mathbb{R}} e^{ix \xi} \,d \mu_n(x) 
\quad\text{and}\quad
\phi(\xi)
=\int_{\mathbb{R}} e^{ix \xi} \,d \mu(x). 
\end{gather*}
If $\phi_n(\xi) \to \phi(\xi)$ holds as $n \to\infty$ uniformly in $\xi \in [-R,R]$ for any $R>0$, then the measure $\mu_n$ converges weakly to $\mu$ as $n \to\infty$, that is, $\mu_n(A) \to \mu(A)$ as $n \to\infty$ for any continuity set $A$ of $\mathbb{R}$. 
\end{lemma}

\begin{proof}
This is just L\'{e}vy's continuity theorem; see \cite[Section 3]{JessenWintner1935} for example. 
\end{proof}

\begin{lemma}\label{lem:3.10}
Let $(\mu_n)_n$ be a sequence of probability measures on $(\mathbb{R}, \mathcal{B}(\mathbb{R}))$, and let $\mu$ be a probability measure on $(\mathbb{R}, \mathcal{B}(\mathbb{R}))$.  
Suppose that the convolution measure $\mu_1*\cdots*\mu_n$ converges weakly to $\mu$ as $n \to\infty$. 
Then the support of $\mu$ is represented as
\begin{gather*}
\supp \mu
=\lim_{n \to\infty} (\supp \mu_1+\cdots+\supp \mu_n), 
\end{gather*}
where $\lim_{n \to\infty} (A_1+\cdots+A_n)$ is the set of  all points $\alpha \in \mathbb{R}$ such that $\alpha$ has at least one representation $\alpha=\lim_{n \to\infty} (\alpha_1+\cdots+\alpha_n)$ for some sequence $(\alpha_n)_n$ with $\alpha_n \in A_n$ for each $n \geq1$. 
\end{lemma}

\begin{proof}
We have $\supp (\mu_1*\cdots*\mu_n) = \supp \mu_1+\cdots+\supp \mu_n$. 
Hence the inclusion
\begin{gather*}
\supp \mu
\subset \lim_{n \to\infty} (\supp \mu_1+\cdots+\supp \mu_n)
\end{gather*}
follows from the assumption $\mu_1*\cdots*\mu_n \to \mu$ as $n \to\infty$. 
The opposite inclusion is also proved in \cite[Theorem 3]{JessenWintner1935}. 
\end{proof}

\begin{lemma}\label{lem:3.11}
Let $\mu$ be a probability measure on $(\mathbb{R}, \mathcal{B}(\mathbb{R}))$. 
If the characteristic function $\phi$ satisfies 
\begin{gather*}
\int_{-\infty}^{\infty} |\xi|^k |\phi(\xi)| \,d \xi
<\infty
\end{gather*}
for an integer $k \geq0$. 
Then there exists a non-negative $C^k$-function $D$ on $\mathbb{R}$ such that 
\begin{gather*}
\mu(A)
=\int_{A} D(x) \,\frac{dx}{\sqrt{2\pi}}
\end{gather*}
holds for all $A \in \mathcal{B}(\mathbb{R})$. 
\end{lemma}

\begin{proof}
If $k=0$, then the result follows immediately from L\'{e}vy's inversion formula; see \cite[Theorem 26.2]{Billingsley1995}. 
Moreover, the density function $D$ is given by 
\begin{gather}\label{eq:07270221}
D(x)
=\int_{-\infty}^{\infty} \phi(\xi) e^{-ix \xi} \,\frac{d \xi}{\sqrt{2\pi}}. 
\end{gather}
We obtain the result for $k>0$ by differentiating under the integral in \eqref{eq:07270221}.
\end{proof}

\begin{lemma}\label{lem:3.12}
Let $\mu$ and $\nu$ be probability measures on $(\mathbb{R}, \mathcal{B}(\mathbb{R}))$. 
Denote by $\phi$ and $\psi$ the characteristic function of $\mu$ and $\nu$, respectively. 
Put $F(t)=\mu((-\infty, t])$ and $G(t)=\nu((-\infty, t])$. 
Suppose that there exists a non-negative  continuous function $D$ on $\mathbb{R}$ such that 
\begin{gather}\label{eq:07270233}
\nu(A)
=\int_{A} D(x) \,\frac{dx}{\sqrt{2\pi}}
\end{gather}
holds for all $A \in \mathcal{B}(\mathbb{R})$. 
Then we have 
\begin{gather}\label{eq:07310200}
\sup_{t \in\mathbb{R}} |F(t)-G(t)|
\ll \frac{1}{R} \sup_{x \in\mathbb{R}} D(x)
+\int_{-R}^{R} \left|\frac{\phi(\xi)-\psi(\xi)}{\xi}\right| \,d \xi
\end{gather}
for any $R>0$, where the implied constant is absolute. 
\end{lemma}

\begin{proof}
By equality \eqref{eq:07270233}, the distribution function $G$ is represented as 
\begin{gather*}
G(t)
=\int_{-\infty}^{t} D(x) \,\frac{dx}{\sqrt{2\pi}}. 
\end{gather*}
Therefore $G$ is differentiable, and we see that $G'(t) \ll D(t)$ holds. 
Finally, we apply Esseen's inequality \cite[Theorem 7.16]{Tenenbaum2015} for $F$ and $G$ to obtain the conclusion. 
\end{proof}

\subsection{Other preliminary lemmas}\label{sec:3.4}

\begin{lemma}[Mishou--Nagoshi \cite{MishouNagoshi2006b}]\label{lem:3.13}
Let $H$ be a Hilbert space equipped with the inner product $\langle\,\cdot\,,\cdot\,\rangle$ and the norm $\|\cdot\|$. 
Let $(u_n)_n$ be a sequence in $H$ satisfying
\begin{itemize}
\item[$(\mathrm{a})$] 
$\sum_{n=1}^{\infty} \|u_n\|^2<\infty$ 
\item[$(\mathrm{b})$] 
$\sum_{n=1}^{\infty} |\langle u_n,u \rangle|=\infty$ for any $u \in H$ such that $\|u\|=1$. 
\end{itemize}
Then, for any $v \in H$, $k \geq1$ and $\epsilon>0$, there exist an integer $N=N(v,k,\epsilon) \geq k$ such that the inequality
\begin{gather*}
\left\| v-\sum_{n=k}^{N} c_n u_n \right\|
<\epsilon. 
\end{gather*}
holds with some $c_k,\ldots, c_N \in \{1,-1\}$.
\end{lemma}

Let $M$ be a non-negative function on $\mathbb{R}$ with the Fourier transform $\widetilde{M}$. 
According to Ihara--Matsumoto \cite{IharaMatsumoto2011b}, we say that $M$ is a good density function on $\mathbb{R}$ if both $M$ and $\widetilde{M}$ belong to the class $L^1(\mathbb{R}) \cap L^\infty(\mathbb{R})$, and the conditions
\begin{gather*}
\int_{-\infty}^{\infty} M(x) \,\frac{dx}{\sqrt{2\pi}}
=1 
\quad\text{and}\quad
M(x)
=\int_{-\infty}^{\infty} \widetilde{M}(\xi) e^{-ix \xi} \,\frac{d \xi}{\sqrt{2\pi}}
\end{gather*}
are satisfied. 
Then we have the following result. 

\begin{lemma}[Ihara--Matsumoto \cite{IharaMatsumoto2011b}]\label{lem:3.14}
Let $(X_n)_n$ be a sequence of finite sets with probability measures $\omega_n$. 
We take a function $\ell_n: X_n \to \mathbb{R}$ for each $n \geq1$. 
Let $M$ be a good density function on $\mathbb{R}$. 
Suppose that the condition
\begin{gather}\label{eq:07270304}
\lim_{n \to\infty} \sum_{\chi \in X_n} \omega_n \Phi\left(\ell_n(\chi)\right)
=\int_{-\infty}^{\infty} \Phi(x) M(x) \,\frac{dx}{\sqrt{2\pi}}. 
\end{gather}
is satisfied with the function $\Phi(x)=e^{i \xi x}$ for any $\xi \in \mathbb{R}$, and that the convergence is uniform in $\xi \in [-R,R]$ for any  $R>0$. 
Then we have the following results: 
\begin{itemize}
\item[$(\mathrm{i})$]
\eqref{eq:07270304} holds for any $\Phi \in C_b(\mathbb{R})$; 
\item[$(\mathrm{ii})$]
\eqref{eq:07270304} holds for any $\Phi \in C(\mathbb{R})$ with $|\Phi(x)| \leq \phi_0(|x|)$, where $\phi_0(r)$ is a continuous non-decreasing function on $[0,\infty)$ which satisfies $\phi_0(r)>0$ for any $r \in [0,\infty)$, $\phi_0(r) \to\infty$ as $r \to\infty$, and the following two conditions
\begin{align}
&\sum_{\chi \in X_n} \omega_n \phi_0\left(|\ell_n(\chi)|\right)^2
\ll1, 
\label{eq:07311434}\\
&\int_{-\infty}^{\infty} \phi_0(|x|) M(x) \,\frac{dx}{\sqrt{2\pi}}
<\infty. 
\label{eq:07311435}
\end{align}
\end{itemize}
\end{lemma}

\section{The density functions $C_\sigma$ and $K_\sigma$}\label{sec:4}

\subsection{Random Euler products}\label{sec:4.1}
Let $\mathcal{X}=(\mathcal{X}_p)_p$ and $\mathcal{Y}=(\mathcal{Y}_p)_p$ be the sequences of independent random elements on $\{A_\mathfrak{a} \mid \mathfrak{a} \in \mathcal{A}\}$ as in Section \ref{sec:2}. 
We study several properties on the random Euler products $L(s,\mathcal{X})$ and $L(s,\mathcal{Y})$ defined by infinite products \eqref{eq:07271414}. 

\begin{proposition}\label{prop:4.1}
Let $s=\sigma+it$ be a fixed complex number. 
\begin{itemize}
\item[$(\mathrm{i})$]
If $\sigma>1/2$, then the infinite product on $L(s,\mathcal{X})$ converges almost surely. 
\item[$(\mathrm{ii})$]
If $\sigma>2/3$, then the infinite product on $L(s,\mathcal{Y})$ converges almost surely. 
\end{itemize}
\end{proposition}

\begin{proof}
The local components of infinite products \eqref{eq:07271414} are calculated as 
\begin{align*}
\det\left(I-p^{-s} \mathcal{X}_p\right)^{-1} 
&=1+\tr(\mathcal{X}_p) p^{-s}+O(p^{-2\sigma}), \\
\det\left(I-p^{-s} \mathcal{Y}_p\right)^{-1} 
&=1+\tr(\mathcal{Y}_p) p^{-s}+O(p^{-2\sigma}).
\end{align*}
Recall that $\sum_{p} p^{-2\sigma}$ is convergent for $\sigma>1/2$. 
Then it is sufficient to prove the convergences of $\sum_{p} \tr(\mathcal{X}_p) p^{-\sigma}$ for $\sigma>1/2$ and $\sum_{p} \tr(\mathcal{Y}_p) p^{-\sigma}$ for $\sigma>2/3$ almost surely. 
Let $\sigma>1/2$. 
We see that 
\begin{gather*}
\mathbb{E}\left[\tr(\mathcal{X}_p)\right]
=\sum_{\mathfrak{a} \in \mathcal{A}} C_p(\mathfrak{a}) (\alpha_\mathfrak{a}+\beta_\mathfrak{a})
=\frac{p^{-1}}{1+p^{-1}+p^{-2}}, \\ 
\mathbb{E}\left[|\tr(\mathcal{X}_p)|^2\right]
=\sum_{\mathfrak{a} \in \mathcal{A}} C_p(\mathfrak{a}) (\alpha_\mathfrak{a}+\beta_\mathfrak{a})^2
=\frac{1+p^{-1}}{1+p^{-1}+p^{-2}}, 
\end{gather*}
where $(\alpha_\mathfrak{a}, \beta_\mathfrak{a})$ are the eigenvalues of the matrix $A_\mathfrak{a}$. 
Hence $\mathbb{E}\left[\tr(\mathcal{X}_p)\right] \ll p^{-1}$ and $\mathbb{E}\left[|\tr(\mathcal{X}_p)|^2\right] \ll 1$, which yield that 
\begin{gather*}
\sum_{p} \left|\mathbb{E}\left[\tr(\mathcal{X}_p) p^{-\sigma}\right]\right|
\quad\text{and}\quad
\sum_{p} \mathbb{E}\left[|\tr(\mathcal{X}_p) p^{-\sigma}|^2\right]
\end{gather*}
are finite. 
Therefore $\sum_{p} \tr(\mathcal{X}_p) p^{-\sigma}$ converges almost surely by Lemma \ref{lem:3.8}, and the first result follows. 
Let $\sigma>2/3$. 
We also calculate the expected values for $\mathcal{Y}_p$ as 
\begin{gather*}
\mathbb{E}\left[\tr(\mathcal{Y}_p)\right]
=\frac{1-p^{-1/3}}{(1-p^{-5/3})(1+p^{-1})}
\left(p^{-1/3}+p^{-2/3}+O(p^{-1})\right), \\ 
\mathbb{E}\left[|\tr(\mathcal{Y}_p)|^2\right]
=\frac{1-p^{-1/3}}{(1-p^{-5/3})(1+p^{-1})}
\left(1+p^{-1/3}+p^{-2/3}+O(p^{-1})\right)
\end{gather*}
Thus, we obtain $\mathbb{E}\left[\tr(\mathcal{Y}_p) p^{-\sigma}\right] \ll p^{-\sigma-1/3}$ and $\mathbb{E}\left[|\tr(\mathcal{Y}_p)|^2\right] \ll 1$ in this case. 
Hence we obtain the second result by applying Lemma \ref{lem:3.8} again. 
\end{proof}

By Lemma \ref{prop:4.1}, $L(s,\mathcal{X})$ is a random variable for $\RE(s)>1/2$, and $L(s,\mathcal{Y})$ is a random variable for $\RE(s)>2/3$. 
We also find that 
\begin{align*}
\log{L}(s,\mathcal{X})
&=\sum_{p} \sum_{m=1}^{\infty} \frac{\tr(\mathcal{X}_p^m)}{m} p^{-ms}, \\
\log{L}(s,\mathcal{Y})
&=\sum_{p} \sum_{m=1}^{\infty} \frac{\tr(\mathcal{Y}_p^m)}{m} p^{-ms} 
\end{align*}
are defined for $\RE(s)>1/2$, and for $\RE(s)>2/3$, respectively. 

\begin{proposition}\label{prop:4.2}
Let $s=\sigma+it$ and $z$ be complex numbers. 
\begin{itemize}
\item[$(\mathrm{i})$]
If $\sigma>1/2$, then the expected value
\begin{gather*}
F_s(z)
=\mathbb{E}\left[\exp(z \log{L}(s,\mathcal{X}))\right]
\end{gather*}
is finite and represented as the infinite product
\begin{gather}\label{eq:07271532}
F_s(z)
=\prod_{p} F_{s,p}(z), 
\end{gather}
where $F_{s,p}(z)=\mathbb{E}\left[\det\left(I-p^{-s} \mathcal{X}_p\right)^{-z} \right]$. 
Furthermore, \eqref{eq:07271532} converges uniformly for $\sigma \geq \sigma_0$ and $|z| \leq R$ with any $\sigma_0 >1/2$ and $R>0$.  
\item[$(\mathrm{ii})$]
If $\sigma>2/3$, then the expected value
\begin{gather*}
G_s(z)
=\mathbb{E}\left[\exp(z \log{L}(s,\mathcal{Y}))\right]
\end{gather*}
is finite and represented as the infinite product
\begin{gather}\label{eq:07271537}
G_s(z)
=\prod_{p} G_{s,p}(z), 
\end{gather}
where $G_{s,p}(z)=\mathbb{E}\left[\det\left(I-p^{-s} \mathcal{Y}_p\right)^{-z} \right]$. 
Furthermore, \eqref{eq:07271537} converges uniformly for $\sigma \geq \sigma_0$ and $|z| \leq R$ with any $\sigma_0 >2/3$ and $R>0$.  
\end{itemize}
\end{proposition}

\begin{proof}
For the finiteness of $F_s(z)$, it is sufficient to show that the expected values
\begin{gather*}
\mathbb{E}\left[\exp(r \RE \log{L}(s,\mathcal{X}))\right]
\quad\text{and}\quad
\mathbb{E}\left[\exp(r \IM \log{L}(s,\mathcal{X}))\right]
\end{gather*}
are finite for any $r \in \mathbb{R}$ since we have 
\begin{align*}
|F_s(z)|^2
&=\left|\mathbb{E}\left[ \exp(z \RE \log{L}(s,\mathcal{X})) \exp(iz \IM \log{L}(s,\mathcal{X})) \right]\right|^2 \\
&\leq \mathbb{E}\left[\exp\left(2\RE(z) \RE \log{L}(s,\mathcal{X})\right)\right]
\cdot \mathbb{E}\left[\exp\left(-2\IM(z) \RE \log{L}(s,\mathcal{X})\right)\right]
\end{align*}
by the Cauchy--Schwarz inequality. 
Recall that the random variable 
\begin{gather}\label{eq:07280318}
\log{L}_y(s,\mathcal{X})
=\sum_{p \leq y} \log\det\left(I-p^{-s} \mathcal{X}_p\right)^{-1} 
\end{gather}
converges to $\log{L}(s,\mathcal{X})$ almost surely as $y \to\infty$ for $\RE(s)>1/2$ and $z \in \mathbb{C}$. 
Thus we have 
\begin{align}\label{eq:07271616}
\mathbb{E}\left[\exp(r \RE \log{L}(s,\mathcal{X}))\right]
&\leq \liminf_{y \to\infty}
\mathbb{E}\left[\exp(r \RE \log{L}_y(s,\mathcal{X}))\right] \\
&=\liminf_{y \to\infty} 
\prod_{p \leq y} \mathbb{E}\left[\exp(r \RE \log\det\left(I-p^{-s} \mathcal{X}_p\right)^{-1}\right] \nonumber
\end{align}
by applying Fatou's lemma. 
Here, we also use the independence of $\mathcal{X}=(\mathcal{X}_p)_p$ in the second equality.  
Furthermore, we derive
\begin{gather*}
\exp(r \RE \log\det\left(I-p^{-s} \mathcal{X}_p\right)^{-1}
=1
+r \tr(\mathcal{X}_p) \RE(p^{-s})
+O_r(p^{-2\sigma})
\end{gather*}
by the Taylor series expansion. 
Since $\mathbb{E}\left[\tr(\mathcal{X}_p)\right] \ll p^{-1}$, we find that the formula
\begin{gather}\label{eq:07271656}
\mathbb{E}\left[\exp(r \RE \log\det\left(I-p^{-s} \mathcal{X}_p\right)^{-1}\right]
=1+O_r(p^{-\sigma-1}+p^{-2\sigma})
\end{gather}
holds. 
Note that $\sum_{p} p^{-\sigma-1}$ and $\sum_{p} p^{-2\sigma}$ are convergent for $\sigma>1/2$. 
Thus the infinite product 
\begin{gather*}
\prod_{p} \mathbb{E}\left[\exp(r \RE \log\det\left(I-p^{-s} \mathcal{X}_p\right)^{-1}\right]
\end{gather*}
is convergent as well. 
Hence we conclude that $\mathbb{E}\left[\exp(r \RE \log{L}(s,\mathcal{X}))\right]$ is finite by \eqref{eq:07271616}. 
One can show that $\mathbb{E}\left[\exp(r \IM \log{L}(s,\mathcal{X}))\right]$ is also finite in the same line. 
These results yield the finiteness of $F_s(z)$ as described above. 

Then we proceed to the proof of \eqref{eq:07271532}. 
For $\RE(s)>1/2$ and $z \in \mathbb{C}$, we have the inequality $|\log{L}(s,\mathcal{X})-\log{L}_y(s,\mathcal{X})|\leq1$ almost surely if $y$ is large enough. 
Therefore we derive 
\begin{align*}
\left|\exp(z \log{L}_y(s,\mathcal{X}))\right|
& \leq \exp(|z|) \left|\exp(z \log{L}(s,\mathcal{X}))\right|
\end{align*}
almost surely. 
Furthermore, one can show that $\mathbb{E}\left[\left|\exp(z \log{L}(s,\mathcal{X}))\right|\right]$ is finite in a way similar to the finiteness of $F_s(z)$. 
By the dominated convergence theorem, we have 
\begin{align*}
\mathbb{E}\left[\exp(z \log{L}(s,\mathcal{X}))\right]
&=\lim_{y \to\infty} \mathbb{E}\left[\exp(z \log{L}_y(s,\mathcal{X}))\right] \\
&=\prod_{p} \mathbb{E}\left[\det\left(I-p^{-s} \mathcal{X}_p\right)^{-z} \right], 
\end{align*}
where the second inequality is valid due to the independence of $\mathcal{X}=(\mathcal{X}_p)_p$. 
Hence we obtain \eqref{eq:07271532}. 
Let $\sigma \geq \sigma_0$ and $|z| \leq R$ with $\sigma_0>1/2$ and $R>0$. 
Then we obtain 
\begin{gather*}
F_{s,p}(z)
=1+O(p^{-\sigma_0-1}+p^{-2\sigma_0})
\end{gather*}
similarly to \eqref{eq:07271656}, where the implied constant depends only on $\sigma_0$ and $R$. 
Thus the uniform convergence of \eqref{eq:07271532} follows by noting that $\sum_{p} p^{-\sigma_0-1}$ and $\sum_{p} p^{-2\sigma_0}$ are convergent. 

The results for $G_s(z)$ can be proved similarly to $F_s(z)$. 
The difference is just that the expected value $\mathbb{E}\left[\tr(\mathcal{Y}_p)\right]$ is bounded as $\mathbb{E}\left[\tr(\mathcal{Y}_p)\right] \ll p^{-1/3}$. 
Therefore we omit the proof. 
\end{proof}

\begin{proposition}\label{prop:4.3}
Let $s=\sigma+it$ and $z$ be complex numbers. 
\begin{itemize}
\item[$(\mathrm{i})$]
If $1/2<\sigma<1$, then there exists an absolute constant $c_1>0$ such that 
\begin{gather*}
|F_s(z)|
\leq \exp\left( \frac{c_1}{(2\sigma-1)(1-\sigma)} \frac{(|z|+3)^{1/\sigma}}{\log(|z|+3)} \right).
\end{gather*}
\item[$(\mathrm{i})$]
If $2/3<\sigma<1$, then there exists an absolute constant $c_2>0$ such that 
\begin{gather*}
|G_s(z)|
\leq \exp\left( \frac{c_2}{(3\sigma-2)(1-\sigma)} \frac{(|z|+3)^{1/\sigma}}{\log(|z|+3)} \right).
\end{gather*}
\end{itemize}
\end{proposition}

\begin{proof}
Let $H_r(z)$ and $H_r(z;\mathfrak{a})$ be as in \eqref{eq:07261516} and \eqref{eq:07272133}. 
We have $|H_r(z)| \leq (|z|+1)^r$ by \eqref{eq:07261502}. 
Furthermore, relation	\eqref{eq:07251802} yields $|H_r(z;\mathfrak{a})| \leq (r+1)(|z|+1)^r$. 
Thus we deduce from \eqref{eq:07272133} the inequality
\begin{gather}\label{eq:07272215}
\left|\det(I-w A_\mathfrak{a})^{-z}-1-z\tr(A_\mathfrak{a})w \right|
\leq \sum_{r=2}^{\infty} |H_r(z;\mathfrak{a})| |w|^r
\leq 8 (|z|+1)^2 |w|^2
\end{gather}
for $|w| \leq (2|z|+2)^{-1}$ since $H_0(z;\mathfrak{a})=1$ and $H_1(z;\mathfrak{a})=z\tr(A_\mathfrak{a})$. 
Let $s=\sigma+it$ with $1/2<\sigma<1$ and $z \in \mathbb{C}$. 
If we put $Q=(4|z|+4)^{1/\sigma}$, then $|p^{-s}| \leq (2|z|+2)^{-1}$ is satisfied for prime numbers $p \geq Q$. 
Therefore we derive 
\begin{gather*}
\left|\det(I-p^{-s} A_\mathfrak{a})^{-z}-1-z\tr(A_\mathfrak{a})p^{-s} \right|
\leq 8 (|z|+1)^2 p^{-2\sigma}
\end{gather*}
for any $p \geq Q$ by \eqref{eq:07272215}. 
This deduces the formula
\begin{gather*}
F_{s,p}(z)
=\sum_{\mathfrak{a} \in \mathcal{A}} C_p(\mathfrak{a}) \det(I-p^{-s} A_\mathfrak{a})^{-z}
=1+\mu+E, 
\end{gather*}
where $\mu$ and $E$ satisfy
\begin{align*}
\mu
=\frac{z}{1+p^{-1}+p^{-2}} p^{-s-1}
\quad\text{and}\quad
|E|
\leq 8 (|z|+1)^2 p^{-2\sigma}. 
\end{align*}
For $p \geq Q$, we have $|\mu| \leq 1/8$ and $|E| \leq 1/2$. 
Recall that $\log(1+w) =w+O(|w|^2)$ holds uniformly in $|w| \leq 5/8$. 
We obtain 
\begin{gather}\label{eq:07272311}
\log{F}_{s,p}(z)
=\mu+E+O\left(|\mu|^2+|E|^2 \right)
\ll (|z|+1) p^{-2\sigma},  
\end{gather}
where the implied constant is absolute. 
By the prime number theorem, we obtain that 
\begin{gather}\label{eq:07272306}
\sum_{p \geq Q} \log|F_{s,p}(z)|
\ll \frac{1}{2\sigma-1} \frac{(|z|+3)^{1/\sigma}}{\log(|z|+3)}
\end{gather}
with an absolute implied constant. 
Let $p<Q$. 
We have 
\begin{gather*}
\log\det(I-p^{-s} A_\mathfrak{a})^{-1}
=\log(1-\alpha_\mathfrak{a} p^{-s})^{-1} 
+\log(1-\beta_\mathfrak{a} p^{-s})^{-1}
\ll p^{-\sigma},  
\end{gather*}
where $(\alpha_\mathfrak{a}, \beta_\mathfrak{a})$ are the eigenvalues of $A_\mathfrak{a}$. 
Hence there exists an absolute constant $c>0$ such that 
\begin{gather*}
|F_{s,p}(z)|
\leq \sum_{\mathfrak{a} \in \mathcal{A}} C_p(\mathfrak{a}) \left|\det(I-p^{-s} A_\mathfrak{a})^{-z}\right|
\leq \exp\left( c |z| p^{-\sigma} \right), 
\end{gather*}
is satisfied. 
Then, we again apply the prime number theorem to obtain 
\begin{gather}\label{eq:07272307}
\sum_{p>Q} \log|F_{s,p}(z)|
\ll \frac{1}{1-\sigma} \frac{(|z|+3)^{1/\sigma}}{\log(|z|+3)}. 
\end{gather}
Combining \eqref{eq:07272306} and \eqref{eq:07272307}, we finally arrive at
\begin{gather*}
\prod_{p} |F_{s,p}(z)|
\leq \exp\left( \frac{c_1}{(2\sigma-1)(1-\sigma)} \frac{(|z|+3)^{1/\sigma}}{\log(|z|+3)} \right)
\end{gather*}
with an absolute constant $c_1>0$. 
Hence we obtain the result for $F_s(z)$ by \eqref{eq:07271532}. 

To show the result for $G_s(z)$, we note that $\log{G}_{s,p}(z) \ll (|z|+1) p^{-\sigma-1/3}$ holds for $p \geq Q$ instead of \eqref{eq:07272311}. 
This deduces 
\begin{gather*}
\sum_{p \geq Q} \log|G_{s,p}(z)|
\ll \frac{1}{3\sigma-2} \frac{(|z|+3)^{1/\sigma}}{\log(|z|+3)}
\end{gather*}
by applying the prime number theorem, where the implied constant is absolute.  
The remaining work is similar to the case of $F_s(z)$. 
\end{proof}

\subsection{Construction of the density functions}\label{sec:4.2}
Let $\sigma>1/2$ be a real number. 
We define probability measures $\mu_\sigma$ and $\mu_{\sigma,p}$ as 
\begin{gather}\label{eq:07281421}
\begin{aligned}
\mu_\sigma(A)
&=\mathbb{P}\left(\log{L}(\sigma,\mathcal{X}) \in A\right), \\
\mu_{\sigma, p}(A)
&=\mathbb{P}\left(\log\det\left(I-p^{-\sigma} \mathcal{X}_p\right)^{-1}  \in A\right)
\end{aligned}
\end{gather}
for $A \in \mathcal{B}(\mathbb{R})$, where $p$ is a prime number. 
In a similar way, we define for $\sigma>2/3$ probability measures $\nu_\sigma$ and $\nu_{\sigma,p}$ as 
\begin{align*}
\nu_\sigma(A)
&=\mathbb{P}\left(\log{L}(\sigma,\mathcal{Y}) \in A\right), \\
\nu_{\sigma, p}(A)
&=\mathbb{P}\left(\log\det\left(I-p^{-\sigma} \mathcal{Y}_p\right)^{-1} \in A\right). 
\end{align*}
Denote by $\phi_\sigma$ and $\psi_\sigma$ be the characteristic functions of $\mu_\sigma$ and $\nu_\sigma$, respectively. 

\begin{proposition}\label{prop:4.4}
\begin{itemize}
\item[$(\mathrm{i})$]
Let $\sigma>1/2$ be a real number. 
There exists a positive constant $c_1(\sigma)$ depending only on $\sigma$ such that 
\begin{gather*}
|\phi_\sigma(\xi)|
\leq \exp\left(-c_1(\sigma) \frac{|\xi|^{1/\sigma}}{\log|\xi|}\right)
\end{gather*}
holds for any $\xi \in \mathbb{R}$ with $|\xi| \geq3$. 
\item[$(\mathrm{i})$]
Let $\sigma>2/3$ be a real number. 
There exists a positive constant $c_2(\sigma)$ depending only on $\sigma$ such that 
\begin{gather*}
|\psi_\sigma(\xi)|
\leq \exp\left(-c_2(\sigma) \frac{|\xi|^{1/\sigma}}{\log|\xi|}\right)
\end{gather*}
holds for any $\xi \in \mathbb{R}$ with $|\xi| \geq3$. 
\end{itemize}
\end{proposition}

\begin{proof}
Recall that $H_0(z;\mathfrak{a})=1$, $H_1(z;\mathfrak{a})=z\tr(A_\mathfrak{a})$ and  
\begin{gather*}
H_2(z;\mathfrak{a})
=z\tr(A_\mathfrak{a})+\frac{z^2}{2}\tr(A_\mathfrak{a})^2
\end{gather*}
for the coefficients $H_r(z;\mathfrak{a})$ as in \eqref{eq:07272133}. 
Then, similarly to \eqref{eq:07272215}, we obtain
\begin{align}\label{eq:07280155}
&\left|\det(I-w A_\mathfrak{a})^{-z}
-1
-z\tr(A_\mathfrak{a})w
-z\tr(A_\mathfrak{a})w^2-\frac{z^2}{2}\tr(A_\mathfrak{a})^2w^2\right| \\
&\leq \sum_{r=3}^{\infty} |H_r(z;\mathfrak{a})| |w|^r \nonumber
\leq 10 (|z|+1)^3 |w|^3
\end{align}
for $|w| \leq (2|z|+2)^{-1}$. 
Let $\sigma>1/2$. 
The characteristic function $\phi_\sigma$ is represented as
\begin{gather*}
\phi_\sigma(\xi)
=\mathbb{E}\left[\exp(i \xi \log{L}(\sigma,\mathcal{X}))\right]
=F_\sigma(i \xi)
\end{gather*}
by using the function $F_s(z)$ of Proposition \ref{prop:4.2}. 
Therefore, formula \eqref{eq:07271532} yields 
\begin{gather}\label{eq:07280236}
\phi_\sigma(\xi)
=\prod_{p} F_{\sigma,p}(i \xi),  
\end{gather}
where $F_{\sigma,p}(i \xi)$ is represented as 
\begin{align*}
F_{\sigma,p}(i \xi)
=\sum_{\mathfrak{a} \in \mathcal{A}} C_p(\mathfrak{a}) \det(I-p^{-\sigma} A_\mathfrak{a})^{-i \xi}. 
\end{align*}
Put $Q=(4|\xi|+4)^{1/\sigma}$ as in the proof of Proposition \ref{prop:4.2}. 
Since $p^{-\sigma} \leq (2|z|+2)^{-1}$ holds for $p \geq Q$, we deduce from \eqref{eq:07280155} that 
\begin{gather}\label{eq:07280206}
F_{\sigma,p}(i \xi)
=1+i \mu_1+i \mu_{2,1}-\frac{1}{2} \mu_{2,2}+E
\end{gather}
for $p \geq Q$, where $\mu_1$, $\mu_{2,1}$ and $\mu_{2,2}$ are real numbers given by
\begin{align*}
\mu_1
&=\frac{\xi}{1+p^{-1}+p^{-2}} p^{-\sigma-1},\\
\mu_{2,1}
&=\frac{\xi}{1+p^{-1}+p^{-2}} (1+p^{-1}) p^{-2\sigma},\\
\mu_{2,2}
&=\frac{\xi^2}{1+p^{-1}+p^{-2}} (1+p^{-1}) p^{-2\sigma},
\end{align*}
and $E$ is evaluated as $|E| \leq 10(|\xi|+1)^3 p^{-3\sigma}$. 
For $p \geq Q$, we have $|\mu_1| \leq1/8$, $|\mu_{2,1}| \leq1/16$, $|\mu_{2,2}| \leq1/16$ and $|E| \leq5/32$. 
Hence we deduce from \eqref{eq:07280206} the asymptotic formula 
\begin{gather*}
\log{F}_{\sigma,p}(i \xi)
=i \mu_1+i \mu_{2,1}-\frac{1}{2} \mu_{2,2}
+O\left(|\mu_1|^2+|\mu_{2,1}|^2+|\mu_{2,1}|^2+|E|^2\right) 
\end{gather*}
for $p \geq Q$, where the implied constant is absolute. 
Since $\mu_1$, $\mu_{2,1}$ and $\mu_{2,2}$ are real, we further obtain
\begin{align*}
\log|F_{\sigma,p}(i \xi)|
&=-\frac{1}{2} \mu_{2,2}
+O\left(|\mu_1|^2+|\mu_{2,1}|^2+|\mu_{2,1}|^2+|E|^2\right) \\
&=-\frac{\xi^2}{2} \frac{1+p^{-1}}{1+p^{-1}+p^{-2}} p^{-2\sigma}
+O\left(\xi^2 p^{-2\sigma-2}+|\xi|^3 p^{-3\sigma}\right). 
\end{align*}
Thus, there exists an absolute constant $A>0$ such that 
\begin{gather}\label{eq:07280225}
\log|F_{\sigma,p}(i \xi)|
\leq-\frac{\xi^2}{4} p^{-2\sigma} 
+A \xi^2 p^{-2\sigma-2}
+A |\xi|^3 p^{-3\sigma}
\end{gather}
holds for $p \geq Q$. 
Put $Q(M)=(M |\xi|)^{1/\sigma}$ with $M \geq1$. 
Then we have $Q(M) \geq (3M)^{1/\sigma}$ for $|\xi| \geq3$. 
We take a constant $M_\sigma \geq 6$ depending only on $\sigma$ so that 
\begin{gather*}
Q(M_\sigma)^{-2} 
\leq \frac{1}{16A} 
\quad\text{and}\quad 
A |\xi| Q(M_\sigma)^{-\sigma}
\leq \frac{1}{16A}
\end{gather*}
are satisfied.
We have also $Q(M_\sigma) \geq Q$ due to $M_\sigma \geq6$. 
Therefore, inequality \eqref{eq:07280225} yields 
\begin{gather*}
\log|F_{\sigma,p}(i \xi)|
\leq-\frac{\xi^2}{8} p^{-2\sigma} 
\end{gather*}
for $p \geq Q(M_\sigma)$. 
By the prime number theorem, we obtain 
\begin{gather*}
\sum_{p \geq Q(M_\sigma)} \log|F_{\sigma,p}(i \xi)|
\leq -\frac{\xi^2}{8} \sum_{p \geq Q(M_\sigma)} p^{-2\sigma}
\leq -c_1(\sigma) \frac{|\xi|^{1/\sigma}}{\log|\xi|},
\end{gather*}
where $c_1(\sigma)$ is a positive constant depending only on $\sigma$. 
Thus, the inequality
\begin{gather*}
\prod_{p \geq Q(M_\sigma)} |F_{\sigma,p}(i \xi)|
\leq\exp\left(-c_1(\sigma) \frac{|\xi|^{1/\sigma}}{\log|\xi|}\right) 
\end{gather*}
follows. 
Note that $|F_{\sigma,p}(i \xi)|\leq1$ holds for any prime number $p$. 
Hence we obtain the desired inequality for $|\phi_\sigma(\xi)|$ by formula \eqref{eq:07280236}. 
The result for $\psi_\sigma$ can be shown in the same line. 
\end{proof}

Let $\sigma>1/2$ be a real number. 
For all $k \geq0$, we obtain 
\begin{gather*}
\int_{-\infty}^{\infty} |\xi|^k |\phi(\xi)| \,d \xi
<\infty
\end{gather*}
by Proposition \ref{prop:4.4} $(\mathrm{i})$. 
Hence we deduce from Lemma \ref{lem:3.11} the existence of a non-negative $C^\infty$-function $C_\sigma$ on $\mathbb{R}$ such that 
\begin{gather*}
\mu_\sigma(A)
=\int_{A} C_\sigma(x) \,\frac{dx}{\sqrt{2\pi}}
\end{gather*}
holds for all $A \in \mathcal{B}(\mathbb{R})$.
The function $C_\sigma$ is given by
\begin{gather*}
C_\sigma(x)
=\int_{-\infty}^{\infty} \phi_\sigma(\xi) e^{-ix \xi} \,\frac{d \xi}{\sqrt{2\pi}} 
\end{gather*}
by formula \eqref{eq:07270221}. 
Furthermore, it can be constructed as an infinite convolution of Schwartz distributions as follows. 
Let $C_{\sigma,p}$ denote the Schwartz distribution on $\mathbb{R}$ such that
\begin{align*}
\int_{\mathbb{R}} \Phi(x) C_{\sigma,p}(x) \,\frac{dx}{\sqrt{2\pi}}
=\mathbb{E}\left[\Phi\left(\log\det\left(I-p^{-\sigma} \mathcal{X}_p\right)^{-1}\right)\right]  
\end{align*} 
for each prime number $p$. 
In other words, we put 
\begin{gather*}
C_{\sigma,p}(x)
=\sqrt{2\pi} \sum_{\mathfrak{a} \in \mathcal{A}} C_p(\mathfrak{a}) 
\delta\left(x-\log\det\left(I-p^{-\sigma} A_\mathfrak{a}\right)^{-1}\right), 
\end{gather*}
where $\delta$ is the Dirac distribution on $\mathbb{R}$. 
Denote $n$-th prime number by $p_n$. 
Then we obtain the identity
\begin{gather*}
\int_{\mathbb{R}} \Phi(x) (C_{\sigma,p_1}*\cdots*C_{\sigma,p_n})(x) \,\frac{dx}{\sqrt{2\pi}}
=\mathbb{E}\left[\Phi\left(\log{L}_y(\sigma,\mathcal{X})\right)\right]  
\end{gather*}
with $y=p_n$, where $\log{L}_y(\sigma,\mathcal{X})$ is the random variable defined by \eqref{eq:07280318}. 
Letting $n \to\infty$, we have
\begin{align*}
\int_{\mathbb{R}} \Phi(x) \lim_{n \to\infty} (C_{\sigma,p_1}*\cdots*C_{\sigma,p_n})(x) \,\frac{dx}{\sqrt{2\pi}}
&=\mathbb{E}\left[\Phi\left(\log{L}(\sigma,\mathcal{X})\right)\right] \\
&=\int_{\mathbb{R}} \Phi(x) C_\sigma(x) \,\frac{dx}{\sqrt{2\pi}}. 
\end{align*}
Hence we conclude that the density function $C_\sigma$ is given by 
\begin{gather*}
C_\sigma(x)
=\ast_p\, C_{\sigma,p}(x), 
\end{gather*}
where $\ast_p$ stands for the infinite convolution over the all prime numbers. 
We also obtain the density function $K_\sigma$ for $\sigma>2/3$ such that 
\begin{gather*}
\nu_\sigma(A)
=\int_{A} K_\sigma(x) \,\frac{dx}{\sqrt{2\pi}}
\end{gather*}
holds by Lemma \ref{lem:3.11} and Proposition \ref{prop:4.4} $(\mathrm{ii})$. 
It is a non-negative $C^\infty$-function, and we have the infinite convolution representation
\begin{gather*}
K_\sigma(x)
=\ast_p\, K_{\sigma,p}(x), 
\end{gather*}
where $K_{\sigma,p}$ is the Schwartz distribution on $\mathbb{R}$ defined as
\begin{gather*}
K_{\sigma,p}(x)
=\sqrt{2\pi} \sum_{\mathfrak{a} \in \mathcal{A}} K_p(\mathfrak{a}) 
\delta\left(x-\log\det\left(I-p^{-\sigma} A_\mathfrak{a}\right)^{-1}\right). 
\end{gather*}
We also remark that the $M$-function $M_\sigma$ of Theorem \ref{thm:1.3} is constructed as a similar infinite convolution of Schwartz distributions; see \cite[Section 3]{IharaMatsumoto2011a}.

\subsection{Analytic properties of the density functions}\label{sec:4.3}
In the remaining part of Section \ref{sec:4} is devoted to the proof of Theorem \ref{thm:2.1}. 
We prove only the properties of the density function $C_\sigma$ since one can obtain the corresponding properties of $K_\sigma$ by following similar arguments.

\subsubsection{Proof of Theorem \ref{thm:2.1} $(\mathrm{i})$}\label{sec:4.3.1}
Let $\sigma>1/2$ be a real number. 
Denote by $\mu_\sigma$ and $\mu_{\sigma,p}$ the probability measures defined as \eqref{eq:07281421}. 
Since $\mathcal{X}=(\mathcal{X}_p)_p$ is independent, we have 
\begin{gather*}
\mu_{\sigma,p_1}*\cdots*\mu_{\sigma,p_n}(A)
=\mathbb{P}\left(\log{L}_y(\sigma,\mathcal{X}) \in A\right)
\end{gather*}
with $y=p_n$, where $p_n$ is the $n$-th prime number, and $\log{L}_y(\sigma,\mathcal{X})$ is defined as \eqref{eq:07280318}. 
Therefore we see that $\mu_{\sigma,p_1}*\cdots*\mu_{\sigma,p_n} \to \mu_\sigma$ as $n \to\infty$ by Proposition \ref{prop:4.1}. 
We deduce from Lemma \ref{lem:3.10} the identity
\begin{gather}\label{eq:07281429}
\supp \mu_\sigma
=\lim_{n \to\infty} (\supp \mu_{\sigma,p_1}+\cdots+\supp \mu_{\sigma,p_n}). 
\end{gather}
Note that $\supp C_\sigma$ is equal to $\supp \mu_\sigma$. 
Thus Theorem \ref{thm:2.1} $(\mathrm{i})$ follows if the right-hand side of \eqref{eq:07281429} is $\mathbb{R}$. 
To prove this, we apply Lemma \ref{lem:3.13} for $H=\mathbb{R}$ equipped with the usual inner product $\langle{u,v}\rangle=uv$. 
Let $u_n=p_n^{-\sigma}$ with $1/2<\sigma \leq1$. 
Then conditions $(\mathrm{a})$ and $(\mathrm{b})$ in Lemma \ref{lem:3.13} are easily checked. 
Take a real number $\epsilon>0$ arbitrarily. 
Note that we have 
\begin{gather*}
\sum_{n=k}^{N} \sum_{m=2}^{\infty} \frac{2}{m} p_n^{-m \sigma}
\leq B_\sigma k^{1-2\sigma}
\end{gather*}
for any $N \geq k \geq1$, where $B_\sigma>0$ is a constant depending only on $\sigma$. 
Then we take an integer $k=k(\sigma,\epsilon) \geq 2$ so that $B_\sigma k^{1-2\sigma}<\epsilon$ is satisfied. 
Let $x$ be an arbitrary real number. 
By Lemma \ref{lem:3.13}, we obtain an integer $N=N(\sigma,\epsilon,x) \geq k$ and numbers $c_k, \ldots, c_N \in\{1,-1\}$ such that 
\begin{gather}\label{eq:07281442}
\left|\left(x-\sum_{n=1}^{k-1} \sum_{m=1}^\infty \frac{2}{m} p_n^{-m\sigma}\right)
-\sum_{n=k}^{N} c_n p_n^{-\sigma}\right|
<\epsilon
\end{gather}
is satisfied. 
For $1\leq n \leq N$, we choose a symbol $\mathfrak{a}_n \in \mathcal{A}$ as follows. 
For $1 \leq n<k$, take $\mathfrak{a}_n=(111)$. 
For $k \leq n \leq N$, take $\mathfrak{a}_n=(1^21)$ if $c_n=1$ and $\mathfrak{a}_n=(3)$ if $c_n=-1$. 
Then \eqref{eq:07281442} deduces 
\begin{gather*}
\left|x-\sum_{n=1}^{k-1} \sum_{m=1}^{\infty} \frac{\tr(A_{\mathfrak{a}_n}^m)}{m} p_n^{-m\sigma}
-\sum_{n=k}^{N} \tr(A_{\mathfrak{a}_n}) p_n^{-\sigma}\right|
<\epsilon. 
\end{gather*}
By the choice of $k=k(\sigma,\epsilon)$, we finally obtain 
\begin{gather*}
\left|x-\sum_{n=1}^{N} \log\det\left(I-p_n^{-\sigma} A_{\mathfrak{a}_n}\right)^{-1}\right|
<2\epsilon, 
\end{gather*}
which means that $x$ belongs to $\lim_{N \to\infty} (\supp \mu_{\sigma,p_1}+\cdots+\supp \mu_{\sigma,p_N})$. 
Hence the desired result follows.

\subsubsection{Proof of Theorem \ref{thm:2.1} $(\mathrm{ii})$}\label{sec:4.3.2}
Let $\sigma>1$ be a real number. 
By \eqref{eq:07281429}, it is sufficient to show that 
\begin{gather}\label{eq:07281457}
\lim_{n \to\infty} (\supp \mu_{\sigma,p_1}+\cdots+\supp \mu_{\sigma,p_n})
\subset [-R,R]
\end{gather}
with some $R>0$. 
We have 
\begin{gather*}
\supp \mu_{\sigma,p}
=\left\{\log\det\left(I-p^{-\sigma} A_{\mathfrak{a}}\right) ~\middle|~ \mathfrak{a} \in \mathcal{A}\right\}
\end{gather*}
for every prime number $p$. 
Since $\sigma>1$, we have
\begin{gather*}
\left|\sum_{p} \log\det\left(I-p^{-\sigma} A_{\mathfrak{a}_p}\right) \right|
\leq \sum_{p} \sum_{m=1}^{\infty} \frac{|\tr(A_{\mathfrak{a}_p}^m)|}{m} p^{-m\sigma}
\leq 2\log\zeta(\sigma)
\end{gather*}
for any choices of $\mathfrak{a}_p \in \mathcal{A}$.
Hence we obtain \eqref{eq:07281457} with $R=2\log{\zeta(\sigma)}$.

\subsubsection{Proof of Theorem \ref{thm:2.1} $(\mathrm{iii})$}\label{sec:4.3.3}
Let $\sigma>1/2$ and $a>0$. 
We see that the equality
\begin{gather*}
\int_{-\infty}^{\infty} e^{ax} C_\sigma(x) \,\frac{dx}{\sqrt{2\pi}}
=\mathbb{E}\left[\exp(a \log{L}(\sigma,\mathcal{X}))\right]
\end{gather*}
is valid. 
Hence Theorem \ref{thm:2.1} $(\mathrm{iii})$ is just a consequence of Proposition \ref{prop:4.2}.

\section{Complex moments of $L(\sigma,\rho_K)$}\label{sec:5}
The goal of this section is to complete the proofs of Theorems \ref{thm:2.2} and \ref{thm:2.3} and their corollaries. 
We fix the branch of $\log{L}(s,\rho_K)$ as in Section \ref{sec:3.1}. 
Then we define 
\begin{gather*}
g_z(s,\rho_K)
=L(s,\rho_K)^z
\end{gather*}
for $z \in \mathbb{C}$, which is holomorphic for $s \in G_K$. 
If $\RE(s)>1$, then we obtain 
\begin{gather}\label{eq:07281545}
g_z(s,\rho_K)
=\sum_{n=1}^{\infty} d_z(n,\rho_K) n^{-s}
\end{gather}
by \eqref{eq:07261755}, where $d_z(n,\rho_K)$ is the multiplicative function of Section \ref{sec:3.2}. 
The function $g_z(s,\rho_K)$ is an analogue of the $g$-function introduced by Ihara--Matsumoto \cite{IharaMatsumoto2011b}.

\subsection{Approximation of the $g$-function}\label{sec:5.1}
Throughout this subsection, we write 
\begin{gather*}
c
=\max\{1-\sigma,0\}+\frac{1}{\log{X}} 
\quad\text{and}\quad
\kappa
=\frac{1}{\log\log{X}} 
\end{gather*}
for $\sigma>1/2$ and $X \geq3$. 
We define the functions $g_z^+(\sigma,\rho_K;Y)$ and $g_z^-(\sigma,\rho_K;Y)$ as 
\begin{gather*}
g_z^\pm(\sigma,\rho_K;Y)
=\frac{1}{2\pi i} \int_{L^\pm} g_z(\sigma+w,\rho_K) \Gamma(w) Y^w \,dw
\end{gather*}
for $Y \geq1$, where the integral contours $L^\pm$ are described as follows. 
Let $L^+$ be the vertical line $\RE(w)=c$. 
Let $L^-=L_1+\cdots+L_5$ be the contour given by connecting the points $c-i\infty$, $c-i(\log{X})^2$, $-\kappa-i(\log{X})^2$, $-\kappa+i(\log{X})^2$, $c+i(\log{X})^2$, $c+i\infty$, in order. 
The connection between $g_z(\sigma,\rho_K)$ and $g_z^\pm(\sigma,\rho_K;Y)$ is as follows.

\begin{proposition}\label{prop:5.1}
Let $\sigma_1$ be a real number for which \eqref{eq:07251505} holds, and suppose that $\sigma>\sigma_1+\kappa$ is satisfied. 
Take a cubic field $K \in L_3^\pm(X) \setminus E(X)$, where $E(X)$ is defined as \eqref{eq:07252237}. 
Then the equality 
\begin{gather*}
g_z(\sigma,\rho_K)
=g_z^+(\sigma,\rho_K;Y)-g_z^-(\sigma,\rho_K;Y)
\end{gather*}
holds for any $z \in \mathbb{C}$ and $Y \geq1$. 
\end{proposition}

\begin{proof}
We have $\RE(\sigma+w)>1$ for $w \in L^+$ by the choice of $c$. 
Hence formula \eqref{eq:07281545} yields that $|g_z(\sigma+w,\rho_K)|$ is bounded for $w \in L^+$. 
Furthermore, $|\Gamma(w)|$ is rapidly decreasing on $L^+$ as $|\IM(w)| \to\infty$ by \eqref{eq:07281548}. 
As a result, $g_z(\sigma+w,\rho_K) \Gamma(w) Y^w$ is absolutely integrable on the contour $L^+$. 
Then we shift the contour to $L^-$. 
Remark that the function $g_z(s,\rho_K)$ is holomorphic in the region $\RE(s)>\sigma_1$, $|\IM(s)|<(\log{X})^3$ if $K \in L_3^\pm(X) \setminus E(X)$. 
Hence we do not encounter any pole of the integrand except for the simple pole at $w=0$ while shifting the contour. 
The residue at $w=0$ is equal to $g_z(\sigma,\rho_K)$. 
Therefore, we obtain 
\begin{align*}
&\frac{1}{2\pi i} \int_{L^+} g_z(\sigma+w,\rho_K) \Gamma(w) Y^w \,dw \\
&=\frac{1}{2\pi i} \int_{L^-} g_z(\sigma+w,\rho_K) \Gamma(w) Y^w \,dw
+g_z(\sigma,\rho_K)
\end{align*}
as desired. 
\end{proof}

Then, we show several properties on the functions $g_z^+(\sigma,\rho_K;Y)$ and $g_z^-(\sigma,\rho_K;Y)$. 

\begin{proposition}\label{prop:5.2}
Let $\sigma>1/2$ be a real number. 
Take a cubic field $K \in L_3^\pm(X)$ arbitrarily. 
Then the function $g_z^+(\sigma,\rho_K;Y)$ is represented as 
\begin{gather*}
g_z^+(\sigma,\rho_K;Y)
=\sum_{n=1}^{\infty} d_z(n,\rho_K) n^{-\sigma} e^{-n/Y}
\end{gather*}
for any $z \in \mathbb{C}$ and $Y \geq3$. 
If $Y=X^\eta$ with $\eta>0$, then we have also 
\begin{gather}\label{eq:07281608}
g_z^+(\sigma,\rho_K;Y)
\ll X^\eta
\end{gather}
for any $z \in \mathbb{C}$ such that $|z| \leq R_\sigma(X)$, where $R_\sigma(X)$ is given by \eqref{eq:07281607}. 
The implied constant in \eqref{eq:07281608} depends only on $\sigma$ and $\eta$. 
\end{proposition}

\begin{proof}
Recall that we have $\RE(\sigma+w)>1$ for $w \in L^+$. 
Hence formula \eqref{eq:07281545} yields 
\begin{align*}
g_z^+(\sigma,\rho_K;Y)
&=\frac{1}{2\pi i} \int_{\RE(w)=c} g_z(\sigma+w,\rho_K) \Gamma(w) Y^w \,dw \\
&=\sum_{n=1}^{\infty} d_z(n,\rho_K) n^{-\sigma} 
\frac{1}{2\pi i} \int_{\RE(w)=c} \Gamma(w) (n/Y)^{-w} \,dw \\
&=\sum_{n=1}^{\infty} d_z(n,\rho_K) n^{-\sigma} e^{-n/Y}. 
\end{align*}
Put $k=\lfloor{|z|}\rfloor+1$. 
By Lemmas \ref{lem:3.6} and \ref{lem:3.7}, we further obtain
\begin{gather*}
|g_z^+(\sigma,\rho_K;Y)|
\leq \sum_{n=1}^{\infty} d_{2k}(n) n^{-\sigma} e^{-n/Y}
\ll (Y^{1-\sigma}+1) (C \log{Y})^{2k+1}, 
\end{gather*}
where the implied constant is absolute. 
Hence \eqref{eq:07281608} follows if we let $Y=X^\eta$ with $\eta>0$ and $|z| \leq R_\sigma(X)$. 
\end{proof}

\begin{proposition}\label{prop:5.3}
Let $\sigma_1$ be a real number for which \eqref{eq:07251505} holds, and suppose that $\sigma \geq \sigma_1+3\kappa$ is satisfied. 
Take a cubic field $K \in L_3^\pm(X) \setminus E(X)$, where $E(X)$ is defined as \eqref{eq:07252237}. 
If $Y=X^\eta$ with $\eta>0$, then there exists a constant $b(\eta)>0$ depending only on $\eta$ such that 
\begin{gather*}
g_z^-(\sigma,\rho_K;Y)
\ll \exp\left(-\frac{\eta}{2} \frac{\log{X}}{\log\log{X}}\right)
\end{gather*}
holds for any $z \in \mathbb{C}$ satisfying $|z| \leq b(\eta) R_\sigma(X)$, where $R_\sigma(X)$ is given by \eqref{eq:07281607}, and the implied constant is absolute. 
\end{proposition}

\begin{proof}
We divide the integral contour $L^-$ into $L_1,L_2,\ldots,L_5$ as above. 
Then we have $\RE(\sigma+w) \geq 1+(\log{X})^{-1}$ on $L_1$ and $L_5$. 
Hence formula \eqref{eq:07252233} is available to deduce
\begin{gather*}
\log{L}(\sigma+w,\rho_K) 
\ll \log\log{X}
\end{gather*}
for $w \in L_1 \cup L_5$. 
If $w$ lies on $L_2 \cup L_3 \cup L_4$, then we have $\RE(\sigma+w) \geq \sigma_1+2\kappa$ and $|\IM(\sigma+w)| \leq (\log{X})^2$. 
Since $K \in L_3^\pm(X) \setminus E(X)$, we have
\begin{gather*}
\log{L(\sigma+w,\rho_K)}
\ll (\log\log{X}) (\log{X})^{(1-\sigma)/(1-\sigma_1)} 
+\log\log{X}
\end{gather*}
by Lemma \ref{lem:3.3} in this case. 
As a result, we obtain for any $w \in L^-$ the upper bound
\begin{gather}\label{eq:07281702}
g_z(\sigma+w,\rho_K)
\ll \exp\left(A b(\eta) \frac{\log{X}}{\log\log{X}}\right)
\end{gather}
due to $|z| \leq b(\eta) R_\sigma(X)$, where $A>0$ is an absolute constant. 
Note that $Y^w$ satisfies 
\begin{gather}\label{eq:07281703}
Y^w
\ll
\begin{cases}
\displaystyle{X^\eta} 
& \text{if $w$ lies on $L_1$, $L_2$, $L_4$, $L_5$}, \\
\displaystyle{\exp\left(-\eta \frac{\log{X}}{\log\log{X}}\right)} 
& \text{if $w$ lies on $L_3$}. 
\end{cases}
\end{gather}
Applying \eqref{eq:07281548}, we also obtain
\begin{gather}\label{eq:07281704}
\frac{1}{2\pi i} \int_{L_j} \Gamma(w) \,dw
\ll
\begin{cases}
\exp\left(-(\log{X})^2\right)
& \text{if $j=1,2,4,5$}, \\
\log\log{X}
& \text{if $j=3$}. 
\end{cases}
\end{gather}
Combining \eqref{eq:07281702}, \eqref{eq:07281703} and \eqref{eq:07281704}, we obtain the desired estimate on $g_z^-(\sigma,\rho_K;Y)$ if we choose $b(\eta)>0$ as a suitably small constant. 
\end{proof}

As an application of Propositions \ref{prop:5.1}, \ref{prop:5.2} and \ref{prop:5.3}, we obtain the following asymptotic formula for the $g$-function $g_z(\sigma,\rho_K)$. 

\begin{corollary}\label{cor:5.1}
Let $\sigma_1$ be a real number for which \eqref{eq:07251505} holds. 
Take a cubic field $K \in L_3^\pm(X) \setminus E(X)$, where $E(X)$ is defined as \eqref{eq:07252237}. 
If $Y=X^\eta$ with $\eta>0$, then there exists a constant $0<b(\eta) \leq1$ depending only on $\eta$ such that 
\begin{gather*}
g_z(\sigma,\rho_K)
=\sum_{n=1}^{\infty} d_z(n,\rho_K)n^{-\sigma}e^{-n/Y}
+O\left(\exp\left(-\frac{\eta}{2}\frac{\log{X}}{\log\log{X}}\right)\right)
\end{gather*}
holds for $\sigma>\sigma_1$ with $z \in \mathbb{C}$ satisfying $|z| \leq b(\eta) R_\sigma(X)$, where $R_\sigma(X)$ is given by \eqref{eq:07281607}, and the implied constant depends only on $\sigma$ and $\eta$. 
\end{corollary}

One can improve the upper bound on $g_z^-(\sigma,\rho_K;Y)$ of Proposition \ref{prop:5.3} if GRH is assumed. 
Indeed, we obtain the following result. 

\begin{proposition}\label{prop:5.4}
Assume GRH. 
Let $\sigma>1/2$ be a real number. 
Choose a positive real number $\kappa_1(\sigma)$ depending on $\sigma$ such that $\sigma-\kappa_1(\sigma)>1/2$ is satisfied. 
Take a cubic field $K \in L_3^\pm(X)$ arbitrarily. 
If $Y=X^\eta$ with $\eta>0$, then there exists an absolute constant $A \geq1$ such that 
\begin{gather*}
g_z^-(\sigma,\rho_K;Y)
\ll X^{-\eta \kappa_1(\sigma)} \exp\left(A \frac{\log{X}}{\log\log{X}}\right)
\end{gather*}
holds for any $z \in \mathbb{C}$ satisfying $|z| \leq \tilde{R}_\sigma(X)$, where $\tilde{R}_\sigma(X)$ is given by \eqref{eq:07301552} with $a(\sigma)=\min(2\sigma-2\kappa_1(\sigma)-1, 2/3)>0$. 
Here, the implied constant depends only on the choice of $\kappa_1(\sigma)$. 
\end{proposition}

\begin{proof}
Denote by $\tilde{L}$ the contour $\tilde{L}_1+\cdots+\tilde{L}_5$ which is obtained by connecting the points $c-i \infty$, $c-i(\log{X})^2$, $-\kappa_1(\sigma)-i(\log{X})^2$, $-\kappa_1(\sigma)+i(\log{X})^2$, $c+i(\log{X})^2$, $c+i \infty$, in order. 
Then we have 
\begin{align*}
g_z^-(\sigma,\rho_K;Y)
&=\frac{1}{2\pi i} \int_{L^-} g_z(\sigma+w,\rho_K) \Gamma(w) Y^w \,dw \\
&=\frac{1}{2\pi i} \int_{\tilde{L}} g_z(\sigma+w,\rho_K) \Gamma(w) Y^w \,dw
\end{align*}
since GRH yields the absence of any pole of $g_z(\sigma+w,\rho_K)$ while shifting the contour. 
We have 
\begin{gather*}
\log{L}(\sigma+w,\rho_K) 
\ll \log\log{X}
\end{gather*}
for $w \in \tilde{L}_1 \cup \tilde{L}_5$. 
Furthermore, we deduce from Lemma \ref{lem:3.2} that the upper bound
\begin{gather*}
\log{L}(\sigma+w,\rho_K) 
\ll \frac{(\log{X})^{2-2\sigma+2\kappa_1(\sigma)}}{\log\log{X}}
+\log\log{X}
\end{gather*}
is valid for $w \in \tilde{L}_2 \cup \tilde{L}_3 \cup \tilde{L}_4$. 
Then, by the choice of $a(\sigma)$, we obtain
\begin{gather*}
g_z(\sigma+w,\rho_K)
\ll \exp\left(A \frac{\log{X}}{\log\log{X}}\right)
\end{gather*}
for any $z, w \in \mathbb{C}$ with $|z| \leq \tilde{R}_\sigma(X)$ and $w \in \tilde{L}$. 
Similarly to \eqref{eq:07281703} and \eqref{eq:07281704}, $Y^w$ and the integral on $\Gamma(w)$ are estimated as
\begin{align*}
Y^w
&\ll
\begin{cases}
\displaystyle{X^\eta} 
& \text{if $w$ lies on $\tilde{L}_1$, $\tilde{L}_2$, $\tilde{L}_4$, $\tilde{L}_5$}, \\
\displaystyle{X^{-\eta \kappa_1(\sigma)}} 
& \text{if $w$ lies on $\tilde{L}_3$},  
\end{cases} \\
\frac{1}{2\pi i} \int_{\tilde{L}_j} \Gamma(w) \,dw
&\ll
\begin{cases}
\exp\left(-(\log{X})^2\right)
& \text{if $j=1,2,4,5$}, \\
1
& \text{if $j=3$}. 
\end{cases}
\end{align*}
Hence we obtain the conclusion by these estimates. 
\end{proof}

By using Propositions \ref{prop:5.1}, \ref{prop:5.2} and \ref{prop:5.4}, we obtain another asymptotic formula for $g_z(\sigma,\rho_K)$ under GRH. 

\begin{corollary}\label{cor:5.2}
Assume GRH. 
Let $\sigma>1/2$ be a real number. 
Choose a positive real number $\kappa_1(\sigma)$ depending on $\sigma$ such that $\sigma-\kappa_1(\sigma)>1/2$ holds. 
Take a cubic field $K \in L_3^\pm(X)$ arbitrarily. 
If $Y=X^\eta$ with $\eta>0$, then there exists a constant $0<\tilde{b}(\eta) \leq1$ depending only on $\eta$ such that 
\begin{gather*}
g_z(\sigma,\rho_K)
=\sum_{n=1}^{\infty} d_z(n,\rho_K)n^{-\sigma}e^{-n/Y}
+O\left(X^{-\eta \kappa_1(\sigma)} \exp\left(A \frac{\log{X}}{\log\log{X}}\right)\right)
\end{gather*}
holds for $\sigma>1/2$ with $z \in \mathbb{C}$ satisfying $|z| \leq \tilde{b}(\eta) \tilde{R}_\sigma(X)$, where $\tilde{R}_\sigma(X)$ is given by \eqref{eq:07301552} with $a(\sigma)=\min(2\sigma-2\kappa_1(\sigma)-1, 2/3)>0$. 
Here, the constant $A \geq1$ is absolute, and the implied constant depends only on $\sigma$ and $\eta$. 
\end{corollary}

\subsection{Calculation of the complex moments}\label{sec:5.2}
Let $H_r(z;\mathfrak{a})$ be the coefficients in power series \eqref{eq:07272133}. 
Recall that $d_z(n,\rho_K)$ is defined as the multiplicative function satisfying $d_z(p^r,\rho_K)=H_r(z;\mathfrak{a})$ if $K$ satisfies a local specification $\mathfrak{a}$ at $p$. 
In this subsection, we calculation the $z$-th moment 
\begin{gather*}
M_{z, \sigma}^{\pm}(X)
=\sum_{K \in L_3^\pm(X) \setminus E_\sigma(X)} L(\sigma,\rho_K)^z 
=\sum_{K \in L_3^\pm(X) \setminus E_\sigma(X)} g_z(\sigma,\rho_K), 
\end{gather*}
where $E_\sigma(X)$ is the subset defined by $E_\sigma(X)=E(X)$ for $\sigma_1<\sigma<1$ and $E_\sigma(X)=\emptyset$ for $\sigma \geq1$ as in Section \ref{sec:2}. 
The following proposition is a key tool for the calculation of $M_{z, \sigma}^{\pm}(X)$. 

\begin{proposition}\label{prop:5.5}
Let $\alpha$ and $\beta$ be constants for which \eqref{eq:07251548} holds. 
For $z \in \mathbb{C}$, we denote by $\lambda_z(n)$, $\mu_z(n)$ and $\nu_z(n)$ the multiplicative functions satisfying 
\begin{align*}
\lambda_z(p^r)
&=\sum_{\mathfrak{a} \in \mathcal{A}} C_p(\mathfrak{a}) H_r(z;\mathfrak{a}),\\
\mu_z(p^r)
&=\sum_{\mathfrak{a} \in \mathcal{A}} K_p(\mathfrak{a}) H_r(z;\mathfrak{a}),\\
\nu_z(p^r)
&=\sum_{\mathfrak{a} \in \mathcal{A}} \left|H_r(z;\mathfrak{a})\right|, 
\end{align*}
respectively. 
Then we have the asymptotic formula
\begin{align*}
&\sum_{K \in L_3^\pm(X)} d_z(n,\rho_K) \\
&=C^\pm \frac{1}{12\zeta(3)} X \lambda_z(n)
+K^\pm \frac{4\zeta(1/3)}{5\Gamma(2/3)^3 \zeta(5/3)} X^{5/6} \mu_z(n)
+O\left(X^{\alpha+\epsilon} \nu_z(n) n^\beta \right)
\end{align*}
for each $\epsilon>0$, where the implied constant depends only on $\epsilon$.
\end{proposition}

\begin{proof}
The result for $n=1$ follows from \eqref{eq:07251548} with $\supp \mathcal{S}=\emptyset$. 
Thus we consider the case $n>1$ below. 
Let $n=p_1^{r_1} \cdots p_k^{r_k}$ be the prime factorization of $n$. 
For $\underline{\mathfrak{a}}=(\mathfrak{a}_1,\ldots,\mathfrak{a}_k) \in \mathcal{A}^k$, we denote by $\mathcal{S}(\underline{\mathfrak{a}})$ the local specifications of cubic fields such that $\supp \mathcal{S}(\underline{\mathfrak{a}})=\{p_1,\ldots,p_k\}$ and $\mathcal{S}_{p_j}=\mathfrak{a}_j$ for $j=1,\ldots,k$. 
Then 
\begin{gather*}
d_z(n,\rho_K)
=H_{r_1}(z;{\mathfrak{a}_1}) \cdots H_{r_k}(z;{\mathfrak{a}_k})
\end{gather*}
holds if $K$ satisfies the local specifications $\mathcal{S}(\underline{\mathfrak{a}})$. 
Dividing the condition $K \in L_3^\pm(X)$ into $K \in L_3^\pm(X, \mathcal{S}(\underline{\mathfrak{a}}))$ for $\underline{\mathfrak{a}}=(\mathfrak{a}_1,\ldots,\mathfrak{a}_k) \in \mathcal{A}^k$, we have
\begin{align*}
\sum_{K \in L_3^\pm(X)} d_z(n,\rho_K)
&=\sum_{\underline{\mathfrak{a}} \in \mathcal{A}^k} 
\sum_{K \in L_3^\pm(X, \mathcal{S}(\underline{\mathfrak{a}}))} d_z(n,\rho_K) \\
&=\sum_{\underline{\mathfrak{a}} \in \mathcal{A}^k} 
N_3^\pm(X,\mathcal{S}(\underline{\mathfrak{a}})) \prod_{j=1}^{k} H_{r_j}(z;\mathfrak{a}_j). 
\end{align*}
Then we apply \eqref{eq:07251548} to obtain
\begin{align*}
N_3^\pm(X,\mathcal{S}(\underline{\mathfrak{a}}))
&=C^\pm \frac{1}{12\zeta(3)} X \prod_{j=1}^{k} C_p(\mathfrak{a}_j)
+K^\pm \frac{4\zeta(1/3)}{5\Gamma(2/3)^3\zeta(5/3)} X^{5/6} \prod_{j=1}^{k} K_p(\mathfrak{a}_j) \\
&\quad+O_\epsilon\left(X^{\alpha+\epsilon} \prod_{j=1}^{k} p_j^{\beta}\right).
\end{align*}
From the above, we deduce
\begin{align*}
\sum_{K \in L_3^\pm(X)} d_z(n,\rho_K)
&=C^\pm \frac{1}{12\zeta(3)} X 
\sum_{\underline{\mathfrak{a}} \in \mathcal{A}^k} \prod_{j=1}^{k} C_p(\mathfrak{a}_j) H_{r_j}(z;\mathfrak{a}_j) \\
&\quad+K^\pm \frac{4\zeta(1/3)}{5\Gamma(2/3)^3\zeta(5/3)} X^{5/6} 
\sum_{\underline{\mathfrak{a}} \in \mathcal{A}^k} \prod_{j=1}^{k} K_p(\mathfrak{a}_j) H_{r_j}(z;\mathfrak{a}_j) \\
&\quad+O_\epsilon\left(X^{\alpha+\epsilon} 
\sum_{\underline{\mathfrak{a}} \in \mathcal{A}^k} \prod_{j=1}^{k} |H_{r_j}(z;\mathfrak{a}_j)| p_j^{\beta}\right).
\end{align*}
Hence we obtain the conclusion by noting that $\lambda_z(n)$ satisfies
\begin{gather*}
\lambda_z(n)
=\prod_{j=1}^{k} \left(\sum_{\mathfrak{a}_j \in \mathcal{A}} C_{p_j}(\mathfrak{a}_j) H_{r_j}(z;\mathfrak{a}_j) \right)
=\sum_{\underline{\mathfrak{a}} \in \mathcal{A}^k} \prod_{j=1}^{k} C_p(\mathfrak{a}_j) H_{r_j}(z;\mathfrak{a}_j), 
\end{gather*}
and that similar equalities are valid for $\mu_z(n)$ and $\nu_z(n)$. 
\end{proof}

Let $\mathcal{X}=(\mathcal{X}_p)_p$ and $\mathcal{Y}=(\mathcal{Y}_p)_p$ denote the sequences of independent random elements on $\{A_\mathfrak{a} \mid \mathfrak{a} \in \mathcal{A}\}$ as in Section \ref{sec:2}. 
Furthermore, we define $F_s(z)$ and $G_s(z)$ as 
\begin{gather*}
F_s(z)
=\mathbb{E}\left[\exp(z \log{L}(s,\mathcal{X}))\right]
\quad\text{and}\quad
F_s(z)
=\mathbb{E}\left[\exp(z \log{L}(s,\mathcal{Y}))\right]
\end{gather*}
as in Proposition \ref{prop:4.2}. 
If we define $F_{s,p}(z)=\mathbb{E}\left[\det\left(I-p^{-s} \mathcal{X}_p\right)^{-z} \right]$ and $G_{s,p}(z)=\mathbb{E}\left[\det\left(I-p^{-s} \mathcal{Y}_p\right)^{-z} \right]$ for any prime number $p$, then we have by \eqref{eq:07272133} the identities
\begin{align*}
F_{s,p}(z)
=\sum_{\mathfrak{a} \in \mathcal{A}} C_p(\mathfrak{a}) \det(I-p^{-s} A_\mathfrak{a})^{-z} 
=\sum_{r=0}^{\infty} \lambda_z(p^r)p^{-rs} \\
G_{s,p}(z)
=\sum_{\mathfrak{a} \in \mathcal{A}} K_p(\mathfrak{a}) \det(I-p^{-s} A_\mathfrak{a})^{-z} 
=\sum_{r=0}^{\infty} \mu_z(p^r)p^{-rs}. 
\end{align*}
Hence we derive
\begin{gather}\label{eq:07301318}
F_s(z)
=\prod_{p} \left(\sum_{r=0}^{\infty} \lambda_z(p^r) p^{-rs}\right)
=\sum_{n=1}^{\infty} \lambda_z(n) n^{-s}
\end{gather}
for $\RE(s)>1/2$ by formula \eqref{eq:07271532}, and similarly, 
\begin{gather}\label{eq:07301319}
G_s(z)
=\prod_{p} \left(\sum_{r=0}^{\infty} \mu_z(p^r) p^{-rs}\right)
=\sum_{n=1}^{\infty} \mu_z(n) n^{-s}
\end{gather}
for $\RE(s)>2/3$ by formula \eqref{eq:07271537}.

\subsubsection{Proof of Theorem \ref{thm:2.2}}\label{sec:5.2.1}
Let $\sigma_1$ be a real number for which \eqref{eq:07251505} holds. 
Let $Y=X^\eta$ with $\eta>0$ chosen later. 
By Corollary \ref{cor:5.1}, there exists a constant $0<b(\eta) \leq1$ depending only on $\eta$ such that 
\begin{gather}\label{eq:07292205}
\sum_{K \in L_3^\pm(X) \setminus E(X)} g_z(\sigma,\rho_K)
=S_1-S_2+O\left(X \exp\left(-\frac{\eta}{2}\frac{\log{X}}{\log\log{X}}\right)\right)
\end{gather}
holds for $z \in \mathbb{C}$ satisfying $|z| \leq b(\eta) R_\sigma(X)$, where 
\begin{align}
S_1
&=\sum_{K \in L_3^\pm(X)} \sum_{n=1}^{\infty} d_z(n,\rho_K) n^{-\sigma} e^{-n/Y}, 
\label{eq:07301700}\\
S_2
&=\sum_{K \in E(X)} \sum_{n=1}^{\infty} d_z(n,\rho_K) n^{-\sigma} e^{-n/Y}. 
\nonumber
\end{align}
Here, $E(X)$ is the subset of $L_3^\pm(X)$ defined by \eqref{eq:07252237}, and $R_\sigma(X)$ is as in \eqref{eq:07281607}. 

The main term comes from $S_1$. 
Recall that \eqref{eq:07251548} holds for $\alpha=7/9$ and $\beta=16/9$ due to Taniguchi--Thorne \cite{TaniguchiThorne2013}. 
By Proposition \ref{prop:5.5}, we obtain
\begin{align*}
S_1
&=C^\pm \frac{1}{12\zeta(3)} X 
\left(\sum_{n=1}^{\infty} \lambda_z(n) n^{-\sigma} e^{-n/Y}\right) \\
&\quad+K^\pm \frac{4\zeta(1/3)}{5\Gamma(2/3)^3\zeta(5/3)} X^{5/6} 
\left(\sum_{n=1}^{\infty} \mu_z(n) n^{-\sigma} e^{-n/Y} \right) \\
&\quad+O\left(X^{7/9+\epsilon} 
\sum_{n=1}^{\infty} \nu_z(n) n^{-\sigma+16/9} e^{-n/Y}\right) 
\end{align*}
for any $\epsilon>0$. 
Note that formula \eqref{eq:07301318} yields
\begin{gather*}
\sum_{n=1}^{\infty} \lambda_z(n) n^{-\sigma} e^{-n/Y}
=\frac{1}{2\pi i} \int_{\RE(w)=c} F_{\sigma+w}(z) \Gamma(w) Y^w \,dw
\end{gather*}
for $\sigma>\sigma_1$ and $c>0$. 
The function $F_{\sigma+w}(z)$ is a holomorphic function in $w$ on the half-plane $\RE(w)>1/2-\sigma$ by Proposition \ref{prop:4.2} $(\mathrm{i})$. 
Hence, by shifting the integral contour to left, we obtain 
\begin{gather}\label{eq:07301349}
\sum_{n=1}^{\infty} \lambda_z(n) n^{-\sigma} e^{-n/Y}
=F_\sigma(z)
+\frac{1}{2\pi i} \int_{\RE(w)=-\kappa_1} F_{\sigma+w}(z) \Gamma(w) Y^w \,dw
\end{gather}
with $0<\kappa_1<\sigma-1/2$. 
We choose the parameter $\kappa_1$ as 
\begin{gather*}
\kappa_1
=
\begin{cases}
\frac{1}{2} \min(\sigma-\sigma_1, 1-\sigma) 
& \text{if $\sigma_1<\sigma<1$}, \\
\sigma-1+(\log\log{X})^{-1} 
& \text{if $\sigma\geq1$}
\end{cases} 
\end{gather*}
so as to keep $1/2<\RE(\sigma+w)<1$ on the vertical line $\RE(w)=-\kappa_1$. 
Then we apply Proposition \ref{prop:4.3} $(\mathrm{i})$ to obtain
\begin{gather*}
F_{\sigma+w}(z)
\ll \exp\left( c_1(\sigma) \frac{(|z|+3)^{1/(\sigma-\kappa_1)}}{\log(|z|+3)} \right)
\end{gather*}
for $\RE(w)=-\kappa_1$, where $c_1(\sigma)$ is a positive constant depending only on $\sigma$. 
Therefore, there exists a small constant $b_\sigma(\eta)$ with $0<b_\sigma(\eta) \leq b(\eta)$ such that 
\begin{gather*}
F_{\sigma+w}(z)
\ll \exp\left( \frac{\eta}{2} \frac{\log{X}}{\log\log{X}} \right)
\end{gather*}
holds for any $z,w \in \mathbb{C}$ satisfying $|z| \leq b_\sigma(\eta) R_\sigma(X)$ and $\RE(w)=-\kappa_1$. 
For $Y=X^\eta$, we have 
\begin{gather*}
Y^w
\ll \exp\left( -\eta \frac{\log{X}}{\log\log{X}} \right)
\end{gather*}
on the line $\RE(w)=-\kappa_1$. 
Furthermore, the integral on $\Gamma(w)$ is estimated as
\begin{gather*}
\frac{1}{2\pi i} \int_{\RE(w)=-\kappa_1} \Gamma(w) \,dw
\ll \log\log{X}
\end{gather*}
by \eqref{eq:07281548}. 
The above implied constants depend at most on $\sigma$ and $\eta$. 
Combining these estimates, we deduce from \eqref{eq:07301349} the asymptotic formula
\begin{gather}\label{eq:07301431}
\sum_{n=1}^{\infty} \lambda_z(n) n^{-\sigma} e^{-n/Y}
=F_\sigma(z)
+O\left(\exp\left( -\frac{\eta}{2} \frac{\log{X}}{\log\log{X}} \right)\right)
\end{gather}
for $\sigma>\sigma_1$. 
Next, we obtain $|\mu_z(n)| \leq d_{2k}(n)$ and $|\nu_z(n)| \leq d_{2k}(n)$ with $k=\lfloor{|z|}\rfloor+1$ by Lemma \ref{lem:3.7}. 
Thus Lemma	\ref{lem:3.6} yields 
\begin{gather*}
\sum_{n=1}^{\infty} \mu_z(n) n^{-\sigma} e^{-n/Y}
\ll ({Y}^{1-\sigma}+1) (C \log{Y})^{2k+1}
\ll X^{\eta}, \\
\sum_{n=1}^{\infty} \nu_z(n) n^{-\sigma+16/9} e^{-n/Y}
\ll({Y}^{25/9-\sigma}+1) (C \log{Y})^{2k+1}
\ll X^{25\eta/9} 
\end{gather*}
for $\sigma>\sigma_1$. 
Making the parameter $\eta>0$ small enough, we then obtain
\begin{gather}
\sum_{n=1}^{\infty} \mu_z(n) n^{-\sigma} e^{-n/Y}
\ll X^{1/6} \exp\left( -\frac{\eta}{2} \frac{\log{X}}{\log\log{X}} \right), 
\label{eq:07301432}\\
\sum_{n=1}^{\infty} \nu_z(n) n^{-\sigma+16/9} e^{-n/Y}
\ll X^{2/9-\epsilon} \exp\left( -\frac{\eta}{2} \frac{\log{X}}{\log\log{X}} \right). 
\label{eq:07301433}
\end{gather}
By \eqref{eq:07301431}, \eqref{eq:07301432} and \eqref{eq:07301433}, we calculate the first term $S_1$ as
\begin{gather*}
S_1
=C^\pm \frac{1}{12\zeta(3)} X 
\int_{-\infty}^{\infty} e^{zx} C_\sigma(x) \,\frac{dx}{\sqrt{2\pi}} 
+O\left(X \exp\left( -\frac{\eta}{2} \frac{\log{X}}{\log\log{X}} \right)\right) 
\end{gather*}
for $\sigma>\sigma_1$ and $z \in \mathbb{C}$ such that $|z| \leq b_\sigma(\eta) R_\sigma(X)$ since $F_\sigma(z)$ is represented as
\begin{gather*}
F_\sigma(z)
=\mathbb{E}\left[\exp(z \log{L}(\sigma,\mathcal{X}))\right]
=\int_{-\infty}^{\infty} e^{zx} C_\sigma(x) \,\frac{dx}{\sqrt{2\pi}} 
\end{gather*}
by using the density function $C_\sigma$ of Theorem \ref{thm:2.1}. 

The second term $S_2$ is estimated as follows. 
First, we obtain by Proposition \ref{prop:5.2} the upper bound
\begin{gather*}
\max_{K \in L_3^\pm(X)} \left| \sum_{n=1}^{\infty} d_z(n,\rho_K) n^{-\sigma} e^{-n/Y} \right|
\ll X^\eta 
\end{gather*}
for $\sigma>\sigma_1$ and $|z| \leq R_\sigma(X)$. 
Furthermore, by the definition of $E(X)$, zero-density estimate \eqref{eq:07251505} deduces
\begin{gather}\label{eq:07301528}
\# E(X)
\ll \sum_{K \in L_3^\pm(X)} N(\sigma_1, (\log{X})^3; \rho_K)
\ll X^{1-\delta}
\end{gather}
with an absolute constant $\delta>0$. 
Hence we have 
\begin{gather*}
S_2
\ll \# E(X) \max_{K \in L_3^\pm(X)} \left| \sum_{n=1}^{\infty} d_z(n,\rho_K) n^{-\sigma} e^{-n/Y} \right|
\ll X^{1-\delta/2}
\end{gather*}
if we let $\eta<\delta/2$. 
As a result, we deduce from \eqref{eq:07292205} the formula
\begin{align}\label{eq:07301444}
&\sum_{K \in L_3^\pm(X) \setminus E(X)} g_z(\sigma,\rho_K) \\
&=C^\pm \frac{1}{12\zeta(3)} X 
\int_{-\infty}^{\infty} e^{zx} C_\sigma(x) \,\frac{dx}{\sqrt{2\pi}} 
+O\left(X\exp\left(-\delta' \frac{\log{X}}{\log\log{X}}\right)\right) \nonumber
\end{align}
for $\sigma>\sigma_1$ and $|z| \leq b_\sigma(\eta) R_\sigma(X)$, where $\delta'$ is an absolute positive constant. 

Note that the left-hand side of \eqref{eq:07301444} is equal to $M_{z, \sigma}^{\pm}(X)$ for $\sigma_1<\sigma<1$ since we have $E_\sigma(X)=E(X)$ in this case. 
Hence Theorem \ref{thm:2.2} follows if $\sigma_1<\sigma<1$. 
To complete the proof, we consider the case $\sigma \geq1$. 
We see that 
\begin{gather*}
M_{z, \sigma}^{\pm}(X)
=\sum_{K \in L_3^\pm(X) \setminus E(X)} g_z(\sigma,\rho_K)
+\sum_{K \in E(X)} g_z(\sigma,\rho_K)
\end{gather*}
since we have $E_\sigma(X)=\emptyset$ in this case. 
Thus Theorem \ref{thm:2.2} is deduced from if we evaluate the second sum as
\begin{gather}\label{eq:07301529}
\sum_{K \in E(X)} g_z(\sigma,\rho_K)
\ll X\exp\left(-\delta' \frac{\log{X}}{\log\log{X}}\right)
\end{gather}
with some absolute constant $\delta'>0$. 
Recall that $|\log{L}(\sigma,\rho_K)|\leq 2\log{\zeta(\sigma)}$ holds for $\sigma>1$. 
Thus, there exists a positive constant $b_\sigma(\delta)$ for $\delta>0$ such that 
\begin{gather*}
|g_z(\sigma,\rho_K)|
\leq \exp\left(\frac{\delta}{2} \frac{\log{X}}{\log\log{X}}\right)
\end{gather*}
holds for $|z| \leq b_\sigma(\delta) R_\sigma(X)$ in the case $\sigma>1$. 
Furthermore, we see that this is true even for $\sigma=1$ by using Lemma \ref{lem:3.4}. 
Together with \eqref{eq:07301528}, we finally obtain upper bound \eqref{eq:07301529} for $\sigma \geq1$. 
Hence the proof of Theorem \ref{thm:2.2} is completed.

\subsubsection{Proof of Theorem \ref{thm:2.3}}\label{sec:5.2.2}
Let $Y=X^\eta$ with a parameter $\eta>0$ determined later. 
First, we remark that the subset $E(X)$ of \eqref{eq:07252237} is empty for any $\sigma_1>1/2$ under GRH. 
Hence $E_\sigma(X)$ remains empty for $\sigma>1/2$, and thus, there exists a constant $0<\tilde{b}(\eta) \leq1$ such that
\begin{gather}\label{eq:07302135}
M_{z, \sigma}^{\pm}(X)
=S_1
+O\left(X^{1-\eta \kappa_1(\sigma)} \exp\left(A \frac{\log{X}}{\log\log{X}}\right)\right) 
\end{gather}
holds for $|z| \leq \tilde{b}(\eta) \tilde{R}_\sigma(X)$ by Corollary \ref{cor:5.2}, where $S_1$ is as in \eqref{eq:07301700}, and $\kappa_1(\sigma)$ is a real number chosen later such that $1/2<\sigma-\kappa_1(\sigma)<1$. 
By Proposition \ref{prop:5.5}, we calculate $S_1$ as
\begin{align}\label{eq:07302040}
S_1
&=C^\pm \frac{1}{12\zeta(3)} X 
\left(\sum_{n=1}^{\infty} \lambda_z(n) n^{-\sigma} e^{-n/Y}\right) \\
&\quad+K^\pm \frac{4\zeta(1/3)}{5\Gamma(2/3)^3\zeta(5/3)} X^{5/6} 
\left(\sum_{n=1}^{\infty} \mu_z(n) n^{-\sigma} e^{-n/Y} \right) \nonumber\\
&\quad+O\left(X^{\alpha+\epsilon} 
\sum_{n=1}^{\infty} \nu_z(n) n^{-\sigma+\beta} e^{-n/Y}\right), \nonumber
\end{align}
where $\alpha$ and $\beta$ are constants such that $0<\alpha<5/6$. 
We deduce from \eqref{eq:07301318} and \eqref{eq:07301319} the identities 
\begin{align*}
\sum_{n=1}^{\infty} \lambda_z(n) n^{-\sigma} e^{-n/Y}
&=F_\sigma(z)
+\frac{1}{2\pi i} \int_{\RE(w)=-\kappa_1(\sigma)} F_{\sigma+w}(z) \Gamma(w) Y^w \,dw, \\
\sum_{n=1}^{\infty} \mu_z(n) n^{-\sigma} e^{-n/Y}
&=G_\sigma(z)
+\frac{1}{2\pi i} \int_{\RE(w)=-\kappa_2(\sigma)} G_{\sigma+w}(z) \Gamma(w) Y^w \,dw 
\end{align*}
for any $\sigma>2/3$ similarly to \eqref{eq:07301349}, where $\kappa_2(\sigma)$ is a real number chosen later such that $2/3<\sigma-\kappa_2(\sigma)<1$. 
Then, we apply Proposition \ref{prop:4.3} $(\mathrm{i})$, $(\mathrm{ii})$ to evaluate $F_{\sigma+w}(z)$ on $\RE(w)=-\kappa_1(\sigma)$ and $G_{\sigma+w}(z)$ on $\RE(w)=-\kappa_2(\sigma)$. 
Let $z$ be a complex number satisfying $|z| \leq \tilde{b}_\sigma(\eta) \tilde{R}_\sigma(X)$ with some constant $\tilde{b}_\sigma(\eta)$ such that $0<\tilde{b}_\sigma(\eta) \leq \tilde{b}(\eta)$. 
We have 
\begin{gather*}
F_{\sigma+w}(z)
\ll \exp\left(c_1(\sigma) \tilde{b}_\sigma(\eta) \frac{\log{X}}{\log\log{X}}\right)
\end{gather*}
on $\RE(w)=-\kappa_1(\sigma)$ by noting that $a(\sigma)\leq 2\sigma-2\kappa_1(\sigma)-1$, and 
\begin{gather*}
G_{\sigma+w}(z)
\ll \exp\left(c_2(\sigma) \tilde{b}_\sigma(\eta) \frac{\log{X}}{\log\log{X}}\right)
\end{gather*}
on $\RE(w)=-\kappa_2(\sigma)$ by noting that $a(\sigma)\leq 2/3$. 
Here, $c_1(\sigma)$ and $c_2(\sigma)$ are positive constants depending only on $\sigma$. 
These upper bounds yield
\begin{align*}
\sum_{n=1}^{\infty} \lambda_z(n) n^{-\sigma} e^{-n/Y}
&=F_\sigma(z)
+O\left(X^{-\eta \kappa_1(\sigma)} 
\exp\left(\frac{\log{X}}{\log\log{X}}\right)\right), \\
\sum_{n=1}^{\infty} \mu_z(n) n^{-\sigma} e^{-n/Y}
&=G_\sigma(z)
+O\left(X^{-\eta \kappa_2(\sigma)} 
\exp\left(\frac{\log{X}}{\log\log{X}}\right)\right) 
\end{align*}
if we let the parameter $\tilde{b}_\sigma(\eta)$ be small enough. 
Furthermore, we have 
\begin{gather*}
\sum_{n=1}^{\infty} \nu_{z}(n) n^{-\sigma+\beta} e^{-n/Y}
\ll (X^{\eta (1-\sigma+\beta)}+1) \exp\left(\frac{\log{X}}{\log\log{X}}\right) 
\end{gather*}
by Lemma \ref{lem:3.6} since $|\nu_{z}(n)| \leq d_{2k}(n)$ with $k=\lfloor{|z|}\rfloor+1$. 
Hence we deduce from \eqref{eq:07302040} that $S_1$ is calculated as
\begin{gather}\label{eq:07302136}
S_1
=C^\pm \frac{1}{12\zeta(3)} X F_\sigma(z)
+K^\pm \frac{4\zeta(1/3)}{5\Gamma(2/3)^3\zeta(5/3)} X^{5/6} G_\sigma(z)
+E, 
\end{gather}
where the error term $E$ satisfies
\begin{align*}
E
&\ll \left( X^{1-\eta \kappa_1(\sigma)} 
+X^{5/6-\eta \kappa_2(\sigma)}
+X^{\alpha+\eta (1-\sigma+\beta)+\epsilon}
+X^{\alpha+\epsilon} \right) 
\exp\left(A \frac{\log{X}}{\log\log{X}}\right) \\
&\ll 
X^{5/6+\epsilon} 
\left( X^{-\eta \kappa_1(\sigma)+1/6} 
+X^{-\eta \kappa_2(\sigma)}
+X^{-(5/6-\alpha)+\eta (1-\sigma+\beta)}
+X^{-(5/6-\alpha)} \right). 
\end{align*}
with the implied constants depend at most on $\sigma$, $\eta$ and $\epsilon$. 
Then we show that it is evaluated as $E \ll X^{5/6-\delta}$ with some $\delta=\delta(\sigma)>0$ as follows. 
\begin{itemize}
\item[$(\mathrm{a})$]
Let $\sigma \geq 1+\beta$. 
In this case, we choose the constants $\kappa_1(\sigma)$ and $\kappa_2(\sigma)$ satisfying $1/2<\sigma-\kappa_1(\sigma)<1$ and $2/3<\sigma-\kappa_2(\sigma)<1$ arbitrarily. 
Then, we take a real number $\eta$ so that $\eta \kappa_1(\sigma) \geq 1/3$ is satisfied. 
In this setting, we obtain
\begin{gather*}
E
\ll X^{5/6+\epsilon} \left(X^{-1/6}+X^{-\eta \kappa_2(\sigma)}+{X}^{-(5/6-\alpha)}\right)
\ll X^{5/6+\epsilon-\delta}
\end{gather*}
with some constant $\delta=\delta(\sigma)>0$ depending only on $\sigma$. 
\item[$(\mathrm{b})$]
Let $\max(\sigma_2,2/3) < \sigma < 1+\beta$, where $\sigma_2$ is the real number of \eqref{eq:07302054}. 
In this case, we need to choose $\kappa_1(\sigma)$ and $\eta$ more carefully. 
Note that the inequality 
\begin{gather*}
\sigma-1
<\frac{1-\sigma+\beta}{5-6\alpha}
<\sigma-\frac{1}{2}
\end{gather*}
holds by the definition $\sigma_2$. 
Thus we choose a constant $\kappa_1(\sigma)>0$ so that 
\begin{gather*}
\frac{1-\sigma+\beta}{5-6\alpha}
<\kappa_1(\sigma)
<\sigma-\frac{1}{2}
\end{gather*}
is satisfied. 
Therefore we obtain
\begin{gather*}
\frac{1}{2}<\sigma-\kappa_1(\sigma)<1
\quad\text{and}\quad
-\frac{(1-\alpha) \kappa_1(\sigma)}{1-\sigma+\beta+\kappa_1(\sigma)}+\frac{1}{6}<0
\end{gather*}
by the choice of $\kappa_1(\sigma)$. 
Then, we choose a real number $\eta$ as 
\begin{gather*}
\eta
=\frac{1-\alpha}{1-\sigma+\beta+\kappa_1(\sigma)}>0
\end{gather*}
to satisfy $X^{-\eta \kappa_1(\sigma)+1/6}=X^{-(5/6-\alpha)+\eta(1-\sigma+\beta)}$. 
Finally, we take $\kappa_2(\sigma)>0$ such that $2/3<\sigma-\kappa_2(\sigma)<1$ arbitrarily. 
From the above, we conclude 
\begin{gather*}
E
\ll X^{5/6+\epsilon} \left(X^{-\eta \kappa_1(\sigma)+1/6}+X^{-\eta \kappa_2(\sigma)}++X^{-(5/6-\alpha)}\right)
\ll X^{5/6+\epsilon-\delta}, 
\end{gather*}
where $\delta=\delta(\sigma)>0$ is a constant depending only on $\sigma$. 
\end{itemize}
Therefore, we obtain $E \ll X^{5/6-\delta/2}$ in both cases if we let $\epsilon=\delta/2$. 
Recall that $F_\sigma(z)$ and $G_\sigma(z)$ are represented as
\begin{gather*}
F_\sigma(z)
=\int_{-\infty}^{\infty} e^{zx} C_\sigma(x) \,\frac{dx}{\sqrt{2\pi}} 
\quad\text{and}\quad
G_\sigma(z)
=\int_{-\infty}^{\infty} e^{zx} K_\sigma(x) \,\frac{dx}{\sqrt{2\pi}}. 
\end{gather*}
Hence we derive Theorem \ref{thm:2.3} by \eqref{eq:07302135} and \eqref{eq:07302136}.

\subsection{Applications of the class number formula}\label{sec:5.3}
Let $r$ be a real number. 
For an integer $d>0$, we define 
\begin{gather*}
f_{r,\pm}(d)
=\sum_{\substack{K \in L_3^\pm(X) \\ |d_K|=d}} \left(\frac{h_K R_K}{\sqrt{|d_K|}}\right)^r. 
\end{gather*}
Then we obtain 
\begin{align}\label{eq:07302222}
\sum_{K \in L_3^\pm(X)} (h_K R_K)^r
&=\sum_{0<d \leq X} f_{r,\pm}(d) d^{r/2} \\
&=F_{r,\pm}(X) X^{r/2} 
-\frac{r}{2} \int_{1}^{X} F_{r,\pm}(y) y^{r/2-1} \,dy \nonumber
\end{align}
by the partial summation, where the function $F_{r,\pm}(y)$ is given by
\begin{gather*}
F_{r,\pm}(y)
=\sum_{0<d\leq y} f_{r,\pm}(d)
=\sum_{K \in L_3^\pm(y)} \left(\frac{h_K R_K}{\sqrt{|d_K|}}\right)^r. 
\end{gather*}

\subsubsection{Proof of Corollary \ref{cor:2.1}}\label{sec:5.3.1}
By formula \eqref{eq:07250142}, we have 
\begin{gather*}
F_{r,\pm}(y)
=\frac{1}{(D^\pm)^r} \sum_{K \in L_3^\pm(y)} L(1,\rho_K)^r. 
\end{gather*}
Hence we deduce from Theorem \ref{thm:2.2} the formula
\begin{gather}\label{eq:07302221}
F_{r,\pm}(y)
=C^\pm(r) y
+O\left(y \exp\left(-\delta \frac{\log{y}}{\log\log{y}}\right)\right)
\end{gather}
with an absolute constant $\delta>0$, where we write
\begin{gather*}
C^\pm(r)
=C^\pm \frac{1}{12\zeta(3)} 
\frac{1}{(D^\pm)^r} \int_{-\infty}^{\infty} e^{rx} C_1(x) \,\frac{dx}{\sqrt{2\pi}},  
\end{gather*}
and the implied constant depends only on $r$. 
Let $r>-2$ be a fixed real number. 
Inserting \eqref{eq:07302221} to formula \eqref{eq:07302222}, we derive
\begin{align*}
&\sum_{K \in L_3^\pm(X)} (h_K R_K)^r \\
&=C^\pm(r) X^{r/2+1}
-\frac{r}{2} \int_{1}^{X} C^\pm(r) y^{r/2} \,dy 
+O\left(X^{r/2+1} \exp\left(-\delta \frac{\log{X}}{\log\log{X}}\right)\right) \\
&=C^\pm(r) \frac{2}{r+2} X^{r/2+1}+O\left(X^{r/2+1} \exp\left(-\delta \frac{\log{X}}{\log\log{X}}\right)\right)
\end{align*}
as desired.

\begin{remark}\label{rem:5.1}
If we let $r=1$ in Corollary \ref{cor:2.1}, then we have 
\begin{gather*}
\sum_{K \in L_3^\pm(X)} h_K R_K
=\frac{C^\pm}{D^\pm} \frac{F_1(1)}{18\zeta(3)} X^{3/2}
+O\left(X^{3/2} \exp\left(-\delta \frac{\log{X}}{\log\log{X}}\right)\right), 
\end{gather*}
where $=\mathbb{E}\left[\exp(z \log{L}(s,\mathcal{X}))\right]$ as in Proposition \ref{prop:4.2} $(\mathrm{i})$. 
Here, we have $C^+/D^+=1/4$ and $C^-/D^-=3/(2\pi)$. 
Furthermore, we obtain $F_1(1)=\prod_{p} F_{1,p}(1)$ by formula \eqref{eq:07271532}, where $F_{1,p}(1)$ is calculated as
\begin{align*}
F_{1,p}(1)
=&\frac{1}{1+p^{-1}+p^{-2}}
\biggl[\frac{1}{6}(1-p^{-1})^{-2}+\frac{1}{2}(1-p^{-2})^{-1} \\
&\qquad\qquad\qquad\quad
+\frac{1}{3}(1+p^{-1}+p^{-2})^{-1}+\frac{1}{p}(1-p^{-1})^{-z}+\frac{1}{p^2}\biggr]\\
=&(1-p^{-3})^{-2}(1-p^{-2})^{-1}(1+p^{-2}-2p^{-3}-2p^{-4}+2p^{-6}+p^{-7}-p^{-8}). 
\end{align*} 
Thus we derive the asymptotic formula
\begin{gather*}
\sum_{K \in L_3^\pm(X)} h_K R_K
=\frac{C^\pm}{D^\pm} 4 c X^{3/2}
+O\left(X^{3/2}\exp\left(-\delta\frac{\log{X}}{\log\log{X}}\right)\right)
\end{gather*}
with $c>0$ described in Theorem \ref{thm:1.2}. 
\end{remark}

\subsubsection{Proof of Corollary \ref{cor:2.2}}\label{sec:5.3.2}
If we assume that \eqref{eq:07251548} holds with some $\alpha$ and $\beta$ such that $3\alpha+\beta<5/2$, then Theorem \ref{thm:2.3} is available at $\sigma=1$. 
Thus we obtain 
\begin{gather*}
F_{r,\pm}(y)
=C^\pm(r) y
+K^\pm(r) y^{5/6}
+O\left(y^{5/6} \exp\left(-\delta \frac{\log{y}}{\log\log{y}}\right)\right)
\end{gather*}
similarly to \eqref{eq:07302221}, where $C^\pm(r)$ is as in the proof of Corollary \ref{cor:2.1}, and 
\begin{gather*}
K^\pm(r)
=K^\pm \frac{4\zeta(1/3)}{5\Gamma(2/3)^3\zeta(5/3)} 
\frac{1}{(D^\pm)^r} \int_{-\infty}^{\infty} e^{rx} K_1(x) \,\frac{dx}{\sqrt{2\pi}}. 
\end{gather*}
By \eqref{eq:07302222}, we have 
\begin{align*}
\sum_{K \in L_3^\pm(X)} (h_K R_K)^r 
&=C^\pm(r) X^{r/2+1}
-\frac{r}{2} \int_{1}^{X} C^\pm(r) y^{r/2} \,dy \\
&\quad+K^\pm(r) X^{r/2+5/6}
-\frac{r}{2} \int_{1}^{X} K^\pm(r) y^{r/2-1/6} \,dy \\
&\quad+O\left(X^{r/2+5/6} \exp\left(-\delta \frac{\log{X}}{\log\log{X}}\right)\right) \\
&=C^\pm(r) \frac{2}{r+2} X^{r/2+1}
+K^\pm(r) \frac{5}{3r+5} X^{r/2+5/6} \\
&\quad+O\left(X^{r/2+5/6} \exp\left(-\delta \frac{\log{X}}{\log\log{X}}\right)\right)
\end{align*}
if $r/2+5/6>0$, i.e.\ $r>-5/3$. 
Hence we obtain the desired result.

\section{Completion of the proofs}\label{sec:6}
The remaining work is to complete the proofs of Theorems \ref{thm:2.4} and \ref{thm:2.5}. 
We begin with showing the following corollary of Theorem \ref{thm:2.2} which is used in the proof of Theorem \ref{thm:2.4}. 

\begin{corollary}\label{cor:6.1}
Let $k$ be a positive integer. 
Let $\sigma_1$ be a real number for which \eqref{eq:07251505} holds. 
Then there exists an absolute constant $\delta>0$ such that 
\begin{align*}
&\frac{1}{N_3^\pm(X)} \sum_{K \in L_3^\pm(X) \setminus E_\sigma(X)} (\log{L}(\sigma,\rho_K))^k \\
&=\int_{-\infty}^{\infty} x^k C_1(x) \,\frac{dx}{\sqrt{2\pi}}
+O\left(\exp\left(-\delta \frac{\log{X}}{\log\log{X}}\right)\right)
\end{align*}
holds for $\sigma>\sigma_1$, where the implied constant depends only on $k$ and $\sigma$. 
Here, $E_\sigma(X)$ is the subset defined by $E_\sigma(X)=E(X)$ for $\sigma_1<\sigma<1$ and $E_\sigma(X)=\emptyset$ for $\sigma \geq1$ as in Section \ref{sec:2}. 
\end{corollary}

\begin{proof}
We define two holomorphic functions $f$ and $g$ as 
\begin{align*}
f(z)
&=\frac{1}{N_3^\pm(X)} \sum_{K \in L_3^\pm(X) \setminus E_\sigma(X)} \exp(z \log{L}(\sigma,\rho_K)), \\
g(z)
&=\int_{-\infty}^{\infty} e^{zx} C_1(x) \,\frac{dx}{\sqrt{2\pi}}
\end{align*}
for $z \in \mathbb{C}$. 
Let $\epsilon>0$. 
Then it is deduced from Theorem \ref{thm:2.2} that 
\begin{gather*}
|f(z)-g(z)|
\leq A(\sigma) \exp\left(-\delta \frac{\log{X}}{\log\log{X}}\right)
\end{gather*}
holds in the range $|z|<\epsilon$, where $A(\sigma)$ is a positive constant depending on $\sigma$. 
Hence, as a consequence of Cauchy's integral formula, we obtain
\begin{align*}
&\frac{1}{N_3^\pm(X)} \sum_{K \in L_3^\pm(X) \setminus E_\sigma(X)} (\log{L}(\sigma,\rho_K))^k 
-\int_{-\infty}^{\infty} x^k C_1(x) \,\frac{dx}{\sqrt{2\pi}} \\
&=\frac{d^k}{dz^k} (f(z)-g(z))\bigg|_{z=0}
\ll_{k,\sigma} \exp\left(-\delta \frac{\log{X}}{\log\log{X}}\right), 
\end{align*}
which completes the proof. 
\end{proof}

\subsection{Proof of Theorem \ref{thm:2.4}}\label{sec:6.1}
We apply Lemma \ref{lem:3.12} for the probability measures $\mu$ and $\nu$ defined as
\begin{align*}
\mu(A)
&=\frac{\# \{K \in L_3^\pm(X) \setminus E_\sigma(X) \mid \log{L}(\sigma,\rho_K) \in A \}}
{\#(L_3^\pm(X) \setminus E_\sigma(X))}, \\
\nu(A)
&=\int_{A} C_\sigma(x) \,\frac{dx}{\sqrt{2\pi}}
\end{align*}
for $A \in \mathcal{B}(\mathbb{R})$. 
Note that the characteristic function of $\mu$ is calculated as
\begin{align*}
\phi(\xi)
&=\frac{1}{\#(L_3^\pm(X) \setminus E_\sigma(X))} 
\sum_{K \in L_3^\pm(X) \setminus E_\sigma(X)} L(\sigma,\rho_K)^{i \xi}, \\
&=\frac{1}{N_3^\pm(X)} M_{i \xi,\sigma}^\pm(X)
+O(X^{-\delta}) 
\end{align*}
since $\# E_\sigma(X) \ll X^{1-\delta}$ by zero-density estimate \eqref{eq:07251505}. 
Furthermore, the characteristic function of $\nu$ is given by
\begin{gather*}
\psi(\xi)
=\int_{-\infty}^{\infty} e^{i \xi x} C_\sigma(x) \,\frac{dx}{\sqrt{2\pi}}. 
\end{gather*}
Thus Theorem \ref{thm:2.2} yields 
\begin{gather}\label{eq:07310138}
\phi(\xi)-\psi(\xi)
\ll \exp\left(-\delta \frac{\log{X}}{\log\log{X}}\right)
\end{gather}
for $|\xi| \leq b_\sigma R_\sigma(X)$, where the implied constant depends only on $\sigma$. 
Here, $b_\sigma$ is a positive constant, and $R_\sigma(X)$ is given by \eqref{eq:07281607}. 
Put $R=b_\sigma R_\sigma(X)$ and $r=(\log{X})^{-2}$. 
Then we obtain 
\begin{align}\label{eq:07310156}
\int_{r}^{R} \left|\frac{\phi(\xi)-\psi(\xi)}{\xi}\right| \,d \xi
&\ll \exp\left(-\delta \frac{\log{X}}{\log\log{X}}\right) 
\int_{r}^{R} \frac{1}{\xi} \,d \xi \\
&\ll \exp\left(-\frac{\delta}{2} \frac{\log{X}}{\log\log{X}}\right) \nonumber
\end{align}
by \eqref{eq:07310138}. 
We have also 
\begin{gather}\label{eq:07310157}
\int_{-R}^{-r} \left|\frac{\phi(\xi)-\psi(\xi)}{\xi}\right| \,d \xi
\ll \exp\left(-\frac{\delta}{2} \frac{\log{X}}{\log\log{X}}\right). 
\end{gather}
Let $-r<\xi<r$. 
Then, we recall that $e^{i \theta}=1+O(|\theta|)$ holds uniformly for $\theta \in \mathbb{R}$. 
Therefore $\phi(\xi)$ is approximated as
\begin{align*}
\phi(\xi)
&=1+O\left(\frac{1}{N_3^\pm(X)} 
\sum_{K \in L_3^\pm(X) \setminus E_\sigma(X)} |\xi| |\log{L}(\sigma,\rho_K)|\right) \\
&=1+O\left(|\xi| \left(\frac{1}{N_3^\pm(X)} 
\sum_{K \in L_3^\pm(X) \setminus E_\sigma(X)} (\log{L}(\sigma,\rho_K))^2\right)^{1/2}\right) \\
&=1+O\left(|\xi| \left(\int_{-\infty}^{\infty} x^2 C_1(x) \,\frac{dx}{\sqrt{2\pi}}\right)^{1/2}\right)
\end{align*}
by the Cauchy--Schwarz inequality and Corollary \ref{cor:6.1} with $k=2$. 
In a similar way, we derive
\begin{gather*}
\psi(\xi)
=1+O\left(|\xi| \left(\int_{-\infty}^{\infty} x^2 C_1(x) \,\frac{dx}{\sqrt{2\pi}}\right)^{1/2}\right). 
\end{gather*}
Then we see that $\phi(\xi)-\psi(\xi)$ is estimated as $\phi(\xi)-\psi(\xi) \ll |\xi|$, which yields 
\begin{gather}\label{eq:07310158}
\int_{-r}^{r} \left|\frac{\phi(\xi)-\psi(\xi)}{\xi}\right| \,d \xi
\ll \int_{-r}^{r} \,d \xi
\ll (\log{X})^{-2}. 
\end{gather}
Combining \eqref{eq:07310156}, \eqref{eq:07310157} and \eqref{eq:07310158}, we obtain that the integral in the right-hand side of \eqref{eq:07310200} is estimated as
\begin{gather*}
\int_{-R}^{R} \left|\frac{\phi(\xi)-\psi(\xi)}{\xi}\right| \,d \xi
\ll (\log{X})^{-2}. 
\end{gather*}
Recall that $\sup_{x \in\mathbb{R}} C_\sigma(x) \ll_\sigma 1$. 
As a result, we drive by Lemma \ref{lem:3.12} the upper bound
\begin{gather*}
\sup_{t \in\mathbb{R}} |F(t)-G(t)|
\ll \frac{1}{R_\sigma(X)} 
+(\log{X})^{-2}
\ll \frac{1}{R_\sigma(X)}, 
\end{gather*}
where $F(t)=\mu((-\infty, t])$ and $G(t)=\nu((-\infty, t])$, and the implied constant depends only on $\sigma$. 
To end the proof, we see that the quantity $D_\sigma^\pm(X; a)$ is represented as 
\begin{align*}
D_\sigma^\pm(X; a) 
&=\frac{\# \left\{ K \in L_3^\pm(X) ~\middle|~ L(\sigma,\rho_K) \leq e^a \right\}}{N_3^\pm(X)}
-G(a) \\
&= F(a)-G(a)+E, 
\end{align*}
where 
\begin{gather*}
E
=\frac{\# \left\{ K \in L_3^\pm(X) ~\middle|~ L(\sigma,\rho_K) \leq e^a \right\}}{N_3^\pm(X)} 
-\frac{\# \{K \in L_3^\pm(X) \setminus E_\sigma(X) \mid L(\sigma,\rho_K) \leq e^a \}}
{\#(L_3^\pm(X) \setminus E_\sigma(X))}. 
\end{gather*}
We have $E \ll X^{-\delta}$ by using $\# E_\sigma(X) \ll X^{1-\delta}$. 
Hence we obtain the conclusion
\begin{gather*}
\sup_{a \in \mathbb{R}} \left|D_\sigma^\pm(X; a)\right|
\ll \sup_{a \in\mathbb{R}} |F(a)-G(a)|
+X^{-\delta}
\ll \frac{1}{R_\sigma(X)}. 
\end{gather*}

\subsection{Proof of Theorem \ref{thm:2.5}}\label{sec:6.2}
Let $A_\sigma(X)$ be the set of \eqref{eq:07252224} for $\sigma_1<\sigma<1$ and $A_\sigma(X)=\emptyset$ for $\sigma \geq1$. 
Then we have $\# A_\sigma(X) \ll X^{1-\delta}$ with some $\delta>0$ by zero-density estimate \eqref{eq:07251505}. 
Therefore, limit formula \eqref{eq:07310245} follows if we obtain
\begin{gather}\label{eq:07310306}
\lim_{X \to\infty} \frac{1}{\#(L_3^\pm(X) \setminus A_\sigma(X))} 
\sideset{}{'}\sum_{K \in L_3^\pm(X)} \Phi\left(\log{L}(\sigma,\rho_K)\right) 
=\int_{-\infty}^{\infty} \Phi(u) C_\sigma(x) \,\frac{dx}{\sqrt{2\pi}}. 
\end{gather}
Let $\mu_{\sigma,X}$ and $\mu_\sigma$ be the probability measures on $(\mathbb{R},\mathcal{B}(\mathbb{R}))$ defined as
\begin{align*}
\mu_{\sigma,X}(A)
&=\frac{\{K \in L_3^\pm(X) \setminus A_\sigma(X) \mid \log{L}(\sigma,\rho_K) \in A \}}
{\#(L_3^\pm(X) \setminus A_\sigma(X))}, \\
\mu_\sigma(A)
&=\int_{A} C_\sigma(x) \,\frac{dx}{\sqrt{2\pi}}. 
\end{align*}
Then the characteristic functions are given by
\begin{align*}
\phi_{\sigma,X}(\xi)
&=\frac{1}{\#(L_3^\pm(X) \setminus A_\sigma(X))} 
\sum_{K \in L_3^\pm(X) \setminus A_\sigma(X)} L(\sigma,\rho_K)^{i \xi} \\
&=\frac{1}{N_3^\pm(X)} M_{i \xi, \sigma}^\pm(X)
+O(X^{-\delta}), \\
\phi_\sigma(\xi)
&=\int_{-\infty}^{\infty} e^{i \xi x} C_\sigma(x) \,\frac{dx}{\sqrt{2\pi}}. 
\end{align*}
Therefore, Theorem \ref{thm:2.2} yields that $\phi_{\sigma,X}(\xi) \to \phi_\sigma(\xi)$ holds as $X \to\infty$ uniformly in $\xi \in [-R,R]$ for any $R>0$. 
By Lemma \ref{lem:3.9}, we obtain \eqref{eq:07310306} for $\Phi(x)=1_A(x)$ with any continuity set $A$ of $\mathbb{R}$. 
Hence we obtain Theorem \ref{thm:2.5} in the case $\Phi \in I(\mathbb{R})$. 

Next, we prove the result in the case of continuous test functions by applying Lemma \ref{lem:3.14}. 
Let $X_n=L_3^\pm(n) \setminus A_\sigma(n)$ and $\omega_n=1/\# X_n$. 
We take $\ell_n: X_n \to \mathbb{R}$ as 
\begin{gather*}
\ell_n(K)
=\log{L}(\sigma, \rho_K). 
\end{gather*}
Then we check that the assumptions of Lemma \ref{lem:3.14} are satisfied. 
Recall that the density function $C_\sigma$ is a non-negative continuous function on $\mathbb{R}$ such that 
\begin{gather*}
\int_{-\infty}^{\infty} C_\sigma(x) \,\frac{dx}{\sqrt{2\pi}}
=1. 
\end{gather*}
The Fourier transform $\widetilde{C}_\sigma$ is also a continuous function and satisfies 
\begin{gather*}
\int_{-\infty}^{\infty} |\widetilde{C}_\sigma(\xi)| \,\frac{d \xi}{\sqrt{2\pi}}
<\infty
\end{gather*}
by Proposition \ref{prop:4.4} $(\mathrm{i})$. 
Therefore both $C_\sigma$ and $\widetilde{C}_\sigma$ belong to the class $L^1(\mathbb{R}) \cap L^\infty(\mathbb{R})$. 
Furthermore, we have 
\begin{gather*}
C_\sigma(x)
=\int_{-\infty}^{\infty} \widetilde{C}_\sigma(\xi) e^{-ix \xi} \,\frac{d \xi}{\sqrt{2\pi}}
\end{gather*}
by inversion formula \eqref{eq:07270221}. 
Hence $C_\sigma$ is a good density function in the sense of Ihara--Matsumoto \cite{IharaMatsumoto2011b}. 
Since we have 
\begin{gather*}
\sum_{K \in X_n} \omega_n \exp\left(i \xi \ell_n(K)\right)
=\frac{1}{N_3^\pm(X)} M_{i \xi, \sigma}^\pm(X)
+O(X^{-\delta}), 
\end{gather*}
Theorem \ref{thm:2.2} yields that \eqref{eq:07270304} is satisfied with $\Phi(x)=e^{i \xi x}$ for any $\xi \in \mathbb{R}$, and that the convergence is uniform in $\xi \in [-R,R]$ for any $R>0$. 
Then we prove that \eqref{eq:07270304} holds with general test functions as follows.  
\begin{itemize}
\item
Let $\sigma>1$ and $\Phi \in C(\mathbb{R})$. 
We define a continuous function $\phi_0$ on $[0,\infty)$ as 
\begin{gather*}
\phi_0(r)
=\max_{|x|\leq r} |\Phi(x)|, 
\end{gather*}
which is non-decreasing and satisfies $\phi_0(r)>0$ for any $r \in [0,\infty)$ and $\phi_0(r) \to\infty$ as $r\to\infty$. 
To obtain \eqref{eq:07270304} in this case, it is sufficient to check that conditions \eqref{eq:07311434} and \eqref{eq:07311435} are satisfied. 
We have $|\log{L}(\sigma,\rho_K)| \leq 2\log\zeta(\sigma)$ for $\sigma>1$. 
Thus, there exists a constant $A_\sigma>0$ depending only on $\sigma$ such that 
\begin{gather*}
\phi_0\left(|\ell_n(K)|\right)\leq\phi_0(A_\sigma)
\end{gather*}
holds for any $K\in{X}_n$, which implies that condition \eqref{eq:07311434} is satisfied. 
Then, we recall that $C_\sigma$ is compactly supported on $\mathbb{R}$ by Theorem \ref{thm:2.1} $(\mathrm{ii})$. 
Hence condition \eqref{eq:07311435} is also satisfied. 
By Lemma \ref{lem:3.14} $(\mathrm{ii})$, we conclude that \eqref{eq:07270304} holds for $\Phi \in C(\mathbb{R})$ since we obtain $|\Phi(x)| \leq \phi_0(|x|)$ by the definition of $\phi_0$. 
\item 
Let $\sigma=1$ and $\Phi \in C^{\exp}(\mathbb{R})$. 
In this case, we define the function $\phi_0$ as 
\begin{gather*}
\phi_0(r)
=C e^{ar}
\end{gather*}
for $r \in [0,\infty)$, where $C$ is a positive constant. 
Then we have 
\begin{gather*}
\sum_{K \in X_n} \omega_n \phi_0\left(|\ell_n(K)|\right)^2
=\frac{C}{N_3^\pm(X)} M_{2a, 1}^\pm(X)
\end{gather*}
since $E_1(X)=A_1(X)=\emptyset$. 
Therefore, Theorem \ref{thm:2.2} implies that condition \eqref{eq:07311434} is satisfied. 
Note that condition \eqref{eq:07311435} is also valid due to Theorem \ref{thm:2.1} $(\mathrm{iii})$. 
Hence Lemma \ref{lem:3.14} $(\mathrm{ii})$ yields that \eqref{eq:07270304} holds for any continuous function $\Phi$ satisfying $\Phi(x) \ll  e^{a|x|}$. 
\item 
Let $\sigma_1<\sigma<1$ and $\Phi \in C_b(\mathbb{R})$.  
By Lemma \ref{lem:3.14} $(\mathrm{i})$, we obtain \eqref{eq:07270304} in this case. 
\item 
Let $\sigma_1<\sigma<1$ and $\Phi \in C^{\exp}(\mathbb{R})$. 
We take $\phi_0(r)=C e^{ar}$ with a constant $C>0$ as before. 
Then we have 
\begin{gather*}
\sum_{\chi \in X_n} \omega_n \phi_0\left(|\ell_n(\chi)|\right)^2
=\frac{C}{N_3^\pm(X)} M_{2a, \sigma}^\pm(X)
\end{gather*}
since $E_\sigma(X)=A_\sigma(X)=\emptyset$ under GRH. 
Hence we deduce from Theorem \ref{thm:2.2} that condition \eqref{eq:07311434} holds. 
As described before, condition \eqref{eq:07311435} is valid. 
As a result, Lemma \ref{lem:3.14} $(\mathrm{ii})$ yields formula \eqref{eq:07270304} in this case. 
\end{itemize}

We finally note that formula \eqref{eq:07270304} implies \eqref{eq:07310306} by the setting of $X_n$, $\omega_n$ and $\ell_n$. 
Hence we complete the proof of Theorem \ref{thm:2.5}.

\renewcommand{\thesection}{\Alph{section}}
\setcounter{section}{1}
\setcounter{equation}{0}

\section*{Appendix: Results on the logarithmic derivative}\label{sec:a}
In the above sections, we proved several results on values $\log{L}(\sigma,\rho_K)$. 
We obtain similar results for $(L'/L)(\sigma,\rho_K)$, which are listed in this section. 
Note that the value $(L'/L)(1,\rho_K)$ is connected with the Euler--Kronecker constant $\gamma_K$ defined by
\begin{gather*}
\gamma_K
=\lim_{s \to1} \left(\frac{\zeta'_K}{\zeta_K}(s)+\frac{1}{s-1}\right). 
\end{gather*}
By definition, $\gamma_\mathbb{Q}$ is equal to Euler's constant $\gamma=0.577\ldots$. 
Thus we obtain 
\begin{gather*}
\frac{L'}{L}(1,\rho_K)
=\gamma_K-\gamma
\end{gather*}
for any non-Galois cubic field $K$ since $L(s,\rho_K)=\zeta_K(s)/\zeta(s)$. 
The proofs for the results in this section are omitted unless we need special remarks arising from the difference between $\log{L}(\sigma,\rho_K)$ and $(L'/L)(\sigma,\rho_K)$. 
Let $\mathcal{X}=(\mathcal{X}_p)_p$ and $\mathcal{Y}=(\mathcal{Y}_p)_p$ be the sequences of independent random elements on $\{A_\mathfrak{a} \mid \mathfrak{a} \in \mathcal{A}\}$ as in Section \ref{sec:2}. 
Then we define the random Euler products $L(s,\mathcal{X})$ and $L(s,\mathcal{Y})$ as in \eqref{eq:07271414}. 

\begin{theorem}\label{thm:A.1}
For $\sigma>1/2$, there exists a non-negative $C^\infty$-function $\mathcal{C}_\sigma$ such that 
\begin{gather*}
\mathbb{P}\left( \frac{L'}{L}(\sigma,\mathcal{X}) \in A \right)
=\int_{A} \mathcal{C}_\sigma(x) \,\frac{dx}{\sqrt{2\pi}}
\end{gather*}
holds for all $A \in \mathcal{B}(\mathbb{R})$. 
Furthermore, it satisfies the following properties.  
\begin{itemize}
\item[$(\mathrm{i})$]
If $1/2<\sigma\leq1$, we have $\supp \mathcal{C}_\sigma=\mathbb{R}$, that is, $\mathcal{C}_\sigma(x)$ is not identically zero in any interval on $\mathbb{R}$. 
\item[$(\mathrm{ii})$]
If $\sigma>1$, the function $\mathcal{C}_\sigma$ is compactly supported. 
\item[$(\mathrm{iii})$]
Let $\sigma>1/2$. 
Then the integral 
\begin{gather*}
\int_{-\infty}^{\infty} e^{ax} \mathcal{C}_\sigma(x) \,\frac{dx}{\sqrt{2\pi}}
\end{gather*}
is finite for any $a>0$. 
\end{itemize}
Similarly, for $\sigma>2/3$, there exists a non-negative $K^\infty$-function $K_\sigma$ such that 
\begin{gather*}
\mathbb{P}\left( \frac{L'}{L}(\sigma,\mathcal{Y}) \in A \right)
=\int_{A} \mathcal{K}_\sigma(x) \,\frac{dx}{\sqrt{2\pi}}
\end{gather*}
holds for all $A \in \mathcal{B}(\mathbb{R})$. 
Furthermore, it satisfies the following properties. 
\begin{itemize}
\item[$(\mathrm{i'})$]
If $2/3<\sigma\leq1$, we have $\supp \mathcal{K}_\sigma=\mathbb{R}$, that is, $\mathcal{K}_\sigma(x)$ is not identically zero in any interval on $\mathbb{R}$. 
\item[$(\mathrm{ii'})$]
If $\sigma>1$, the function $\mathcal{K}_\sigma$ is compactly supported. 
\item[$(\mathrm{iii'})$]
Let $\sigma>2/3$. 
Then the integral 
\begin{gather*}
\int_{-\infty}^{\infty} e^{ax} \mathcal{K}_\sigma(x) \,\frac{dx}{\sqrt{2\pi}}
\end{gather*}
is finite for any $a>0$. 
\end{itemize}
\end{theorem}

Denote by $\mathcal{E}_\sigma(X)$ the subset of $L_3^\pm(X)$ such that $\mathcal{E}_\sigma(X)=E(X)$ for $\sigma_1<\sigma \leq1$ and $\mathcal{E}_\sigma(X)=\emptyset$ for $\sigma>1$, where $\sigma_1$ is a real number for which \eqref{eq:07251505} holds, and $E(X)$ is defined by \eqref{eq:07252237}. 
Then we define $\mathcal{M}_{z, \sigma}^{\pm}(X)$ as 
\begin{gather*}
\mathcal{M}_{z, \sigma}^{\pm}(X)
=\sum_{K \in L_3^\pm(X) \setminus \mathcal{E}_\sigma(X)} \exp\left(z \frac{L'}{L}(\sigma,\rho_K)\right) 
\end{gather*}
as an analogue of the $z$-th moment $M_{z, \sigma}^{\pm}(X)$ of \eqref{eq:07250226}. 

\begin{theorem}\label{thm:A.2}
Let $\sigma_1$ be a real number for which \eqref{eq:07251505} holds. 
Then there exists an absolute constant $\delta>0$ such that 
\begin{gather}\label{eq:07311539}
\mathcal{M}_{z, \sigma}^{\pm}(X)
=C^\pm \frac{1}{12\zeta(3)} X 
\int_{-\infty}^{\infty} e^{zx} \mathcal{C}_\sigma(x) \,\frac{dx}{\sqrt{2\pi}} 
+O\left(X\exp\left(-\delta \frac{\log{X}}{\log\log{X}}\right)\right)
\end{gather}
holds for $\sigma>\sigma_1$ with $z \in \mathbb{C}$ satisfying $|z| \leq b_\sigma R_{\sigma}(X)$, where $b_\sigma$ is a positive constant, and $R_{\sigma}(X)$ is defined as in \eqref{eq:07281607}. 
The implied constant in \eqref{eq:07311539} depends only on $\sigma$. 
\end{theorem}

Note that the subset $\mathcal{E}_\sigma(X)$ differs from $E_\sigma(X)$ at $\sigma=1$. 
We hereby explain the reason why this difference arises. 
In the same line as \eqref{eq:07301444}, one can show the asymptotic formula 
\begin{align*}
&\sum_{K \in L_3^\pm(X) \setminus E(X)} \exp\left(z \frac{L'}{L}(\sigma,\rho_K)\right) \\
&=C^\pm \frac{1}{12\zeta(3)} X 
\int_{-\infty}^{\infty} e^{zx} \mathcal{C}_\sigma(x) \,\frac{dx}{\sqrt{2\pi}} 
+O\left(X\exp\left(-\delta \frac{\log{X}}{\log\log{X}}\right)\right)
\end{align*}
for $\sigma>\sigma_1$ and $|z| \leq b_\sigma R_\sigma(X)$, where $\delta$ is an absolute positive constant. 
To prove formula \eqref{eq:07311539} at $\sigma=1$ with $\mathcal{E}_1(X)=\emptyset$, one needs to obtain the upper bound
\begin{gather}\label{eq:07311603}
\sum_{K \in E(X)} \exp\left(z \frac{L'}{L}(\sigma,\rho_K)\right)
\ll X\exp\left(-\delta \frac{\log{X}}{\log\log{X}}\right)
\end{gather}
for $|z| \leq b_\sigma R_\sigma(X)$ with $\sigma=1$. 
This is true for $\sigma>1$ due to $(L'/L)(\sigma,\rho_K) \ll 1$ and $\# E(X) \ll X^{1-\delta}$. 
However, we know just 
\begin{gather*}
\frac{L'}{L}(1,\rho_K) \ll \log{X}
\end{gather*}
for $K \in L_3^\pm(X)$ without any assumptions, which prevents us from obtaining \eqref{eq:07311603} at $\sigma=1$. 
A related difficulty on studying the distribution of values $(L'/L)(1,\chi_d)$ was discussed in the Ph.D.\ thesis of Mourtada \cite{Mourtada2013}. 
See also Lamzouri \cite{Lamzouri2015b} for more information about the Euler--Kronecker constants of quadratic fields. 

\begin{theorem}\label{thm:A.3}
Assume GRH and upper bound \eqref{eq:07251548} for each $\epsilon>0$ with some constants $\alpha$ and $\beta$ such that $0<\alpha<5/6$. 
Let $\sigma_2$ be as in \eqref{eq:07302054}. 
Then there exists a constant $\delta=\delta(\sigma)>0$ such that 
\begin{align}\label{eq:07311600}
\mathcal{M}_{z, \sigma}^{\pm}(X)
&=C^\pm \frac{1}{12\zeta(3)} X 
\int_{-\infty}^{\infty} e^{zx} \mathcal{C}_\sigma(x) \,\frac{dx}{\sqrt{2\pi}} \\
&\quad+K^\pm \frac{4\zeta(1/3)}{5\Gamma(2/3)^3\zeta(5/3)} X^{5/6}
\int_{-\infty}^{\infty} e^{zx} \mathcal{K}_\sigma(x) \,\frac{dx}{\sqrt{2\pi}} 
+O\left(X^{5/6-\delta}\right) \nonumber
\end{align}
holds for any $\sigma>\max(\sigma_2,2/3)$ with $z \in \mathbb{C}$ satisfying $|z| \leq \tilde{b}_\sigma \tilde{R}_{\sigma}(X)$, where $\tilde{b}_\sigma$ is a positive constant, and $\tilde{R}_\sigma(X)$ is defined as in \eqref{eq:07301552}. 
The implied constant in \eqref{eq:07311600} depends only on $\sigma$.
\end{theorem}

\begin{corollary}\label{cor:A.1}
There exists an absolute constant $\delta>0$ such that  
\begin{align*}
\sum_{K \in L_3^\pm(X) \setminus E(X)} \exp(z \cdot \gamma_K)
&=C^\pm \frac{1}{12\zeta(3)} X e^{\gamma z} 
\int_{-\infty}^{\infty} e^{rx} \mathcal{C}_1(x) \,\frac{dx}{\sqrt{2\pi}} \\
&\quad+O\left(X \exp\left(-\delta \frac{\log{X}}{\log\log{X}}\right)\right)
\end{align*}
holds for any $z \in \mathbb{C}$ with $|z| \leq b_1 R_1(X)$, where the implied constant is absolute. 
\end{corollary}

\begin{corollary}\label{cor:A.2}
Assume GRH and upper bound \eqref{eq:07251548} with some constants $\alpha$ and $\beta$ such that $3\alpha+\beta<5/2$. 
Then there exists an absolute constant $\delta>0$ such that  
\begin{align*}
\sum_{K \in L_3^\pm(X)} \exp(z \cdot \gamma_K)
&=C^\pm \frac{1}{12\zeta(3)} X e^{\gamma z} 
\int_{-\infty}^{\infty} e^{zx} \mathcal{C}_1(x) \,\frac{dx}{\sqrt{2\pi}} \\
&\quad+K^\pm \frac{4\zeta(1/3)}{5\Gamma(2/3)^3\zeta(5/3)} X^{5/6} e^{\gamma z} 
\int_{-\infty}^{\infty} e^{zx} \mathcal{K}_1(x) \,\frac{dx}{\sqrt{2\pi}} \\
&\quad+O\left(X \exp\left(-\delta \frac{\log{X}}{\log\log{X}}\right)\right)
\end{align*}
holds for any $z \in \mathbb{C}$ with $|z| \leq \tilde{b}_1 \tilde{R}_1(X)$, where the implied constant is absolute. 
\end{corollary}

Define the quantity $\mathcal{D}_\sigma^\pm(X; a)$ as 
\begin{gather*}
\mathcal{D}_\sigma^\pm(X; a)
=\frac{\# \left\{ K \in L_3^\pm(X) ~\middle|~ \frac{L'}{L}(\sigma,\rho_K) \leq a \right\}}{N_3^\pm(X)}
-\int_{-\infty}^{a} \mathcal{C}_\sigma(x) \,\frac{dx}{\sqrt{2\pi}}
\end{gather*}
for $\sigma>1/2$ and $a \in \mathbb{R}$. 

\begin{theorem}\label{thm:A.4}
Let $\sigma_1$ be a real number for which \eqref{eq:07251505} holds. 
Then we obtain
\begin{gather*}
\sup_{a \in \mathbb{R}} \left|\mathcal{D}_\sigma^\pm(X; a)\right|
\ll \frac{1}{R_\sigma(X)}, 
\end{gather*}
where $R_\sigma(X)$ is as in \eqref{eq:07281607}. 
\end{theorem}

Finally, we present a result similar to Theorem \ref{thm:2.5}. 
There exists a difference from the result for $\log{L}(\sigma,\rho_K)$ when we consider the case $\sigma=1$ by the same reason as in Theorem \ref{thm:A.3}. 
For this, we define a subclass of $C(S)$ as
\begin{gather*}
C^{\poly}(\mathbb{R})
=\left\{ \Phi \in{C}(\mathbb{R}) ~\middle|~ \text{$\Phi(x) \ll |x|^a$ with some $a>0$} \right\}. 
\end{gather*}

\begin{theorem}\label{thm:A.5}
Let $\sigma_1$ be a real number for which \eqref{eq:07251505} holds. 
Then the limit formula
\begin{gather*}
\lim_{X \to\infty} \frac{1}{N_3^\pm(X)} 
\sideset{}{'} \sum_{K \in L_3^\pm(X)} \Phi\left(\frac{L'}{L}(\sigma,\rho_K)\right)
=\int_{-\infty}^{\infty} \Phi(u) \mathcal{C}_\sigma(x) \,\frac{dx}{\sqrt{2\pi}}
\end{gather*}
holds in the following cases: 
\begin{itemize}
\item 
$\sigma>1$ and $\Phi \in C(\mathbb{R}) \cup I(\mathbb{R})$; 
\item 
$\sigma=1$ and $\Phi \in C^{\poly}(\mathbb{R}) \cup I(\mathbb{R})$; 
\item 
$\sigma_1<\sigma<1$ and $\Phi \in C_b(\mathbb{R}) \cup I(\mathbb{R})$ without assuming GRH; 
\item 
$\sigma_1<\sigma \leq1$ and $\Phi \in C^{\exp}(\mathbb{R}) \cup I(\mathbb{R})$ if we assume GRH.  
\end{itemize}
\end{theorem}

%\bibliographystyle{amsplain}
%\bibliography{refs}

\providecommand{\bysame}{\leavevmode\hbox to3em{\hrulefill}\thinspace}
\providecommand{\MR}{\relax\ifhmode\unskip\space\fi MR }
\providecommand{\MRhref}[2]{%
  \href{http://www.ams.org/mathscinet-getitem?mr=#1}{#2}
}
\providecommand{\href}[2]{#2}

\end{document}